
\magnification=\magstep1
\input amstex
\documentstyle{amsppt}
\leftheadtext{J. Jer\'onimo-Castro, E. Makai, Jr.}
\rightheadtext{Ball characterizations}
\topmatter
\title 
\centerline{Ball characterizations} 
\centerline{in spaces of constant curvature}
\endtitle
\author J. Jer\'onimo-Castro*, E. Makai, Jr.** \endauthor
\address 
$^*$ 
Facultad de Ingenier\'\i a, Universidad Aut\'onoma de Quer\'etaro, Centro
Uni\-ver\-si\-ta\-rio, 
\newline Cerro de las Campanas s/n C.P. 76010,
Santiago de Quer\'etaro, Qro., ME\-XI\-CO
\newline
$^{**}$ Alfr\'ed R\'enyi Mathematical Institute, Hungarian Academy of Sciences,
\newline
H-1364 Budapest, Pf. 127, HUNGARY
\newline
{\rm{http://www.renyi.hu/\~{}makai}}
\endaddress
\email 
$^*$ jeronimo\@cimat.mx, jesusjero\@hotmail.com
\newline
$^{**}$ makai.endre\@renyi.mta.hu\endemail
\thanks *Research (partially) supported by CONACYT, SNI 38848
\newline
**Research (partially) supported by Hungarian National Foundation for 
Scientific Research, grant nos. T046846, T043520, K68398, K81146\endthanks
\subjclass {\it Mathematics Subject Classification 2010.} 
Primary: 52A55; Secondary 52A20
\endsubjclass
\abstract 
High proved the following theorem. 
If the intersections of any two congruent
copies of a plane convex body are centrally symmetric, then this body is a
circle. 
In our paper we extend the theorem of High to spherical, Euclidean 
and hyperbolic spaces, under some regularity assumptions.
Suppose that in any of these spaces
there is a pair of closed convex sets of class $C^2_+$ with interior points,
different from the whole space,
and the intersections of any congruent copies of these sets are centrally 
symmetric (provided they have non-empty interiors).
Then our sets are congruent balls. Under the same hypotheses, but if we require
only central symmetry of small intersections, then our sets are either congruent
balls, or paraballs, or have as connected components of their boundaries
congruent hyperspheres (and the converse implication also holds). 

Under the same hypotheses, if we require
central symmetry of all compact intersections, then either our
sets are congruent balls or paraballs, or have as connected
components of their boundaries congruent hyperspheres, and either $d \ge 3$,
or $d=2$ and
one of the sets is bounded by one hypercycle, or both sets are congruent 
parallel domains of straight lines, or there are no more compact intersections
than those bounded by two finite hypercycle arcs
(and the converse implication also holds).

We also prove a dual theorem. If in any of these spaces
there is a pair of smooth closed convex sets, such that both of them
have supporting spheres at any of their boundary
points --- for $S^d$ of radius less than $ \pi /2$ --- and the closed convex
hulls of any congruent copies of these sets are centrally symmetric,
then our sets are congruent balls.
\endabstract
\endtopmatter\document

\head 1. Introduction \endhead

We will investigate closed convex sets with non-empty interior 
in $S^d$ ($d$-di\-men\-si\-on\-al sphere), ${\Bbb R}^d$, $H^d$ ($d$-dimensional
hyperbolic space).


R. High proved the following theorem.


\proclaim{Theorem} ([H])
Let $K \subset {\Bbb R}^2$ be a convex body. Then
the following statements are equivalent:
\roster
\item
All intersections $(\varphi K) \cap (\psi K)$, 
having interior points,
where $\varphi, \psi :
{\Bbb R}^2 \to {\Bbb R}^2$ are congruences, are centrally
symmetric.
\item
$K$ is a circle. $ \blacksquare $
\endroster
\endproclaim

It seems that his proof gives the analogous statement, when $\varphi , \psi $
are only allowed to be orientation preserving congruences.


\definition{Problem} 
Describe the pairs of closed convex sets with interior points, in $S^d$,
${\Bbb R}^d$ and $H^d$, different from the whole space,
whose any congruent copies have a centrally symmetric
intersection, provided this intersection has interior points, or have a
centrally symmetric closed convex hull of their unions.
Evidently, two
congruent balls (for $S^d$ of radii at most $\pi /2$), or two parallel slabs
in ${\Bbb R}^d$,
have a centrally symmetric intersection, provided this intersection 
has a non-empty interior, and have a centrally symmetric closed convex hull of 
their unions.
\enddefinition


The authors are indebted to L. Montejano (Mexico City) and G. Weiss (Dresden)
for having turned their interest to characterizations of pairs of convex
bodies with all translated (for ${\Bbb{R}}^d$)
or congruent copies having a centrally or axially
symmetric intersection or convex hull of the union, respectively, 
or with other symmetry properties, e.g., having some affine symmetry.


Central symmetry of a set $X \subset S^d$ with respect to a point $O \in S^d$
is equivalent to central symmetry of $X$ with respect to the point
$-O$ antipodal to $O$. However, the two transformations: central symmetry with
respect to $O$, and central symmetry with respect to $-O$, coincide. In all our
theorems, for the case of $S^d$, we will investigate sets $X \subset S^d$
contained in an open hemisphere, say the southern one. Such a set cannot have
a center of symmetry on the equator, but it may have one in the open southern
or in the open northern hemisphere, and then it has two antipodal centres of
symmetry, one in the open southern, and one in the open northern hemisphere.
In such case we will use the one in the southern hemisphere.


The aim of our paper will be to give partial answers to these
problems.
To exclude trivialities, we always
suppose that {\it{our sets are different from the whole space}},
and also we investigate only such cases, when the 
{\it{intersection has interior points}}.
For $S^d$, ${\Bbb R}^d$ and $H^d$, where $d  \ge 2$,
we prove the analogue of the above
theorem under some regularity assumptions
($C^2$ for $S^d$, $C^2$ and having an extreme point for ${\Bbb R}^d$, and
$C^2_+$ for $H^d$, respectively). 

For $S^d$, ${\Bbb R}^d$ and $H^d$, under the above mentioned 
regularity assumptions, we have the following. If {\it{all sufficiently small 
intersections}} of congruent copies of two closed convex sets $K$ and $L$ with
interior points, having a non-empty interior, are centrally symmetric, 
then all connected components of the boundaries of the
two sets are
congruent spheres, paraspheres or hyperspheres.
({\it{``Sufficiently small''}} means here: of sufficiently small diameter.)
Under the same regularity assumptions, if {\it{all intersections}} of
congruent copies of two closed convex sets with
interior points, 

\newpage

having a non-empty interior, are centrally symmetric, 
then they are congruent balls.
There is a question ``between'' the above two questions. Suppose the same
regularity assumptions, and also that {\it{all compact intersections}} are
centrally symmetric. Then there are several possibilities for $K$ and $L$, and
there is a complete description also for this case.

The dual question is the question of centrally symmetric closed convex hull 
of any congruent copies of $K$ and $L$. Under the
hypotheses that both $K$ and $L$ are smooth and at any of their boundary
points have supporting spheres, for $S^d$ of radii less than $ \pi /2 $, 
the only case is two congruent balls, for $S^d$ of radii less than $\pi /2$. 
Observe that for $S^d$, ${\Bbb{R}}^d$ and $H^d$ the hypotheses imply that
any existing sectional curvature of $K$ and $L$ is positive, positive, or
greater than $1$, in the three cases, respectively.

Surveys about characterizations of central symmetry, for convex bodies in
${\Bbb R}^d$, cf. in
[BF], \S 14, pp. 124-127, and, more recently, in [HM], \S 4.

In later papers 
we will give sharper theorems on the one hand about ${\Bbb
R}^d$  (for $d \ge 2$), and on the other hand about $S^2$ and $H^2$. 

In ${\Bbb{R}}^d$ 
we will describe all pairs of closed convex sets with interior points,
different from ${\Bbb{R}}^d$, without any additional hypotheses,
whose any congruent copies have a centrally symmetric intersection (provided
this intersection has interior points). For
$d \ge 2$ these are: (1) two
congruent balls, or (2) two (incongruent)
parallel slabs. (Observe that in Theorem 2 of
this paper the hypothesis about the existence of
an extreme point of $K$ or $L$ excludes the case of two parallel slabs.) 

For the dual question, we will describe in ${\Bbb{R}}^d$ 
all pairs of closed convex sets with interior points, different from
${\Bbb{R}}^d$, without any additional
hypotheses, whose any congruent copies have a union with a
centrally symmetric closed convex hull. For $d \ge 2$ these are: 
(1) $K$ and $L$ are infinite cylinders
over balls of dimensions $2 \le i,j \le d$, having equal radii (this includes
the case of two congruent balls).
(2) one of $K$ and $L$ is an infinite cylinder with dimension of axis
$0 \le i \le d-1$ and with base compact, and the other one is a slab
(this includes the case of two slabs).
The
methods applied for ${\Bbb{R}}^d$ are completely different from those in this
paper. They use some theorems of V. Soltan in \cite{So05}, \cite{So06},
and some other considerations, even for the case of intersections. 

Further, in $S^2$, ${\Bbb{R}}^2$ and $H^2$
we will even describe the pairs of closed convex sets with interior points,
different from the entire space,
whose any congruent copies have a (1) centrally, or (2) axially
symmetric intersection (provided
this intersection 

\newpage

has interior points), under the
hypothesis that (1) for the case of $H^2$, if all connected components of the
boundaries both
of $K$ and $L$ are straight lines, then their numbers are finite, or 
(2) for the case of $H^2$, 
the numbers of the connected components both of $K$ and $L$ are finite. 

Suppose (1). Then $K$ and $L$ are congruent circles, or, for ${\Bbb{R}}^2$,
two (incongruent) parallel strips. (The case of ${\Bbb{R}}^2$ here
also follows from the case of ${\Bbb{R}}^d$ above.)

Suppose (2). Then, 
{\it{for}} $S^2$, 
$K$ and $L$ are (incongruent) circles. {\it{For}} ${\Bbb{R}}^2$ there are 
five cases, each satisfying that any of $K$ and $L$ is a circle, 
a parallel strip or a half-plane. {\it{For}} 
$H^2$ there are a large number of such
cases, each satisfying that all connected components of the boundaries both
of $K$ and $L$ are cycles or straight lines, their curvatures 
depending on the component (this being true
already if all intersections of a sufficiently small diameter are centrally
symmetric, also for $S^2$ and ${\Bbb{R}}^2$ ---
in some analogy with our Theorem 1). Furthermore, if none of
$K$ and $L$ is a circle, then the boundaries of both of them 
have at most two hypercycle or straight line connected 
components, and moreover, if either for $K$ or for $L$
there are two such components, then the respective
set is a parallel domain of a straight line. Some cases are: 
(a) two (incongruent) circles, (b) two paracycles, (c)
two congruent closed convex sets, each 
bounded by one hypercycle, (d) two half-planes, (e) two
congruent parallel domains of lines.
The methods applied for $S^2$, ${\Bbb{R}}^2$ and $H^2$ are
refined versions of the methods applied in this paper.

Still we remark that for the dual problem we cannot give better results for
$d=2$ than Theorem 4 in this paper for $d \ge 2$.

\head 2. New results \endhead

Let $d \ge 2$ be an integer. We investigate the spaces of constant curvature 
$S^d$, ${\Bbb{R}}^d$ and $H^d$. Actually our proofs use absolute geometry, 
i.e., are independent of the parallel axiom. In particular, the case of
${\Bbb{R}}^d$ in our Theorem 1 is not simpler than the general case. Theorem
2 follows from Theorem 1. There the case of $H^d$ requires some
additional considerations.  

As usual, we write ${\text{conv}}\,( \cdot )$,  ${\text{aff}}\,( \cdot )$,
${\text{diam}}\,(\cdot )$, ${\text{cl}}( \cdot )$, 
${\text{int}}\,( \cdot )$, ${\text{bd}}( \cdot )$,
and ${\text{rel\,bd}}\,( \cdot )$ for the convex hull, affine hull, diameter,
closure, interior, boundary and relative boundary (provided it is
understood in which subspace do we consider it) of a set. Further,
${\text{dist}}( \cdot , \cdot )$ denotes distance.

As general hypotheses in all our statements we use
$$
\cases
X {\text{ will be }}S^d, \,\,{\Bbb R}^d {\text{ or }}H^d, {\text{ for }} 
d \ge 2, {\text{ and }} K,L \subsetneqq X {\text{ will be closed}}\\
{\text{convex sets with interior points, and }} \varphi , \psi : X \to X, 
{\text{ sometimes}} \\
{\text{with indices, will be orientation preserving congruences.}}
\endcases
\tag *
$$

\newpage

Further, we will need the following weakening of the $C^2$ property.
$$
\cases
{\text{Let for each }} x \in {\text{bd}}\,K, {\text{ and each }} 
y \in {\text{bd}}\,L,{\text{ there exist an }} 
\varepsilon _1 (x) > 0, \\
{\text{and an }} \varepsilon _1 (y)>0, {\text{ such that }}
K {\text{ and }}L {\text{ contain balls of radius }} \varepsilon _1 (x) \\
 {\text{and }} \varepsilon _1 (y), {\text{ containing }}
x {\text{ and }} y {\text{ in their boundaries, respectively.}} 
\endcases
\tag A
$$
Moreover, we will need the following property, which together with \thetag{A}
is a weakening of the $C^2_+$ property.
$$
\cases
{\text{Let for each }} x \in {\text{bd}}\,K, {\text{ and each }} 
y \in {\text{bd}}\,L,{\text{ there exist an}} \\
\varepsilon _2 (x) > 0 {\text{ and }} \varepsilon _2 (y) > 0, 
{\text{ such that the set of points of }} K \\
{\text{and }} L,
{\text{ lying at a distance at most }} 
\varepsilon _2 (x) {\text{ and }} \varepsilon _2 (y) {\text{ from }} x \\
{\text{and from }} y, 
{\text{ is contained in a
ball }} B
{\text{ (for }} X=S^d, \,\,{\Bbb R}^d) {\text{ or}} \\
{\text{in a convex set }} B {\text{ bounded by a
hypersphere (for }} X = H^d), \\
{\text{with bd\,}}B {\text{ having sectional curvatures at least }} 
\varepsilon _2 (x) {\text{ and}} \\
\varepsilon _2 (y), 
{\text{ and with bd\,}}B {\text{ containing }} x {\text{ or }} y, 
{\text{ respectively.}} 
\endcases
\tag B
$$
Clearly (A) implies
smoothness and (B) implies strict convexity, respectively.
Observe that both in (A) and (B) $\varepsilon _i( x ) > 0$ and 
$\varepsilon _i(y) >
0$ can be decreased, 
and then (A) and (B) remain valid.


The following Theorem 1 will be the basis of our considerations for the case
of intersections.
Observe that in Theorem 1, (2), for ${\Bbb R}^d$ and $H^d$,
hyperplanes cannot occur, by the
hypothesis about the existence of
an extreme point of $K$ or $L$, and by $C^2_+$ (or by (B)), respectively.
By the same reason, in Theorem 2, for ${\Bbb R}^d$ parallel strips cannot
occur.


\proclaim{Theorem 1}
Let $X$ be $S^d$, ${\Bbb{R}}^d$ or $H^d$, and let $K,L$ and $\varphi, \psi $
be as in (*).
Let us assume $C^2$ for $K$ and $L$ (actually $C^2$ can be weakened to
\thetag{A}). 
For $X={\Bbb R}^d$ assume additionally 
that one of $K$ and $L$ has an extreme point. For
$X=H^d$ assume $C^2_+$ for $K$ and $L$ (actually $C^2_+$ can be weakened
to \thetag{A} and \thetag{B}).
Then the following statements are equivalent.
\roster
\item
There exists some $\varepsilon = \varepsilon (K,L) > 0$, such that
for each $\varphi , \psi $, for 
which
${\text{\rm{int}}}$
\newline
$\left( (\varphi K) \cap (\psi L) \right) \ne \emptyset $ 
and ${\text{\rm{diam}}} \left( (\varphi K) \cap (\psi L) \right)
< \varepsilon $, 
we have that $(\varphi K) \cap (\psi L)$
is centrally symmetric.

\newpage

\item
The connected components of the boundaries of both $K$ and $L$ are congruent 
spheres (for $X=S^d$ of radius at most $\pi /2$),
or paraspheres, or congruent hyperspheres (for ${\Bbb{R}}^d$ and $H^d$
degeneration to hyperplanes being not admitted). 
For the case of congruent spheres or
paraspheres we have that either $K$ and $L$ are congruent balls 
(for $X=S^d$ of radius at most $\pi /2$), or they are paraballs. 
\endroster
\endproclaim


\proclaim{Theorem 2}
Let $X$ be $S^d$, ${\Bbb{R}}^d$ or $H^d$, and let $K,L$ and $\varphi , \psi $
be as in (*).
Let us assume $C^2$ for $K$ and $L$ (actually $C^2$ can be weakened to
\thetag{A}). For $X={\Bbb R}^d$ assume additionally 
that one of $K$ and $L$ has an extreme point.  For
$X=H^d$ assume $C^2_+$ for $K$ and $L$ 
(actually $C^2_+$ can be weakened to \thetag{A} and \thetag{B}).
Then the following statements are equivalent.
\roster
\item
For each $\varphi , \psi $, for 
which\,\,\,{\rm{int}}$\left( (\varphi K) \cap (\psi L) \right) \ne
\emptyset $ (here we may suppose additionally that 
$(\varphi K) \cap (\psi L)$ has at most one infinite point),
we have that $(\varphi K) \cap (\psi L)$ is centrally symmetric.
\item
$K$ and $L$ are two congruent balls, and, for $X=S^d$, 
their common radius is at most $\pi /2$.
\endroster
\endproclaim


Observe that in Theorem 1, (1)  we considered only small intersections
with non-empty interiors,
in Theorem 2, (1) all intersections with non-empty interiors (or in brackets,
additionally having at most one infinite point). There
is a third possibility, a condition ``between'' these two conditions: namely
all compact intersections. This will be done in the following theorem.


\proclaim{Theorem 3}
Let $X$ be $S^d$, ${\Bbb{R}}^d$ or $H^d$, and let $K,L$ and $\varphi , \psi $
be as in (*).
Let us assume $C^2$ for $K$ and $L$ (actually $C^2$ can be weakened to
\thetag{A}). For $X={\Bbb R}^d$ assume additionally 
that one of $K$ and $L$ has an extreme point.  For
$X=H^d$ assume $C^2_+$ for $K$ and $L$ 
(actually $C^2_+$ can be weakened to \thetag{A} and \thetag{B}).
Then the following statements are equivalent.
\roster
\item
For each $\varphi , \psi $, for 
which\,\,\,{\rm{int}}$\left( (\varphi K) \cap (\psi L) \right) \ne
\emptyset $ and $(\varphi K) \cap (\psi L)$ is compact,
we have that $(\varphi K) \cap (\psi L)$ is centrally symmetric.
\item
$K$ and $L$ are either 

\noindent (a) two congruent balls, and, for $X=S^d$, 
their common radius is at most $\pi /2$, or

\noindent (b) two paraballs, or 

\noindent (c) the connected components of the boundaries of both $K$ and $L$ 
are congruent 
hyperspheres (degeneration to hyperplanes being not admitted), and either

($\alpha $)
$d \ge 3$, or

($\beta $)
$d=2$, and either

\newpage

\hskip.3cm ($ \beta '$)
one of $K$ and $L$ is bounded by one hypercycle, or 

\hskip.3cm ($\beta ''$)
$K$ and $L$ are congruent parallel domains of straight lines, or

\hskip.3cm ($\beta '''$)
there are no more compact intersections $(\varphi K) \cap (\psi L)$
than those bounded by two finite hypercycle arcs.
\endroster
\endproclaim


We observe that given $K$ and $L$ we cannot in general 
decide whether $(\beta ''')$
holds for them or not. So this is not such an explicit description as the
other cases in Theorem 3, (2). 

As an example, suppose that both $K$ and $L$ have two connected components
of their boundaries, $K_1,K_2$, and $L_1,L_2$, say. Let $K_1,K_2$, and
$L_1,L_2$ have no common infinite points (they have no common finite points).
Let the first and last
points of $K_1$, in the positive sense, be $k_{11}$ and $k_{12}$, and those of
$K_2$ be $k_{21}$ and $k_{22}$. Then the straight lines $k_{11}k_{21}$ and
$k_{12}k_{22}$ intersect each other at a point $O_K \in H^2$, 
and these lines make two
opposite angles $\alpha _K \in (0, \pi )$, with
their respective angular domains containing $K_1$ and $K_2$. (Then 
$K_1 \cup K_2$ is centrally symmetric with respect to $O_K$.) In an
analogous way we define the angle $\alpha _L$. 
Then we claim that 
$$
\alpha _K + \alpha _L  > \pi \Longrightarrow \lnot (\beta ''')
\tag C
$$

In fact,
we may choose $\varphi O_K = \psi O_L = 0$. Then $\varphi $ and $\psi $ are 
determined up to some rotations,
which we can choose so that the images by $\varphi $ and by $\psi $ of the 
above described, altogether four, {\it{open}}  
angular domains of angles $\alpha _K$ and $\alpha _L$ cover $S^1$. Then
$(\varphi K) \cap (\psi L)$ is compact and is not bounded by two finite
hypercycle arcs. Hence
$(\beta ''')$ does not hold, {\it{and \thetag{C} is shown}}. 
Maybe in \thetag{C} we have actually an equivalence?

\vskip.2cm


Observe that the hypotheses of the following Theorem 4 imply compactness of $K$
and $L$. 
Moreover, for $S^d$, ${\Bbb{R}}^d$, or $H^d$ they imply
that any existing sectional curvature both of $K$ and of $L$ is greater than
$0$, $0$, or $1$, respectively,
which for $H^d$ is a serious geometric restriction.

The convex hull of a set $Y \subset H^d$ 
is defined as for ${\Bbb{R}}^d$ (or one can
use the collinear model). For $Y \subset S^d$, 
since we will use only sets $Y$ with interior points, we will call $Y$
{\it{convex}}, if for any two non-antipodal points of $Y$ the unique
smaller great circle arc connecting them belongs to $Y$. Then for any two
antipodal points $\pm x \in Y$ there is a point $y \in Y$ such that $y \ne \pm
x$, and then the smaller large circle arcs ${\widehat{xy}}$ and 
${\widehat{(-x)y}}$ lie in $Y$. So also in the antipodal case there is at
least one half large 

\newpage

circle arc connecting $x$ and $-x$ in $Y$. 
The convex hull, or closed convex hull of a set $Y \subset S^d$ is defined
using this definition of convexity in $S^d$. For $X$ being $S^d$,
${\Bbb{R}}^d$ or $H^d$, and $Y \subset X$, we write ${\text{conv}}\,Y$ and
${\text{cl\,conv}}\,Y$ for the convex hull, and for the closed convex hull of
$Y$, respectively.

We say that 
{\it{a set $Y$ in $S^d$, ${\Bbb{R}}^d$ or $H^d$ has at its boundary point $x$ 
a supporting sphere if
there exists a ball containing $Y$, for $S^d$ of radius at most $ \pi /2$, 
such that $x$ belongs to the boundary of
this ball}}, which boundary is called the {\it{supporting sphere}}.


\proclaim{Theorem 4} 
Let $X$ be $S^d$, ${\Bbb{R}}^d$ or $H^d$, and let $K$,
$L$ and $\varphi , \psi $
be as in (*). Let $K$ and $L$ be smooth, and
let both $K$ and $L$ have supporting spheres at any of their
boundary points, for $S^d$ of radius less than $\pi /2$.
Then the following two statements are equivalent:
\roster
\item
For each $\varphi , \psi $, we have that 
${\text{\rm{cl\,conv}}}\left( (\varphi K) \cup (\psi L) \right) $
(where for $S^d$
we may additionally suppose that 
${\text{\rm{diam}}}\,
[{\text{\rm{cl\,conv}}}\left( (\varphi K) \cup (\psi L) \right) ]$ is smaller
than $\pi $, but is arbitrarily close to $\pi $, and for ${\Bbb {R}}^d$ 
and $H^d$ that this diameter is arbitrarily large), is centrally symmetric.
\item
$K$ and $L$ are two congruent balls (for the case of $S^d$ of radius less than
$ \pi /2$).
\endroster
\endproclaim


Observe that in the case of intersections, we had three different equivalent
statements for small, for compact, 
and for all intersections (namely Theorem 1, (2), 
Theorem 3, (2) and Theorem 2 (2)), 
while for the case of closed convex hull of
the union, large convex hulls, or all convex hulls give the same result. 


\definition{Remark} Possibly Theorems 1 and 2 hold for $S^d$ without any
regularity hypotheses, and for
$H^d$ 
only assuming strict convexity (a weakening of (B)). Without supposing
strict convexity Theorem 1 does not hold even for $K, L \subset H^d$ having
analytic boundaries.
Namely, let $1 \le d_1,d_2$ be integers with $d_1 + d_2 < d$. Let $K_0 \subset
H^{d_1}$ and $L_0 \subset H^{d_2}$ be any closed convex sets with nonempty
interiors; their boundaries may be supposed to be analytic. Let $\pi _i:
H^d \to H^{d_i}$ be the orthogonal projection of $H^d$ to $H^{d_i}$ ($H^{d_i}$
considered
as a subspace of $H^d$). Then the closed convex sets with nonempty interiors
$K:=\pi _1^{-1}(K_0)$ and $L:=\pi _2 ^{-1}(L_0)$ 
are unions of some point inverses under the maps $\pi _i$, which 
point inverses are copies of
$H^{d-d_1}$ and $H^{d-d_2}$. Then either $(\varphi K) \cap (\psi L) 
= \emptyset $, or the images by $\varphi $ and $\psi $ of
two such point inverses, which images are copies of $H^{d-d_1}$ and $H^{d-d_2}$,
intersect. In the second case 
by $(d-d_1) + (d-d_2) > d$ these images have a straight line in
common. So (1) of Theorem 1 is
satisfied vacuously. (Even we could have said
``compact 

\newpage

intersections'' or ``line-free intersections'', i.e., ones 
not containing straight lines.) We do not know a similar example when the
dimensions of open portions of $H^{e_1}$ in ${\text{bd}}\,K$ and of $H^{e_2}$
in ${\text{bd}}\,L$ satisfy $e_1 + e_2 \le d$. (Strict convexity of $K$ and
$L$ means $e_1 = e_2 = 0$. This is related to $i$-extreme
or $i$-exposed points of closed convex sets in ${\Bbb{R}}^d$, for $0 \le i \le
d-1$, cf. \cite{Sch}, Ch. 2.1.)
As already mentioned at the end of \S 1, for ${\Bbb{R}}^d$ where $d \ge 2$, 
for Theorem 2 the only additional example is two parallel slabs.   

Possibly for $S^d$ and $H^d$ 
Theorem 4 holds without its hypotheses about smoothness and 
supporting spheres. (For ${\Bbb{R}}^d$ the solution is announced in \S 1.)
\enddefinition


In the proofs of our Theorems we will use some ideas of [H].

\head 3. Preliminaries \endhead

In $S^d$, when saying {\it{ball}} or {\it{sphere}}, 
we always mean one with radius at most $\pi /2$ (thus the ball is convex).
For $S^d$, ${\Bbb{R}}^d$ and 
$H^d$ we denote by $B(x,r)$ the {\it{closed ball of centre $x$
and radius}} $r$. For points $x,y$ in $S^d,{\Bbb{R}}^d$ and 
$H^d$, we write $[x,y]$, $(x,y)$ or line $xy$
for the {\it{closed or open 
segment with end-points $x,y$, or the line passing through the points $x,y$,}}
respectively (these will not be used for $x,y$
antipodal in $S^d$, moreover line $xy$ will not be used for $x=y$) 
and $|xy|$ for the {\it{distance of $x$ and $y$}}. (For
$x=y$ we have $(x,y) = \emptyset $.)
The coordinate planes in ${\Bbb{R}}^d$ will be called $\xi _1 \xi
_2$-{\it{coordinate plane,}} etc.

A closed convex set $K$ in $X=S^d$, ${\Bbb R}^d$, $H^d$ with non-empty interior
is {\it{strictly
convex}} if its boundary does not contain a non-trivial segment.
A boundary point $x$ of this set $K$
is an {\it{extreme point of $K$}}
if it is not in the relative interior of a segment contained in
${\text{bd}}\,K$. 
A boundary point $x$ of this set $K$ is an {\it{exposed point of $K$}}
if the intersection of $K$ and some supporting hyperplane of $K$ 
is the one-point set $\{ x \}$.

For hyperbolic plane geometry we refer to [Ba], [Bo], [L], [P], for 
geometry of hyperbolic
space we refer to [AVS], [C], and for elementary differential geometry we
refer to [St]. 

The space $H^d$ has two usual models, in the interior of the unit ball in
${\Bbb R}^d$, namely the collinear (Caley-Klein) model and the conformal
(Poincar\'e) model. In analogy, we will speak about collinear and conformal
models of $S^d$ in ${\Bbb R}^d$, meaning the ones obtained by central
projection (from the centre), 
or by stereographic projection (from the north pole) to the
tangent hyperplane of $S^d$,
at the south pole, in ${\Bbb R}^{d+1}$ (this being identified with
${\Bbb{R}}^d$). These exist of
course only on the open southern half-sphere, or on $S^d$ minus the north pole,
respectively. Their images are ${\Bbb R}^d$. 

A {\it{paraball in}} $H^d$ is a closed convex set bounded by a parasphere.

\newpage

The {\it{base hyperplane of a hypersphere in}} $H^d$ is the hyperplane, for
which the {\it{hypersphere is a (signed!) 
distance surface}}. It can be given also as the
unique hyperplane, whose infinite points coincide with those of the
hypersphere. 

In the proofs of our theorems by the 
{\it{boundary components of a set}} we will
mean the connected components of the boundary of that set. 

We shortly recall some two-dimensional concepts to be used later. In
$S^2$, $H^2$ there are the following (complete, connected, twice 
differentiable) curves of constant curvature (in $S^2$ meaning geodesic
curvature). In $S^2$ these are 
the circles, of radii $r \in (0, \pi /2]$, with (geodesic) curvature $\cot r
\in [0, \infty )$. In $H^2$, these are circles of radii $r \in (0, \infty )$,
with curvature ${\text{coth}}\,r \in (1, \infty )$, 
paracycles, with curvature $1$, and
{\it{hypercycles, i.e., distance lines, with (signed!)
distance $l > 0$ from their base lines}}
(i.e., the straight lines
that connect their points at infinity), with curvature 
${\text{tanh}}\,l \in (0,1)$,
and straight lines, with curvature $0$. 
Either in $S^2$ or in $H^2$ (and also in ${\Bbb
R}^2$, where we have circles and straight lines), 
each sort of the above curves have
different curvatures, and for one sort, with different $r$ or $l$, they also
have different curvatures.
The common name of these curves 
is, except for straight lines in ${\Bbb R}^2$ and $H^2$, 
{\it{cycles}}. In $S^2$ also a great
circle is called a {\it{cycle}}, but when speaking about straight lines, for
$S^2$ this will mean great circles.
An elementary method for the calculation
of these curvatures for $H^2$ cf. in [V].

\head 4. Proofs of our theorems \endhead

The proof of Theorem 1 will be broken up to several lemmas.

In our proofs there will be chosen several times {\it{sufficiently small
numbers}} $\varepsilon _i>0$. For one $\varepsilon _i$ there may be several
upper bounds.
Whenever there are several $\varepsilon _i$'s, we always
will tell which $\varepsilon _i$ is sufficiently small, 
for which given $\varepsilon _j$.


\proclaim{Lemma 1.1}
Let $X=H^d$. Let $K \subsetneqq H^d$ be a closed convex set with non-empty
interior, such that the connected components $K_i$
of \,${\text{{\rm{bd}}}}\,K$ 
are congruent hyperspheres, with common distance $\lambda > 0$ from
their base hyperplanes $K_{0i}$. Then the hyperplanes $K_{0i}$ bound a
non-empty closed
convex set $K_0$ (possibly with empty interior, and on the other closed
side of each $K_{0i}$ as $K_i$), and $K$ equals the parallel
domain of $K_0$ for distance $\lambda $.
\endproclaim


\demo{Proof}
It will be convenient to use the collinear model. Then the existence,
non-emptyness, closedness and convexity of $K_0$ are evident.

The parallel domain of $K_0$ for distance $\lambda $ 
contains the parallel domain of any $K_{0i}$ for distance $\lambda $.
Consider the parallel domain of $K_0$ for distance $\lambda $, 
which is {\it{closed and convex}}. (This follows from the 
inequality valid for any {\it{Lambert quadrangle,}} i.e., one 
which has three right angles: if $ABCD$ has right angles at $A,B,C$, then
$|AB| < |CD|$, cf. \cite{C}, or \cite{AVS}, p. 68, 3.4.) Thus the 
parallel domain of $K_0$ for
distance $\lambda $ contains all the
hyperspheres $K_i$, hence also their closed convex hull $K$.

\newpage

Conversely, also $K$ contains the parallel domain of $K_0$ for distance
$\lambda $.
Namely, on the one hand $K_0 \subset K$, hence $z \in K_0 \Longrightarrow
z \in K$. On the other hand, let
$z \not\in K_0$; then $z$
is separated from $K_0$ by some hyperplane $K_{0i(z)}$.
Let ${\text{dist}}(z,K_0) \le \lambda $.
Clearly ${\text{dist}}(z,K_0)$ is attained for some point $x \in K_{0i(z)}$
(and $x$ is the orthogonal projection of $z$ to $K_{0i(z)}$). 
Therefore ${\text{dist}}(z,K_{0i(z)}) = {\text{dist}}(z,K_0) \le
\lambda $, and then $z$ (lying outside of the ``facet'' $K_{0i(z)}$ of $K_0$)
lies between $K_{0i(z)}$ and $K_{i(z)}$, hence $z \in K$.
$ \blacksquare $
\enddemo


In the next Lemma 1.2 we use the notations $K,K_i,\lambda ,K_{0i},K_0$ 
from Lemma 1.1, and for $L$
another set satisfying the same properties as $K$ in Lemma 1.1, we use the
analogous notations $L_i,\lambda ,L_{0i},L_0$, as in
Lemma 1.1 for $K$. (The value of $\lambda > 0$ is the same for $K$ and $L$.) 


\proclaim{Lemma 1.2}
Let $X=H^d$ and let $K,K_i,\lambda ,K_{0i},K_0$ and $L,L_i,\lambda ,L_{0i},L_0$
be as written just before this lemma. Let 
$\varphi $ and $\psi $ be orientation preserving
congruences of $H^d$ to itself, 
such that the following hold.

{\rm{(1)}} The hyperplanes $\varphi K_{01}$ and $\psi L_{01}$ 
either have no common finite or infinite point, or have one
common infinite point but no other common finite or infinite point. 

{\rm{(2)}} The sets ${\text{\rm{int}}}(\varphi K_0)$ and 
${\text{\rm{int}}}(\psi L_0)$ lie on the opposite closed sides of $\varphi
K_{01}$ or $\psi L_{01}$, as $\psi L_{01}$ or $\varphi K_{01}$, 
respectively. If one or
both of these sets is/are empty, this requirement is considered as
automatically satisfied for the empty one/s of these sets.

{\rm{(3)}} Let $\varphi K_1$ and $\psi L_1$ denote that connected component 
of \,${\text{\rm{bd}}}(\varphi K)$ or ${\text{\rm{bd}}}(\psi L)$, 
whose base hyperplane is $\varphi K_{01}$ or 
$\psi L_{01}$. If there are two such connected components of 
${\text{\rm{bd}}}(\varphi K)$ or ${\text{\rm{bd}}}(\psi L)$, then
we mean that one which lies on the same side of
$\varphi K_{01}$ or $\psi L_{01}$, as $\psi L_{01}$ or $\varphi K_{01}$,
respectively. 

Then, letting $\varphi K_1^*$ and $\psi L_1^*$ be 
the two closed convex sets 
bounded by the hyperspheres $\varphi K_1$ and $\psi L_1$, we have
$$
(\varphi K) \cap (\psi L) = (\varphi K_1^*) \cap (\psi L_1^*).
$$
\endproclaim


\demo{Proof}
Observe that $\varphi K \subset \varphi K_1^*$ and $\psi L \subset \psi
L_1^*$, hence
$$
(\varphi K) \cap (\psi L) \subset (\varphi K_1^*) \cap (\psi L_1^*).
\tag 1.2.1
$$

For the converse inclusion it suffices to prove
$$
M := (\varphi K_1^*) \cap (\psi L_1^*) \subset \varphi K .
\tag 1.2.2
$$

\newpage

Namely, in the analogous way we prove $M \subset \psi L$, and then 
these two inclusions together will prove 
$$
M \subset (\varphi K) \cap (\psi L) .
\tag 1.2.3
$$

Now we show \thetag{1.2.2}.
The hyperplane $\varphi K_{01}$ cuts $H^d$ into two closed halfspaces
$\varphi H'$ and $\varphi H''$. 
One of them, say, $\varphi H'$ contains $\psi L_{01}$ in its interior.
Then we have two cases: a point $z \in M$ belongs either to $\varphi H'$ or 
to $\varphi H''$. These cases 
which will be settled separately.

Let
$$
z \in M \cap (\varphi H') .
\tag 1.2.4
$$
Then $z \in M \subset \varphi K_1^*$, hence $z \in M \cap (\varphi H') \subset
(\varphi K_1^*) \cap (\varphi H') \subset \varphi K $. Thus
$$
z \in \varphi K .
\tag 1.2.5
$$

Now let
$$
z \in M \cap (\varphi H'') .
\tag 1.2.6
$$ 
Then by $z \in \varphi H''$ we have that
$z$ lies outside of $\psi L_{01}$, with respect to $\psi L_0$ (i.e., on the
side where $\varphi K_0$ lies). That is, 
$$
z \,\,(\in M \subset \psi L_1^*) {\text{ lies between }} 
\psi L_{01} {\text{ and }} \psi L_1, {\text{ hence }} z \in \psi L .
\tag 1.2.7
$$ 
Then Lemma 1.1 (applied to $\psi L$) implies that
$$
{\text{dist}}(z, \psi L_0) \le \lambda .
\tag 1.2.8
$$

Clearly ${\text{dist}}(z, \psi L_0)$ is attained for some point $\psi y \in 
\psi L_{01}$ (and $\psi y$  
is the orthogonal projection of $z$ to $\psi L_{01}$). 
Then 
$z \in \varphi H''$ (cf. \thetag{1.2.6}) and 
$\psi y \in \psi L_{01} \subset \varphi H'$ 
imply that $[z,\psi y]$ intersects $\varphi K_{01}$ at some point 
$\varphi x \in \varphi K_{01}$. Then, also using \thetag{1.2.8},
$$
{\text{dist}}(z, \varphi K_{01}) \le |z(\varphi x)| \le |z(\psi y)| 
= {\text{dist}}(z,\psi L_0) 
\le \lambda .
\tag 1.2.9
$$ 
That is, $z$ lies in the parallel domain of  
$\varphi K_{01}$ for distance $\lambda $, 
and thus also in the parallel domain of $\varphi K_0$
for distance $ \lambda $, 
which equals $\varphi K$ by Lemma 1.1. Thus again
$$
z \in \varphi K ,
\tag 1.2.10
$$

\newpage

ending the proof of the lemma.
$ \blacksquare $
\enddemo


Now we are ready to prove


\proclaim{Lemma 1.3}
{\rm{(2)}} of Theorem 1 implies {\rm{(1)}} of Theorem 1.
\endproclaim


\demo{Proof}
For notational convenience we suppose both $\varphi $ and $\psi $ to be the
identity congruence.

Any intersection (not only a small one)
of two congruent balls, with non-empty interior,
is centrally symmetric, with centre of symmetry the midpoint of the segment
joining their centres.

Any compact intersection (not only a small one)
of two paraballs $K$ and $L$, with non-empty interior,
is centrally symmetric. In fact, the
infinite points of the two paraballs, say, $k$ and $l$, are different, since
else the intersection would not be compact. 
We consider the straight
line $kl$. Let the other points of ${\text{bd}}\,K$ and ${\text{bd}}\,L$ 
on $kl$ be $k'$ and $l'$.
We may suppose that $k' \ne l'$ and that
the order of the points on $kl$ is $k,l',k',k$ (else $K \cap L$ would have an
empty interior). Then the symmetry with respect to the
midpoint of the segment $k'l'$ interchanges $K$ and $L$, hence this midpoint
is a centre of symmetry of $K \cap L$.

There remain the cases when the connected components of the boundaries both
of $K$
and $L$ are congruent hyperspheres,
whose numbers are at least $1$, and at most countably infinite. 

For the case when the boundary components both of $K$ and $L$ 
are congruent hyperspheres, these
hyperspheres are distance surfaces for some distance $\lambda > 0$. 
Replacing these
hyperspheres by their base hyperplanes, we obtain closed convex sets $K_0$ and
$L_0$
(possibly one hyperplane, which has no interior points, but this makes no
difference).
Then by Lemma 1.1 the parallel domain of 
$K_0$ and of $L_0$, at distance $ \lambda $, equals $K$ and $L$, respectively.

Now we
show that two different hypersphere boundary
components of $K$ have a distance at least $2 \lambda $. 
In fact, if $x,y$ belong to
two different boundary components $K_i,K_j$
of $K$, then the segment $[x,y]$ intersects 
the respective base hyperplanes $K_{0i},K_{0j}$
in points $x_1,y_1$, with order $x,x_1,y_1,y$
on $[x,y]$. Then $|xy| \ge |xx_1| + |y_1y| \ge 2 \lambda $.

Now suppose that diam\,$(K \cap L) < 2l$. Observe that ${\text{bd}}\,
(K \cap L) \subset
({\text{bd}}\,K) \cup ({\text{bd}}\,L)$. Thus $K \cap L$ cannot contain points
from different boundary components of $K$, or of $L$. Therefore $K \cap L$
contains points of one boundary component $K_i$ of $K$ and of one boundary
component $L_j$ of $L$. The hyperspheres $K_i$ and $L_j$
bound (uniquely determined)
closed convex sets with interior points, say $K_i^*$ and $L_j^*$, 
containing $K$
and $L$. Then, by Lemma 1.2, $K \cap L =  K_i^* \cap L_j^*$.

That is, we have a compact intersection (with non-empty interior)
of two convex sets $K_i^*$ and $L_j^*$, bounded 
by congruent hyperspheres $K_i$ and $L_j$. Then the sets of infinite points of 
$K_i$ and $L_j$ are disjoint.

\newpage

Considering the collinear model, this implies that
the base hyperplanes $K_{0i}$ and $L_{0j}$ 
have no finite or infinite
points in common. Let us consider the segment realizing the distance of these
hyperplanes. Then the symmetry with respect to its 
midpoint interchanges $K_i^*$ and $ L_j^*$, hence this midpoint
is a centre of symmetry of $K_i^* \cap L_j^* = K \cap L$.
$\blacksquare $
\enddemo


Now we turn to the proof of the implication $(1) \Rightarrow (2)$ in Theorem
1.


We begin with a simple lemma. Observe that by (A) both $K$ and $L$ are smooth.


\proclaim{Lemma 1.4}
Let $K \subsetneqq S^d$ be a smooth
convex body. Then, unless $K$ is a half-sphere, $K$ has an exposed point.
\endproclaim


\demo{Proof}
We consider two cases: 

(1): either ${\text{diam}}\,K < \pi $, or 

(2): ${\text{diam}}\,K= \pi $.

In case (1) here the cone $C \subset {\Bbb{R}}^{d+1}$
with base $K$ and vertex $0$ is a convex body, and
the relative interiors of its generatrices contain no extreme points of $C$. 
However, $0$ is an extreme point of $C$, 
hence it is a limit of exposed points $c_i$ of $C$, cf. \cite{Sch}, 
Theorem 1.4.7, first statement 
(Straszewicz's theorem). These exposed points are in particular
extreme, hence
$0 = \lim c_i$ implies that for sufficiently large $i$ we have $0=c_i$, hence
$0$ is an exposed point of $C$. Therefore 
$K$ is contained in an open half-sphere. Let us suppose that
this half-sphere is the southern half-sphere. Then the collinear model is
defined in a neighbourhood of $K$, and the image $pK$ of $K$ in it
is a compact convex
set in the model ${\Bbb R}^d$ ($p$ maps the open southern half-sphere to
${\Bbb{R}}^d$, which is identified with
the tangent hyperplane of $S^d$ at the south pole). 
Such a set $pK$
has an exposed point $z$ (\cite{Sch}, above cited, second statement), 
thus for some hyperplane $H \subset {\Bbb{R}}^d$
we have $H \cap (pK) = \{ z \} $. Then
$H' := {\text{cl}} \left( (p^{-1}H) \cup (-p^{-1}H) \right)
$ is a hyperplane (large $S^{d-1}$) in $S^d$ such that $H' \cap K = \{ p
^{-1}z \} $. Then $p ^{-1}z$ is an exposed point of $K$. 

In case (2) here $K$ contains two antipodal points of $S^d$, and
we may suppose
that these are $e_{d+1}=(0,...,0, 1)$ and $-e_{d+1}$. 
Since $K$ is smooth at $e_{d+1}$,
therefore we may suppose that it has at $e_{d+1}$
as tangent hyperplane (in $S^d$)
$\{ (\xi _1,...,\xi _d,\xi _{d+1}) $
\newline
$\in S^d \mid \xi _1 = 0 \} $, and $K$ lies on the side 
$\{ (\xi _1,...,\xi _d,\xi _{d+1}) \in S^d \mid \xi _1 \ge 0 \} $ 
of this hyperplane.
For $k \in K \setminus \{ e_{d+1}, -e_{d+1} \} $ both shorter arcs 
${\widetilde{e_{d+1}k}}$ and ${\widetilde{(-e_{d+1})k}}$ lie in $K$. Therefore
$K$ consists of entire half-meridians, connecting $e_{d+1}$ and $-e_{d+1}$.
By the hypothesis about the tangent hyperplane of $K$ at $e_{d+1}$, 
each half-meridian, whose relative interior lies in the
open half-sphere given by $\xi _1 > 0 $, lies entirely in 
$K$. 
Therefore $K$ contains the closed half-sphere given by $\xi _1 \ge 0$. By 
hypothesis we have 
$K \subsetneqq S^d$, therefore $K$ is a half-sphere.
$\blacksquare $
\enddemo


\newpage

\proclaim{Lemma 1.5}
Suppose the hypotheses of Theorem 1, and suppose {\rm{(1)}} of Theorem 1. 
Then the following hold.

{\rm{(1)}} For any $x \in {\text{\rm{bd}}}\,K$ 
and any $y \in {\text{\rm{bd}}}\,L$ all sectional
curvatures exist, and are equal to the same non-negative constant, and in case
of $X={\Bbb{R}}^d$ and $X=H^d$, to the same positive constant.

{\rm{(2)}} 
For any $x \in {\text{\rm{bd}}}\,K$ and any $y \in {\text{\rm{bd}}}\,L$
there exists an $\varepsilon > 0$, such that $B(\varphi x, \varepsilon ) 
\cap \left(
{\text{\rm{bd}}}\,( \varphi K) \right) $ and $B(\psi y, \varepsilon ) \cap 
\left(
{\text{\rm{bd}}}\,( \psi L) \right) $ are rotationally symmetric with respect to
the normal of \,${\text{\rm{bd}}}\,( \varphi K)$ at $\varphi x$, 
and with respect to
the normal of ${\text{\rm{bd}}}\,( \psi L)$ at $\psi y$, respectively.
\endproclaim


\demo{Proof}
{\bf{1.}}
For $S^d$ 
by Lemma 1.4 either both $K$ and $L$ are halfspheres, when the statement of
this lemma is satisfied with sectional curvatures $0$ --- which case we may
further disregard --- 
or, e.g., $K$ has an exposed point $x$ --- which we may suppose.

For ${\Bbb{R}}^d$, by hypothesis,
e.g., $K$ has an extreme point $x$. Then $x$ is an extreme
point of $K \cap B(x,1)$ as well,
hence it is a limit of exposed points $x_i$ of $K \cap B(x,1)$, 
cf. \cite{Sch}, above cited).
For $|xx_i| < 1$ we
have that $x_i$
is an exposed point of $K$ as well. In fact, for, say, 
$x_i=0$ and $K \cap B(x,1)$ 
lying strictly above, say, 
the $\xi _1 \ldots \xi _{d-1}$-coordinate plane, except for
$x_i$, also $K$ lies
strictly above the $\xi _1 \ldots \xi _{d-1}$-coordinate plane, except for
$x_i$. Namely else by convexity of $K$
there would be points of $K \cap B(x,1)$ in any neighbourhood of $x_i$
below or on the $\xi _1 \ldots \xi _{d-1}$-coordinate plane 
and different from $x_i$.

For $H^d$ by $C^2_+$ (or
by hypothesis (B)) all boundary points of $K$ and $L$ are exposed.

Thus in $S^d$, ${\Bbb{R}}^d$ and $H^d$, we have that, e.g.,
$$
K {\text{ has an exposed point }} x \,\,(\in {\text{bd}}\,K) .
\tag 1.5.1
$$

{\bf{2.}}
Let 
$$
\cases
n {\text{ and }} m {\text{ denote the outer unit normals of }} K \\
{\text{and }} L, {\text{ at }} x \in {\text{bd}}\,K {\text{ and }} y \in 
{\text{bd}}\,L, {\text{ respectively,}}
\endcases
\tag 1.5.2
$$
where
$$
y \in {\text{bd}}\,L {\text{ is arbitrary.}}
\tag 1.5.3
$$
(Recall that we have $C^2$, or the weaker \thetag{A}, which still
implies smoothness.)

\newpage

$$
\cases
{\text{Let us choose a point, say, origin }} O \in X, {\text{ and let }}
e_0,f_0 {\text{ be}} \\
{\text{opposite unit vectors in the tangent space of }} X {\text{ at }} O. 
{\text{ Let us}} \\
{\text{choose orientation preserving congruences }} \varphi _0, \psi _0 
{\text{ of }} X, {\text{ such}} \\
{\text{that }} \varphi _0 x = \psi _0 y = O, {\text{ and the images (in 
the tangent bundles)}} \\
{\text{of }} n {\text{ or }} m {\text{ (by the maps induced by }} \varphi _0
{\text{ or }} \psi _0 {\text{ in the tangent}} \\
{\text{bundles) should be }} e_0 {\text{ or }} f_0, {\text{ respectively.}} 
\endcases
\tag 1.5.4
$$

Then $(\varphi _0 K) \cap (\psi _0 L) \supset \{ O
\} $. 
$$
{\text{Let }} g {\text{ be the geodesic from }} O {\text{ in the direction of
}} e_0 {\text{ (equivalently, of }}f_0 ).
\tag 1.5.5
$$
$$
\cases
{\text{Let us move }} \varphi _0 K {\text{ and }}
\psi _0 L {\text{ toward each other, so that their}} \\
{\text{points originally coinciding with }} 
O {\text{ should move on the}} \\
{\text{straight line }} g, 
{\text{to the respective new positions }} O_{\varphi K} {\text{ and }} O_{\psi
L}, \\
{\text{while we allow any rotations of them, independently of each}} \\
{\text{other, about the axis }} g. {\text{ We denote these new 
images by }} \varphi K \\
{\text{and }} \psi L, {\text{ and we denote the
images of }} n {\text{ or }} m {\text{ (by the maps}} \\
{\text{induced by }} \varphi {\text{ or }} \psi  {\text{ in the
tangent bundles) by }} e {\text{ or }} f, \\
{\text{respectively, which are the outer unit normals of }} \varphi K 
{\text{ and}} \\
\psi L, {\text{ at }} \varphi x \in {\text{bd}}(\varphi K) {\text{ and }}
\psi y \in {\text{bd}}(\psi L), {\text{ respectively. Then }} g \\ 
{\text{coincides with the line }} O_{\varphi K}O_{\psi L}, {\text{ and }} 
O_{\varphi K} = \varphi x {\text{ and }} O_{\psi L} = \psi y.
\endcases
\tag 1.5.6
$$
Let the amount of the moving of the points originally coinciding with $O$,  
both for $\varphi _0 K$ and $\psi _0 L$, be a common small
distance 
$$
\cases
|OO_{\varphi K}|=|OO_{\psi L}|=\varepsilon _1 \in (0, \min \{ \varepsilon
_1(x), \varepsilon _1(y) \} /2), \\
{\text{consequently }} O {\text{ is the midpoint of }} [O_{\varphi K},O_{\psi
L}] .
\endcases
\tag 1.5.7
$$
Then by \thetag{A} $O_{\varphi K}$ and $O_{\psi L}$ lie in the
balls of radii $\varepsilon _1 (x)$ and $\varepsilon _1 (y)$ from \thetag{A},
hence by \thetag{A} and \thetag{1.5.7} 
$$
B(O_{\varphi K}, \varepsilon _1) = B(\varphi x, \varepsilon _1) \subset 
\psi L {\text{ and }} 
B(O_{\psi L}, \varepsilon _1) = B(\psi y, \varepsilon _1) \subset 
\varphi K .
\tag 1.5.8
$$

Then 

\newpage

$$
C:=(\varphi K) \cap (\psi L)
\tag 1.5.9
$$ 
has a non-empty interior, and, by
exposedness of $x$ in $K$ and 
convexity of $K$ and $L$, has an
arbitrarily small diameter, {\it{for $\varepsilon _1>0$ sufficiently small}}. 
Whenever its 
diameter is less than $\varepsilon = \varepsilon (K,L) > 0$, then 
it has a {\it{centre of symmetry,}} $c$, say. 

Moreover, by \thetag{1.5.1}, \thetag{1.5.3} and \thetag{1.5.5},
$$
O_{\varphi K} = \varphi x \in {\text{bd}}\,(\varphi K) {\text{ and }} 
O_{\psi L} = \psi y \in {\text{ bd}}\,(\psi L).
\tag 1.5.10
$$
By \thetag{1.5.8} and \thetag{1.5.10} we have 
$$
O_{\varphi K} \in {\text{bd}} \left( (\varphi K) \cap 
(\psi L) \right) = {\text{bd}}\,C
{\text{ and }}
O_{\psi L} \in 
{\text{bd}} \left( (\varphi K) \cap (\psi L) \right) = {\text{bd}}\,C.
\tag 1.5.11
$$
\enddemo


{\centerline{$*$ $*$ $*$}}

\vskip.3cm

We break up the further proof of Lemma 1.5 to several parts, namely, Lemma
1.6 and Corollary 1.7, after proving which we immediately
return to the proof of Lemma 5, and finish it.


\proclaim{Lemma 1.6}
Under the hypotheses of Lemma 1.5, and
with the notations from the proof of Lemma 1.5 above, we have the following.
Either

{\rm{(1)}} 
$X=S^d$, and both $K$ and $L$ are half-spheres, when the statement of Lemma
1.5 is satisfied with sectional curvatures $0$, or

{\rm{(2)}} for $\varepsilon _1 > 0 $ sufficiently small, 
the points $O \in X$, i.e., the origin in $X$, 
and the centre of symmetry $c$ of $C:=(\varphi K) \cap (\psi L)$ coincide. 
(For $S^d$ we mean one of the two antipodal centres of symmetry.)
\endproclaim


\demo{Proof}
First observe that, {\it{for $\varepsilon _1 > 0$
sufficiently small}}, we have by hypothesis $C^2$ (or its weakening \thetag{A}) 
of the theorem that 
$$
B(O, \varepsilon _1) \subset C = (\varphi K) \cap (\psi L) .
\tag 1.6.1
$$

We are going to show that $B(O, \varepsilon _1)$ is the unique ball of
maximal radius, contained in $C$.

We distinguish three cases: $X=S^d$, $X={\Bbb{R}}^d$ and $X=H^d$.

{\bf{1.}}
First we deal with the case of $S^d$. 

By Lemma 1.4 either 

(1) both $K$ and $L$ are halfspheres, when (1) of
this Lemma is satisfied, or, 

(2) e.g., $K$ has an exposed point $x$. 

Further in this proof we deal with this case (2), and we
are going to prove (2) of this lemma in this case (2), for each of $S^d$ (in
{\bf{1}}), ${\Bbb{R}}^d$ (in {\bf{2}}) and $H^d$ (in {\bf{3}}).

\newpage

Let
$\varphi K'$ or $\psi L'$ denote the half-$S^d$ containing $\varphi K$ or
$\psi L$, and containing
$O_{\varphi K} = \varphi x$ or 
$O_{\psi L} = \psi y$ in its boundary, and thus being there
tangent to ${\text{bd}}\,(\varphi K)$, or to ${\text{bd}}\,(\psi L)$, 
respectively. By
$\varphi K \subset \varphi K'$ and $\psi L \subset \psi L'$, we have 
also 
$$
C = (\varphi K) \cap (\psi L) \subset (\varphi K') \cap (\psi L') . 
\tag 1.6.2
$$

We are going to show that
$$
\cases
(\varphi K') \cap (\psi L') {\text{ contains a unique ball}} \\
{\text{of maximal radius, namely }} B(O, \varepsilon _1).
\endcases
\tag 1.6.3 
$$
In fact, we may suppose that $\left( {\text{bd}}\,(\varphi K') \right) 
\cap \left( {\text{bd}}\,(\psi L') \right) $ (a large $S^{d-2}$)
lies in the $\xi _3 \ldots \xi _{d+1}$-coordinate 
plane. Then
any point in  $(\varphi K') \cap (\psi L')$ has the same Euclidean
distances to ${\text{bd}}\,(\varphi K')$ and to ${\text{bd}}\,(\psi L')$ 
as its orthogonal
projection to
the $\xi _1 \xi _2$-coordinate plane has to the orthogonal projections 
of ${\text{bd}}\,(\varphi K')$ and of ${\text{bd}}\,(\psi L')$
to the $\xi _1 \xi _2$-coordinate plane. 
These last projections are
two lines containing the origin and enclosing
an angle $2 \varepsilon _1$, in the $\xi _1 \xi _2$-coordinate plane. 
By elementary geometry, in the 
sector of the unit circle bounded by these two lines, which is the orthogonal
projection of $(\varphi K') \cap (\psi L')$ to the $\xi _1 \xi _2$-coordinate 
plane, 
the maximum of the
distances to these two lines is maximal exactly for the point $O^*$
of this sector
which is the intersection of $S^1$ ($=S^d \cap [\xi _1 \xi _2$-coordinate
plane$]$)
and the (inner) angular bisector of this
sector of circle. However, $O^*$ has exactly one preimage on $S^d$, 
for the above mentioned projection, namely $O$. This proves \thetag{1.6.3}.

By \thetag{1.6.1}, \thetag{1.6.2} and \thetag{1.6.3} we have
$$
\cases
C = (\varphi K) \cap (\psi L) {\text{ contains a unique ball}} \\
{\text{of maximal radius, 
namely }} B(O, \varepsilon _1).
\endcases
\tag 1.6.4
$$

\newpage

Thus the centre of symmetry $c$ of $C$ must coincide with $O$ (or possibly
with $-O$, but also
in that case one of the centres of symmetry is $O$), proving that
unless we have (1) of Lemma 1.6, we have (2) of
Lemma 1.6, for $X=S^d$.

{\bf{2.}}
Now we turn to the case of ${\Bbb{R}}^d$. We write 
$$
O=(0, \ldots ,0,0), 
\,\,O_{\varphi K}=(0, \ldots , 0, - \varepsilon _1) 
{\text{ and }} O_{\psi L}=(0, \ldots
, 0, \varepsilon _1) .
\tag 1.6.5
$$
We recall from \thetag{1.5.6}
that $O_{\varphi K}$ and $O_{\psi L}$ span the
line $g$ from there, thus the opposite unit vectors $e$ and $f$ from
there are parallel to the $\xi _d$-axis. 
Then $e=(0, \ldots , 0,-1)$ and $f = (0, \ldots , 0,1)$.
Then $e$ and $f$ 
(being images of the unit outer normals $n$ of $K$ at $x$ and $m$ of $L$ at $y$)
are the unit
outer normals of $\varphi K$ at $\varphi x$ and of $\psi L$ at $\psi y$. 
(Since we deal with ${\Bbb{R}}^d$, the tangent spaces are obtained by
translation from each other, so we need not care about the difference of
$\varphi _0$ and $\varphi$, and similarly for $\psi _0$ and $\psi $.) Thus
the tangent hyperplanes of $\varphi K$ at $\varphi x$ and of $\psi L$ at $\psi
y$ (which exist by (A)) are parallel to the 
$\xi _1 \ldots \xi _{d-1}$-coordinate hyperplane.

Moreover, the tangent hyperplane of $\varphi K$ at $\varphi x = O_{\varphi K}$ 
is given
by $\xi _d = - \varepsilon _1$ and $\varphi K$ lies (non-strictly) 
above this hyperplane.
Similarly, the 
tangent hyperplane of $\psi L$ at $\psi y = O_{\psi L}$ is given
by $\xi _d = \varepsilon _1$ and $\psi L$ 
lies (non-strictly) below this hyperplane.
Therefore
$$
C = (\varphi K) \cap (\psi L) {\text{ lies in the parallel slab given by }}
- \varepsilon_1 \le \xi _d \le \varepsilon _1 .
\tag 1.6.6
$$
Hence any closed ball contained in $C$ is contained in the parallel slab
from \thetag{1.6.6}, hence has a radius at most $\varepsilon _1$. Moreover, it
has radius equal to $\varepsilon _1$ only if it touches both boundary
hyperplanes of this parallel slab.

Even, by exposedness of $\varphi x=O_{\varphi K}$ 
in $\varphi K$ (cf. \thetag{1.5.1}), 
for some support
hyperplane of $\varphi K$ at $\varphi x=O_{\varphi K}$ --- 
which is unique by (A), and
hence is given by $\xi _d = - \varepsilon _1$ --- we have that $(\varphi K)
\setminus \{ \varphi x \} $ lies strictly inside of this support hyperplane,
i.e.,
$$
\varphi K \subset 
\{ (\xi _1, \ldots , \xi _d) \in {\Bbb{R}}^d \mid \xi _d > - \varepsilon _1 \} 
\cup \{ O_{\varphi K} \} . 
\tag 1.6.7
$$
Hence if some closed ball of radius $\varepsilon _1$ is
contained in $C = (\varphi K) \cap (\psi L)$, then it touches the
hyperplane $\xi _d = - \varepsilon _1$. Also, this ball 
of radius $\varepsilon _1$ lies in $\varphi K$, hence the only point at which
it can touch the hyperplane $\xi _d = - \varepsilon _1$, 
is $O_{\varphi K} = (0, \ldots , 0,
- \varepsilon _1)$. Thus this ball is identical to $B(O, \varepsilon _1)$.
Thus also for ${\Bbb{R}}^d$ we have 

\newpage

$$
\cases
C = (\varphi K) \cap (\psi L) {\text{ contains a unique ball}} \\
{\text{of maximal radius, namely }} B(O, \varepsilon _1)
\endcases
\tag 1.6.8
$$
(like we had for $S^d$ in \thetag{1.6.4}).

Thus the centre of symmetry $c$ of $C$ must coincide with $O$, proving (2) of
Lemma 1.6 for $X={\Bbb{R}}^d$.
 
{\bf{3.}}
Now we turn to the case of $H^d$. 
Then, by hypothesis $C^2_+$ (or its weakening \thetag{B}) of 
the theorem, we have that for $\varepsilon _2 \in (0, \min \{ \varepsilon _2
(x), \varepsilon _2 (y) \} )$,
for a closed $\varepsilon _2$-neighbourhood 
$B(\varphi x,{\varepsilon _2}) \subset H^d$ of $\varphi x$ and 
$B(\psi y,{\varepsilon _2}) \subset H^d$ of $\psi y$ there holds 
$$
(\varphi K) \cap B(\varphi x,{\varepsilon _2}) 
\subset \varphi K'' {\text{ and }} 
(\psi L) \cap B(\psi y,\varepsilon _2) \subset \psi L'',
\tag 1.6.9
$$
where $\varphi K''$ and $\psi L''$ are closed
convex sets bounded by some hyperspheres of sectional curvatures at least
$\varepsilon _2(x)$ and $\varepsilon _2(y)$, respectively, with 
$$
\varphi x \in
{\text{bd}}\,( \varphi K'') {\text{ and }} \psi y \in {\text{bd}}\,( \psi L'').
\tag 1.6.10
$$
Since in \thetag{B} $\varepsilon _2 (x) > 0$ and $\varepsilon _2 (y) > 0$ 
can be decreased, preserving validity of \thetag{B}, therefore for our fixed
$\varphi x \in {\text{bd}}\,( \varphi K)$ and fixed $\psi y \in {\text{bd}}\,(
\psi L)$ we may assume without loss of generality that 
$$
\cases
\varphi K'' {\text{ and }} \psi L'' 
{\text{ are distance surfaces with equal distances }} 
\\ 
\varepsilon ' (x) = \varepsilon ' (y) \in (0, \varepsilon _2) 
{\text{ from their base hyperplanes.}}
\endcases
\tag 1.6.11
$$
(Further, recall from \S 2 that the sectional
curvatures of $\varphi K''$ and $\psi L''$  
and the distance for which they are distance
surfaces are asymptotically equal. The sectional curvatures are
${\text{tanh}}\,\varepsilon ' (x) = {\text{tanh}}\,\varepsilon ' (y)$.
In \S 2 this is stated only for $d=2$, but $\varphi K''$ and $\psi L''$ 
are rotationally
symmetric so all sectional curvatures are equal to that in the two-dimensional
case.) 

By positivity of the sectional curvatures of these hyperspheres we have
exposedness of $\varphi x \in {\text{bd}}\,(\varphi K)$ for $\varphi K$ 
and $\psi y \in {\text{bd}}\,(\psi L)$ for $\psi L$.

Moreover, by \thetag{1.6.10} and \thetag{1.6.11} there hold 
$$
\varphi x \in {\text{bd}}\,(\varphi K'') \subset \varphi
K'' {\text{ and }} \psi y \in
{\text{bd}}\,(\psi L'') \subset \psi L''
\tag 1.6.12
$$ 
and 
$$
\cases
{\text{bd}}\,(\varphi K'') {\text{ and }} {\text{bd}}\,(\psi L'') 
{\text{ have equal positive sectional}} \\
{\text{curvatures at }} \varphi x {\text{ and }} \psi y,
{\text{ which are less than }} {\text{tanh}}\, \varepsilon _2 < \varepsilon _2 .
\endcases
\tag 1.6.13
$$

\newpage

Further, $\varphi K''$ and $\psi L''$ contain
$\varphi x=O_{\varphi K}$ and $\psi y=O_{\psi L}$, 
and are there tangent to ${\text{bd}}\,(\varphi K)$, and to
${\text{bd}}\, (\psi L)$, respectively. 
Then necessarily $\varphi K''$ and $\psi L''$
have there their convex sides towards
${\text{int}}\,(\varphi K)$, or ${\text{int}}\,(\psi L)$,
respectively. 

By \thetag{1.6.9} and \thetag{1.5.6} we have
$$
(\varphi K) 
\cap  B(O_{\varphi K},{\varepsilon _2}) \subset \varphi K'' 
{\text{ and }}
(\psi L) 
\cap B(O_{\psi L},{\varepsilon _2}) \subset \psi L'' .
\tag 1.6.14
$$
Moreover, 
$$
\cases
{\text{if }} \varepsilon _1 {\text{ is sufficiently small for fixed }}
\varepsilon _2, {\text{ we have that }} \\
(\varphi K'') \cap (\psi L'') {\text{ has a sufficiently small diameter.}}
\endcases
\tag 1.6.15
$$
This body $(\varphi K'') \cap (\psi L'')$ is rotationally symmetric
about the line $O_{\varphi K}O_{\psi L}$, 
and also is symmetric with respect to the perpendicular
bisector plane of $[O_{\varphi K},O_{\psi L}]$. 
Its boundary consists of two geodesic $(d-1)$-balls on
${\text{bd}}\,(\varphi K'')$ and ${\text{bd}}\,(\psi L'')$, of centres
$O_{\varphi K}$ and
$O_{\psi L}$, respectively. 
Then also {\it{the (equal) geodesic radii of these two $(d-1)$-balls
are sufficiently small.}} 

As soon as these geodesic radii are less than
$\varepsilon _2$, then all points of these two geodesic
$(d-1)$-balls are at a distance
(in $H^d$) less than $\varepsilon _2$ from their centres
$O_{\varphi K}$ and $O_{\psi L}$. Then by \thetag{1.6.14} these two geodesic
$(d-1)$-balls are disjoint to
${\text{int}}\,(\varphi K)$ and ${\text{int}}\,(\psi L)$, respectively 
(else some points of
them would lie in ${\text{int}}\,(\varphi K'')$ or ${\text{int}}\,(\psi
L'')$, respectively, while they lie on  
${\text{bd}}\,(\varphi K'')$ or ${\text{bd}}\,(\psi L'')$, respectively).
Hence the union of these two geodesic $(d-1)$-balls is disjoint
to the intersection
$\left( {\text{int}}\,(\varphi K) \right) \cap \left( {\text{int}}\,(\psi
L) \right) = {\text{int}} \left( (\varphi K) \cap (\psi L) \right) $. 
Then the radial function of $(\varphi K) \cap (\psi L)$ with respect to $O$
is at most the radial function of $(\varphi K'') \cap (\psi L'')$ 
with respect to $O$. This implies
$$
C = (\varphi K) \cap (\psi L) \subset (\varphi K'') \cap (\psi L'') .
\tag 1.6.16
$$ 

We assert that also for $H^d$ we have that
$$
\cases
C = (\varphi K) \cap (\psi L) {\text{ contains a unique ball}} \\
{\text{of maximal radius, namely }} B(O, \varepsilon _1)
\endcases
\tag 1.6.17
$$
(like we had for $S^d$ in \thetag{1.6.4} and for ${\Bbb{R}}^d$ in
\thetag{1.6.8}).
Observe that for $(\varphi K'') \cap (\psi L'')$ rather than $C$
(cf. \thetag{1.6.16})
this is sufficient to be proved for $d=2$. Namely,
using \thetag{1.6.16}, the
(one-dimensional) axis of rotation $O_{\varphi K}O_{\psi L}$
of $(\varphi K'') \cap (\psi L'')$ and the
centre of a ball of maximal radius contained in  $(\varphi K'') \cap (\psi
L'')$ are contained in a $2$-plane of $H^d$.

Then $(\varphi K'') \cap (\psi L'')$ has as axis of symmetry the orthogonal
bisector line $g^*$ of $[O_{\varphi K},O_{\psi L}]$, 
and $O \in g^*$. Say, $g^*$ is horizontal, 
and $O_{\psi L}$ lies above $O_{\varphi K}$. 
Consider a circle of maximal radius contained in
$(\varphi K'') \cap (\psi L'')$; say, {\it{its centre $x$}}
lies (not strictly) above $g^*$. For contradiction, suppose $x \ne O$.

\newpage

Let $\psi L'''$ be the base line of $\psi L''$ (i.e., $\psi L''$ 
is a distance line for $\psi L'''$). 
Clearly, the straight
line $g = OO_{\psi L}$ is orthogonal to $\psi L''', l^*$ and
$\psi L''$ (these last three curves being distinct, and their
intersections with the straight 
line $g$ follow each other in the given order, from downwards to upwards). 
Let $\pi $ denote the
{\it{orthogonal projection of $H^2$ to $\psi L'''$}}. Let 
$\varrho (x)$ and $\sigma (x)$ denote the points of intersection
of $\left( {\text{bd}}\, [
(\varphi K'') \cap (\psi L'') ] \right) \cap \left( {\text{bd}}\,(\psi
L'') \right) $ and of $g^*$
with the straight line passing through $x$ and orthogonal to $\psi L'''$, 
respectively.

If $x$ lies on the line $OO_{\psi L}$ above $O$, then by \thetag{1.5.7} we
have
$$
|xO_{\psi L}| < |OO_{\psi L}| = \varepsilon _1 .
\tag 1.6.18
$$ 
Else we have 
$$
|x \varrho (x)| \le | \sigma (x) \varrho (x)| = | \pi (x) \varrho (x) | -
| \pi (x) \sigma (x) | .
\tag 1.6.19
$$ 
Here 
$$
| \pi (x) \varrho (x) | = |\pi (O_{\psi L}) O_{\psi L} |
\tag 1.6.20
$$ 
is the distance for which $\psi L''$ is the distance line for
$\psi L'''$. On the other hand, $[\sigma (x), \pi (x)]$ is an edge of 
the Lambert
quadrangle $O \pi (O_{\psi L}) 
\pi (x) \sigma (x)$, which has right angles at its
vertices $O, \pi (O_{\psi L})$ and $ \pi (x)$. (A Lambert quadrangle 
is a quadrangle with three right angles, cf. the proof of Lemma 1.1.) 
For the sides of this Lambert
quadrangle there holds 
$$
| \pi (x) \sigma (x)| > | \pi (O_{\psi L}) O | ,
\tag 1.6.21
$$
cf. \cite{C}, or \cite{AVS}, p. 68, 3.4.
Then
by \thetag{1.6.19}, \thetag{1.6.20} and \thetag{1.6.21} we get
$$
\cases
|x \varrho (x)| \le | \pi (x) \varrho (x) | - | \pi (x) \sigma (x) | < \\
| \pi (O_{\psi L}) O_{\psi L} | - | \pi (O_{\psi L}) O | 
= | OO_{\psi L} | = \varepsilon _1 ,
\endcases
\tag 1.6.22
$$
so \thetag{1.6.17} is proved.

Thus, as in the cases of $S^d$ and ${\Bbb{R}}^d$, also for $H^d$
the centre of symmetry $c$ of $C$ must coincide with $O$, proving (2) of Lemma
1.6 for $X=H^d$.

{\bf{4.}}
Thus the assertion of Lemma 1.6, either (1) or (2), 
is proved for each of $S^d$, ${\Bbb{R}}^d$ and
$H^d$, ending the proof of Lemma 1.6.
$\blacksquare $
\enddemo


\proclaim{Corollary 1.7}
{\rm{(i)}}
Let $X=S^d$. Then under the hypotheses of Lemma 1.5, and

\newpage

with the notations from the proof of Lemma 1.5 above, we have either that

{\rm{(1)}} 
both $K$ and $L$ are
halfspheres (when {\rm{(2)}} of Theorem 1 holds), or that

{\rm{(2)}}
both $K$ and $L$ are strictly convex. 

{\rm{(ii)}}
Let $X={\Bbb{R}}^d$. Then under the hypotheses of Lemma 1.5, and
with the notations from the proof of Lemma 1.5 above, we have that
both $K$ and $L$ are strictly convex. 
\endproclaim


\demo{Proof}
We have either $X=S^d$, and that (i) (1) of this Corollary
holds, which case we further disregard,
or else both for $S^d$ and ${\Bbb{R}}^d$,
recall that in \thetag{1.5.1}
$x \in {\text{bd}}\,K$ was chosen as an exposed point of $K$.
By Lemma 1.6, either 

(1) $X=S^d$, and (i) (1) of this Corollary holds,
which case was disregarded just above, or 

(2) for $\varepsilon _1 > 0 $ sufficiently small, 
$C = (\varphi K) \cap (\psi L)$ is centrally symmetric with
respect to $O$. 
Further in this proof we deal with this case (2).

Recall that $O_{\varphi K} = \varphi x \in 
{\text{bd}}\,(\varphi K)$ 
and $O_{\psi L} = \psi y \in {\text{bd}}\,(\psi L) $ 
are images of each other under this central symmetry,
cf. \thetag{1.5.1}, \thetag{1.5.3},
\thetag{1.5.6} and \thetag{1.5.7}.

Now recall from \thetag{1.5.8} and \thetag{1.5.10} that
$$
O_{\varphi K} \in \left( {\text{bd}}\,(\varphi K) \right)
\cap \left( {\text{int}}\,(\psi L) \right) {\text{ and }}
O_{\psi L} \in \left( {\text{bd}}\,(\psi L) \right) \cap \left( 
{\text{int}}\,(\varphi K) \right) .
\tag 1.7.1
$$
This implies that 
$$
O_{\varphi K},O_{\psi L} 
\in {\text{bd}}\,C = {\text{bd}} \left( (\varphi K) \cap (\psi L)
\right)
\tag 1.7.2
$$
and 
$$
\cases
{\text{for some }} \varepsilon > 0  {\text{ we have that }} 
B(O_{\varphi K}, \varepsilon ) \cap
\left( {\text{bd}}\,(\varphi K) \right) = \\
B(O_{\varphi K}, \varepsilon ) \cap
{\text{bd}} \left( (\varphi K) \cap (\psi L) \right) 
{\text{ and }}
B(O_{\psi L}, \varepsilon ) \cap
{\text{bd}} \left( (\varphi K) \cap (\psi L) \right) \\
= B(O_{\psi L}, \varepsilon ) \cap
\left( {\text{bd}}\,(\psi L) \right) {\text{ are also  
centrally symmetric}} \\
{\text{images of each other with respect to }} O 
\endcases
\tag 1.7.3
$$
(by \thetag{1.5.8} $\varepsilon \in (0, \varepsilon _1)$ suffices for this,
for $\varepsilon _1$ from \thetag{1.5.8}).

Since $x \in {\text{bd}}\,K$ is an exposed point of $K$ (cf. \thetag{1.5.1}), 
also
$O_{\varphi K} = \varphi x$ (cf. \thetag{1.5.6})
is an exposed point of $\varphi K$. By $O_{\varphi K} \in C \subset 
\varphi K$ (Lemma 1.6, (2))
then $O_{\varphi K}$ is an exposed point of $C$. By central symmetry of $C$
with respect to $O$ (cf. Lemma 1.6, (2)), also using \thetag{1.5.7},
then also 
$$
O_{\psi L} {\text{ is an exposed point of }}C .
\tag 1.7.4
$$ 
We claim that then 

\newpage

$$
O_{\psi L} = \psi y {\text{ is an exposed point of }} \psi L {\text{ as well.}}
\tag 1.7.5
$$
Recall that by (A) 
$$
\psi L {\text{ is smooth, hence has a tangent hyperplane }} H {\text{ at }}
O_{\psi L} = \psi y,
\tag 1.7.6
$$
and also recall \thetag{1.5.6}.

Suppose the contrary: $\psi L$ has a point $p$ outside of $H$ or 
on $H$ but different from $O_{\psi L}$. 
Then $\psi L$, being convex, would contain $[O_{\psi L}, p]$, 
thus $\psi L$ would
have a point $q \in [O_{\psi L}, p]$, outside of $H$ or 
on $H$ but different from $O_{\psi L}$, 
and additionally $q$ being arbitrarily close to
$O_{\psi L}$. (This holds even in the case when $X=S^d$ and 
$p$ is antipodal to $O_{\psi L}$. Namely
then we take some point $r \in (\psi L) \setminus \{ O_{\psi L}, p \} $, 
and then
$[O_{\psi L}, r] \cup [r,p] \subset \psi L$ can play the role of $[O_{\psi L}, 
p]$ from
above.)

However, in some neighbourhood of $O_{\psi L}$ we have that
${\text{bd}}\,(\psi L)$ and ${\text{bd}}\,C$ 
(and also $\psi L$ and $C$) coincide (recall $O_{\psi L} \in
{\text{int}}\,(\varphi K)$ from \thetag{1.7.1}). Then $C$ has a tangent plane
at $O_{\psi L}$, which can be defined locally, hence it coincides with $H$. This
tangent plane is the unique supporting plane of $C$ at $O_{\psi L}$, 
and however 
$C \cap [O_{\psi L},p]$ (or $C \cap ([O_{\psi L},r] \cup [r,p])$ for $X=S^d$ and
$p$ antipodal to
$O_{\psi L}$) contains points $q$ outside of $H$ or 
on $H$ but different from $O_{\psi L}$, 
and additionally $q$ being arbitrarily close to
$O_{\psi L}$. Then also $q \,\,(\in C)$ is outside of $H$ or is
on $H$ but is different from $O_{\psi L}$, contradicting \thetag{1.7.4}.
This contradiction ends the proof of our claim \thetag{1.7.5}.

Recapitulating: by {\thetag{1.5.1}} and {\thetag{1.7.5}}
$$
\cases
{\text{exposedness of }} x \in {\text{bd}}\,K {\text{ with respect to }} K
{\text{ implies}} \\
{\text{exposedness of }} y \in {\text{bd}}\,L {\text{ with respect to }} L . 
\endcases
\tag 1.7.7
$$

Recall from \thetag{1.5.1} and \thetag{1.5.3} that $x$ was an exposed point of
$K$ and $y$ was an arbitrary boundary point of $L$. Then by \thetag{1.7.7}
each boundary point $y$ of $L$ is an exposed point of $L$, i.e., $L$ is
strictly convex. 

In particular, $L$ has an exposed point. Changing the roles of $K$ and $L$ we
obtain that also $K$ is strictly convex.
$\blacksquare $
\enddemo


\demo{Proof of Lemma 1.5, {\bf{continuation}}}
Recall \thetag{1.7.3} and consider the line $g$, i.e.. the line containing
$O_{\varphi K}, O, O_{\psi L}$
(cf. \thetag{1.5.4}, \thetag{1.5.6} and
\thetag{1.5.7}). Take some $2$-plane $P$ containing 
the straight line $g$. By \thetag{1.7.3} and $P \ni O$ 
$$
\cases
{\text{for some }} \varepsilon > 0 {\text{ we have that }} 
B(O_{\varphi K}, \varepsilon ) \cap
\left( {\text{bd}}\,(\varphi K) \right) \cap P
\\
{\text{and }} B(O_{\psi L}, \varepsilon ) \cap
\left( {\text{bd}}\,(\psi L) \right) \cap P 
{\text{ are also centrally}} \\
{\text{symmetric images of each other with respect to }} O .
\endcases
\tag 1.5.12
$$

\newpage

Observe that by \thetag{1.6.4}, \thetag{1.6.8} and \thetag{1.6.15}, 
also using \thetag{1.5.6} and \thetag{1.5.7},
for each of
$S^d$, ${\Bbb{R}}^d$ and $H^d$ we have that the segment $[O_{\varphi K},
O_{\psi L}]$ is normal
to ${\text{bd}}\, (\varphi K)$ at $O_{\varphi K}$ and 
to ${\text{bd}}\, (\psi L)$ at $O_{\psi L}$ (${\text{bd}}\,(\varphi K)$ and
${\text{bd}}\,(\psi L)$ are smooth by (A)). 
Therefore both sets in \thetag{1.5.12}
are curves smooth at $O_{\varphi K}$ and at $O_{\psi L}$, respectively.

Therefore,
$$
\cases
{\text{the two curves from \thetag{1.5.12} have, at }} O_{\varphi K} 
{\text{ and }} O_{\psi L},\\
{\text{the same curvatures {\it{(sectional curvatures)}}, if one}} \\
{\text{of them exists, or they do not have curvatures there.}} 
\endcases
\tag 1.5.13
$$
Recall from \thetag{1.5.6} that
$\varphi $ and
$\psi $ were not determined uniquely, but at their definitions it was also
allowed that we applied any rotations to them, about the axis $g$, while $C$
is centrally symmetric with respect to $O$ (cf. Lemma 1.6, (2); recall that
Lemma 1.6, (1) gave $X=S^d$ and $K,L$ being half-spheres, which case was
disregarded at the beginning of the proof of Lemma 1.5). Therefore,
$$
\cases
{\text{for some }} \varepsilon > 0, \,\, B(O_{\varphi K}, \varepsilon )
\cap \left( {\text{bd}}\,(\varphi K) \right) {\text{ and }} \\
B(O_{\psi L}, \varepsilon )
\cap \left( {\text{bd}}\,(\psi L) \right)
{\text{ have }} g {\text{ as axis of rotation.}}
\endcases
\tag 1.5.14
$$
Now observe that $g$ is normal to $\varphi K$ at $\varphi x = O_{\varphi K}$, 
and to $\psi L$ at $\psi y = O_{\psi L}$ (cf. \thetag{1.5.10}), 
by \thetag{1.5.14} and smoothness of $K$ and $L$ (following from (A)).
This proves (2) of Lemma 1.5.

Then 
$$
\cases
{\text{either all sectional curvatures }} \left( {\text{i.e., the curvatures 
of all}} \right. \\
{\text{above curves in \thetag{1.5.12}}}, {\text{ for all }} 
2{\text{-planes }} P {\text{ containing }} \\
\left. g \right) , {\text{ both of }} \varphi K {\text{ and }} \psi L, 
{\text{ at the points }} O_{\varphi K} = \varphi x {\text{ and}} \\
O_{\psi L} = \psi y {\text{ exist and are equal, or all of them do not exist.}} 
\endcases
\tag 1.5.15
$$
Recall that $x$ was an arbitrary exposed point of
$K$ (cf. \thetag{1.5.1})
and $y$ was an arbitrary boundary point of $L$ (cf. \thetag{1.5.3}). 
However, we already know by Corollary 1.7 
that, unless $X=S^d$ and
both $K$ and $L$ are halfspheres of $S^d$ (which case was
disregarded at the beginning of the proof of Lemma 1.5), that $K$ and $L$ are 
strictly convex. (There this is stated only for
$S^d$ and ${\Bbb{R}}^d$. However, for $H^d$ strict convexity of $K$ and $L$ 
follows from the hypotheses of Theorem 1 and of this lemma, 
namely from $C^2_+$, or from (B)).
Hence also $x$ can be any boundary point of $K$, as $y$ can be any boundary
point of $L$, independently of each other.

So either 

\newpage
 
a) all sectional curvatures of both $K$ and $L$ exist, 
at each boundary point $x$ of $K$ and $y$ of $L$, 
and they are equal, namely to some number $\kappa \ge 0$, 
or 

b) they do not exist anywhere. 

However, convex surfaces in ${\Bbb R}^d$
are almost everywhere twice differentiable
(more exactly, the functions having, locally, in a suitable coordinate system,
these graphs, have Taylor series expansions, of second degree, with error term
$o(\| x \| ^2)$ ---
cf. [Sch], pp. 31-32, for ${\Bbb{R}}^d$, that extends to $S^d$ and $H^d$ by
using their collinear models).
This rules out possibility b), so possibility a) holds, as stated in this 
lemma.
Clearly, for ${\Bbb R}^d$ and
$H^d$, the hypotheses of Theorem 1 and of this lemma imply $\kappa > 0$.

{\it{This ends the proof of Lemma 1.5.}}
$\blacksquare $
\enddemo


The later following Lemmas 1.8 and 1.9 will be used not only for the proof of
Theorem 1, but also for the proof of Theorem 4. Therefore the hypotheses 
of Lemmas 1.8 and 1.9 will contain
alternatively (1) of Theorem 1, or (1) of Theorem 4. Because of this
first we have to turn to the proof of Theorem 4, and lead it so far as we
have led the proof of Theorem 1 till now.
Thus we have to prove the necessary analogues of some of Lemmas 1.1 till 1.7, 
including the complete proof of Lemma 1.5, as Lemmas 4.1 till 4.3.
This we do in order to avoid unnecessary repetitions (of Lemmas 1.8 and
1.9). 


\demo{Proof of Theorem 4}
{\bf{1.}}
The implication $(2) \Rightarrow (1)$ of this Theorem
is evident: the midpoint of the (any)
segment connecting the centres of the balls $\varphi K$ and $\psi L$ 
is a centre of symmetry of
${\text{cl\,conv}} \left( (\varphi K) \cap (\psi L) \right)$.

{\bf{2.}}
Now we turn to the proof of the implication $(1) \Rightarrow (2)$ of this
Theorem.

Let $x \in {\text{bd}}\,K$ and $y \in {\text{bd}}\,L$. 
Let $S(x)$ and $S(y)$ denote supporting spheres of $K$ and $L$, at $x$ and
$y$, respectively, of radius less than $\pi /2 $ for $S^d$. 
Observe that increasing the radius of a
supporting sphere at $x$ or $y$, for $S^d$ to a value less than $\pi /2$, 
while retaining their outer unit normals at $x$ or $y$, 
preserves the supporting property in these
points. Therefore we may assume that these two supporting spheres $S(x)$ and
$S(y)$ have equal radii, and this common
radius for $S^d$ is less than $\pi /2$. Even, if we
increase the radius further, for $S^d$ to a value less than $\pi /2$,
we may suppose that 
$$
\cases
{\text{these supporting spheres }} S(x) {\text{ and }} S(y) {\text{ have }}
x {\text{ and }} y {\text{ as}} \\
{\text{the unique common points with }} K {\text{ and }} L ,
{\text{ respectively.}}
\endcases
\tag 4.1
$$
Now we write $K(x)$ and $L(y)$ for the balls bounded by $S(x)$ and $S(y)$,
respectively. The common radius of $K(x)$ and $L(y)$ is denoted by $R$  --- for
the case of $S^d$ we have $R < \pi /2$.

Clearly 

\newpage

$$
\cases
{\text{we may assume for }} S^d {\text{ that }} R < \pi /2 {\text{ is 
arbitrarily close}} \\
{\text{to }} \pi /2 , {\text{ and for }} {\Bbb{R}}^d {\text{ and }} H^d 
{\text{ that }} R {\text{ is arbitrarily large.}}
\endcases
\tag 4.2
$$
$$
\cases
{\text{Let }} B_0 {\text{ be a fixed ball of radius }} R {\text{ in }} S^d, 
{\Bbb{R}}^d {\text{ or }} H^d, 
{\text{ whose}} \\
{\text{centre is denoted by }} O. {\text{ Let us choose orientation 
preserving}} \\
{\text{congruences }} \varphi {\text{ and }} \psi , {\text{such that }} 
\varphi K(x) = \psi L(y) = B_0, {\text{ and}} \\
\varphi (x) {\text{ and }} \psi (y) {\text{ are }} {\text{antipodal points 
of }} {\text{bd}}\,B_0.
\endcases
\tag 4.3
$$

Observe that, like in \thetag{1.5.6}, also here 
$$
\cases
\varphi {\text{ and }} \psi  {\text{ are by their definition not determined
uniquely, but}} \\
{\text{we are allowed to apply any rotation to them, independently}} \\
{\text{of each other, about the axis }} g, {\text{ spanned by }} 
\varphi (x) {\text{ and }} \psi (y).
\endcases
\tag 4.4
$$

Then we have $\varphi K, \psi L \subset B_0$, hence ${\text{cl\,conv}}\left(
( \varphi K) \cup (\psi L) \right) \subset B_0$. 
Moreover, 
$$
\cases
{\text{since the diameter of }} B_0 {\text{ (as a convex body) 
is twice its radius }}\\
R, {\text{ the two points }} \varphi x {\text{ and }} \psi y 
{\text{ form a diametral pair of points}} \\
{\text{in the centrally symmetrical set }} {\text{cl\,conv}}\left( 
( \varphi K) \cup (\psi L) \right) .
\endcases
\tag 4.5
$$ 
By \thetag{4.5} we have that 
$$
{\text{\rm{diam}}}\,[{\text{\rm{cl\,conv}}}
\left( ( \varphi K) \cup (\psi L) \right) ] = 2R.
\tag 4.6
$$

{\centerline{$*$ $*$ $*$}}

\vskip.3cm

\newpage

We break up the further proof of Theorem 4 to several parts, namely, Lemmas
4.1, 4.2 and 4.3. After proving these we continue with the proof of Lemmas 1.8
and 1.9, both being necessary for the proof both 
of Theorems 1 and 4. Then we turn to
prove Theorem 1. The continuation, i.e., finishing of the
proof of Theorem 4 follows after the proof of Theorem 3.
\enddemo

The following Lemma 4.1 will be some analogue of Lemma 1.6, (2), inasmuch
in Lemma 4.1 we
determine the centre of symmetry of our set, which set is now 
${\text{cl\,conv}}\left( ( \varphi K) \cup (\psi L) \right)$ (while
in Lemma 1.6, (2) the set was $(\varphi K) \cap (\psi L)$).


\proclaim{Lemma 4.1}
Supposing the hypotheses of Theorem 4 and (1) of Theorem 4, and with the above
notations, the 
centre of symmetry of the centrally symmetrical set
${\text{\rm{cl\,conv}}}\left( ( \varphi K) \cup (\psi L) \right)$ is 
the centre $O$ of $B_0$. (For $S^d$ we mean one of the two antipodal centres
of symmetry.)
\endproclaim


\demo{Proof}
Observe that $B_0$ is a ball of radius $R$. Then 
$$
\cases
{\text{the diameter of }} B_0, {\text{ in the sense of convex bodies, is }} 
2R, {\text{ and is}} \\
{\text{attained exactly for antipodal pairs of points on its boundary}}.
\endcases
\tag 4.1.1
$$

For ${\Bbb{R}}^d$ we use its usual geometry, for $H^d$ we use its collinear
model, with $O$ at its centre. Thus for $H^d$ the image $B_0'$
of $B_0$ in the collinear model is a Euclidean 
ball with centre $O$, and also for
${\Bbb{R}}^d$ we have by \thetag{4.3} 
that $B_0' := B_0$ is a Euclidean ball with centre $O$.

For $S^d$ we also 
use the collinear model. Namely, supposing that $O$ is the south pole, we use
radial projection $ \pi $ 
from the centre of $S^d$ in ${\Bbb{R}}^{d+1}$ to the model
tangent $d$-plane of $S^d$ in ${\Bbb{R}}^{d+1}$ at the south pole $O$. (This
model exists for the open southern hemisphere.)
Thus, also for $S^d$, the image $B_0'$
of $B_0$ in the collinear model is a Euclidean ball with centre $O$.
Thus, also for $S^d$, for 
simplicity of notation, we will work in this model ${\Bbb{R}}^d$ (similarly as
for $H^d$). 

For points and sets in the models we will apply upper indices $'$. We recall
that

\newpage

$$
\cases
{\text{in }} {\Bbb{R}}^d {\text{ the closed convex hull of the union of
two compact convex}} \\
{\text{sets }} K',L' {\text{ is }} {\text{cl\,conv}}\,( K' \cup L' ) = 
{\text{conv}}\, (K' \cup L' ) = \bigcup \{ [ k',l' ] \mid  \\
k' \in K',\,\,l' \in L' \} = K' \cup L' \cup \left( \bigcup \{ (k',l') \mid k' 
\in K',\,\,l' \in L' \} \right) . 
\endcases
\tag 4.1.2
$$
(Recall from \S 3 
that $(x,y)$ is the open segment with endpoints $x,y$; in particular,
for $x=y$ it is $\emptyset $.)
We apply this to the convex bodies
$$
\cases
K' \subset B_0' {\text{ and }} L' \subset B_0', {\text{ which are the images of 
the}} \\
{\text{sets }} \varphi K {\text{ and }} \psi L {\text{ in the respective 
collinear models}}
\endcases 
\tag 4.1.3
$$
(by collinear model of ${\Bbb{R}}^d$ we mean ${\Bbb{R}}^d$ itself).

By the collinear models, for $S^d$ the statement corresponding to
\thetag{4.1.2} is valid for compact 
convex subsets $\varphi K, \psi L$ of the open southern
hemisphere, in particular for compact convex subsets of $B_0$. For
$H^d$ we work in its collinear model, contained in ${\Bbb{R}}^d$ (as the unit
ball of ${\Bbb{R}}^d$), 
while for ${\Bbb{R}}^d$ by its collinear model we mean ${\Bbb{R}}^d$ itself.
Hence 
$$
\cases
{\text{the statement corresponding to \thetag{4.1.2} is valid for }} \\
S^d, \,\,{\Bbb{R}}^d {\text{ and }} H^d, {\text{ for convex bodies contained 
in }} B_0.
\endcases
\tag 4.1.4
$$
By \thetag{4.1.2} and \thetag{4.1} we have
$$
\cases
[{\text{cl\,conv}} \left( (\varphi K) \cup (\psi L) \right) ] \cap
({\text{bd}}\,B_0) = \left( (\varphi K) \cup (\psi L) \right) \cap
({\text{bd}}\,B_0) \\
= \left( (\varphi K) \cap ({\text{bd}}\,B_0) \right) \cup
\left( (\psi L) \cap ({\text{bd}}\,B_0) \right) = 
\{ \varphi x, \psi y \} .
\endcases
\tag 4.1.5
$$
Then by \thetag{4.5}, \thetag{4.1.1} and \thetag{4.1.5}
$$
\cases
{\text{the diameter of }} {\text{cl\,conv}} \left( (\varphi K) \cup (\psi L) 
\right)  {\text{ is }} 2R, {\text{ and is}}\\
{\text{attained exactly for antipodal pairs of points of }} B_0 {\text{ on }} \\
{\text{bd}} \left( {\text{cl\,conv}} \left( (\varphi K) \cup (\psi L) 
\right) \right), {\text{ i.e., for the unique diametral}} \\
{\text{pair of points of }} {\text{cl\,conv}} \left( (\varphi K) 
\cup (\psi L) \right) , {\text{ i.e., for }} \{ \varphi x , \psi y \} .
\endcases 
\tag 4.1.6
$$
Thus a central symmetry of ${\text{cl\,conv}} \left( (\varphi K) \cup (\psi
L) \right)$ 
preserves the pair of points $\{ \varphi x , \psi y \}$. Hence its centre
of symmetry is the mid-point of the segment $[\varphi x , \psi y]$, i.e., the
point $O$.
(For $S^d$ we mean one of the two antipodal centres of symmetry --- namely the
one in the open southern hemisphere.) 
$\blacksquare $
\enddemo


The following Lemma 4.2 also is an analogue of some step in the proof of
Theorem 1. Namely in Lemma 1.6, (2) 
we had information about the central symmetry, with
respect to the point there denoted also by $O$, of the set $(\varphi K) \cap
(\psi L)$. Then in \thetag{1.7.3} 

\newpage

we could turn from the boundary of this
intersection to the boundaries of $\varphi K$ and $\psi L$. The same will
happen in Lemma 4.2, for the set ${\text{cl\,conv}} \left( (\varphi K) \cup
(\psi L) \right) $.
 

\proclaim{Lemma 4.2}
Supposing the hypotheses of Theorem 4 and (1) of Theorem 4, and using the above
notations,
in a neighbourhood of $\varphi x$ (or of \,$\psi
y$) the sets ${\text{\rm{bd}}}\,(\varphi K)$ (or ${\text{\rm{bd}}}\,(\psi L)$)
and \,${\text{\rm{bd}}}\,[{\text{\rm{cl\,conv}}}
\left( (\varphi K) \cup (\psi L) \right) ]$ 
coincide. In particular, for some $\varepsilon > 0$ we have that $B( \varphi x
, \varepsilon ) \cap \left( {\text{\rm{bd}}}\,( \varphi K) \right) $ and
$B( \psi y
, \varepsilon ) \cap \left( {\text{\rm{bd}}}\,( \psi L) \right) $ are also
centrally symmetric images of each other with respect to $O$.
\endproclaim


\demo{Proof}
{\bf{1.}}
By \thetag{4.1} we have $\varphi x \not\in \psi L$ and $\psi y \not\in \varphi
K$. Thus some neighbourhoods of $\varphi x$ and of $\psi y$ do not intersect
$\psi L$ and $\varphi K$, respectively. Thus there are hyperplanes $P_x$ and
$P_y$ orthogonally intersecting the segment 
$[\varphi x, \psi y]$, in points sufficiently close
to $\varphi x$ or $\psi y$, and having the entire
$\psi L$ or $\varphi K$ on one side (on the side containing $\psi y$ or
$\varphi x$, respectively).

Recall \thetag{4.1.4}, which
will permit us to work further in the proof of this Lemma 
in the collinear model, the Euclidean space ${\Bbb{R}}^d$ or its open unit ball.

{\bf{2.}}
We write the usual
basic unit vectors of ${\Bbb{R}}^d$ as $e_1, \ldots , e_d$. Then 
$$
\cases
B_0' {\text{ is, or can be supposed to be a ball with}} \\
{\text{centre the origin }} 0 {\text{ and with radius }} R', 
\endcases
\tag 4.2.1
$$
and we assume that 
$$
\cases
\pi (\varphi x) = R'e_d, {\text{ and }} \pi (\psi y) = -R'e_d, 
{\text{ where the map }} \pi  {\text{ associates}} \\
{\text{to points of }} S^d, \,\, 
{\Bbb{R}}^d {\text{ and }} H^d {\text{ their images in the respective model.}}
\endcases
\tag 4.2.2
$$ 
(By the collinear model of ${\Bbb{R}}^d$ we mean ${\Bbb{R}}^d$ itself.)
Since $[\pi (\psi y), \pi (\varphi x)] = [-R'e_d, R'$
\newline
$e_d]$ 
contains $0$ (which is, for $H^d$ and $S^d$,
the image of the centre $O$ in the model; also cf. \thetag{4.2.1}), 
therefore the image in the model
of the orthogonally intersecting plane $P_x$ is also an
orthogonally intersecting plane of 
$[-R'e_d, R'e_d]$ in the model ${\Bbb{R}}^d$ or its open unit ball, hence is
given by $\xi _d = R''$ where $R'' < R'$ (and $R' - R''$ is small). 
Therefore, with the notations from \thetag{4.1.3}, 
$$
L' \subset \{ (\xi _1, \ldots ,\xi _d) \mid \xi _d \le R'' \} .
\tag 4.2.3
$$

Observe that the smooth convex body $K' \subset B_0'$ (cf. \thetag{4.1.3})
has at $R'e_d$ the outer unit normal $e_d$,
thus for its support function $h_{K'}: S^{d-1} \to {\Bbb{R}}$ 
we have $h_{K'}(e_d) =
R'$. On the other hand, by \thetag{4.2.3}, with the analogous notation, we have
$h_{L'}(e_d) \le R'' < R'$. Since the
support functions are continuous, therefore 

\newpage

$$
\cases
{\text{for some neighbourhood of }} e_d {\text{ in }} S^{d-1} \\
{\text{we have the inequality }} h_{L'}( \cdot ) < h_{K'}(\cdot ).
\endcases
\tag 4.2.4
$$

Now recall that $K'$ is smooth. Therefore it is even $C^1$, and hence for some
$\varepsilon > 0$, in the open $\varepsilon $-neighbourhood
$U(R'e_d, \varepsilon )$ of $R'e_d$ (in ${\Bbb{R}}^d$) we have that
$({\text{bd}}\,K') \cap U(R'e_d, \varepsilon )$ 
is a connected smooth manifold with
outward unit normals very close to $e_d$. The support function of
${\text{conv}}\,(K' \cup L') =
{\text{cl\,conv}}(K' \cup L')$ (cf. \thetag{4.1.2}) 
is the pointwise maximum of the support functions $h_{K'}$ and
$h_{L'}$. In particular, for points of a subset of $S^{d-1}$
where we have $h_{L'}( \cdot ) < h_{K'} ( \cdot )$, the
support sets of $K'$ and of ${\text{cl\,conv}}\,(K' \cup L')$ coincide. 
This implies by \thetag{4.2.2} and \thetag{4.2.4} that for some $\delta > 0$
$$ 
\cases
({\text{bd}}\,K') \cap U \left( \pi (\varphi x), \delta \right) =
({\text{bd}}\,K') \cap U(R'e_d, \delta ) \\
= [ {\text{bd\,cl\,conv}}\,(K' \cup L') ] \cap U(R'e_d, \delta ) \\
= [ {\text{bd\,cl\,conv}}\,(K' \cup L') ] \cap 
U \left( \pi (\varphi x), \delta \right) .
\endcases
\tag 4.2.5
$$ 
Turning back from the sets in the model to the original sets in $X = S^d,
{\Bbb{R}}^d, H^d$, we obtain 
the statement of the Lemma for $\varphi x$. In fact, $\pi ^{-1} U \left( \pi
(\varphi x), \delta \right) $ is a neighbourhood of $\varphi x$, hence contains
an open ball in $X$ with centre $\varphi x$ and radius some 
$\varepsilon > 0$.

The analogous statement about coincidence of the intersections of some open ball
of centre $\psi y$, 
with ${\text{bd}}\,( \psi L)$ and ${\text{bd\,cl\,conv}} \left( (\varphi K) 
\cup ( \psi L) \right) $, follows analogously.

This proves the first sentence of this Lemma. The second sentence of this
Lemma is an immediate consequence of its first sentence and of Lemma 4.1.
$ \blacksquare $
\enddemo


The following Lemma 4.3 is an analogue of Lemma 1.5.


\proclaim{Lemma 4.3}
Suppose the hypotheses of Theorem 4, and suppose {\rm{(1)}} of Theorem 4. Then
both conclusions {\rm{(1)}} and {\rm{(2)}} 
of Lemma 1.5 hold. Moreover, the constant
sectional curvatures in Lemma 1.5 {\rm{(1)}} 
are positive for $S^d$ and ${\Bbb{R}}^d$, and are greater
than $1$ for $H^d$.
\endproclaim


\demo{Proof} 
In analogy with \thetag{1.7.3}, we have here Lemma 4.2, second sentence. 
As in the proof of
Lemma 1.5, continuation, we consider the line $g$ containing $\varphi x$ and
$\psi y$, hence also $O$, as the midpoint of the segment $[\varphi x, \psi
y]$ (cf. \thetag{4.1.6} --- for $S^d$ we mean the midpoint 
in the open southern hemisphere). 
As in the proof of Lemma 1.5, continuation, we take some $2$-plane $P$
containing the straight line $g$.
Then, in analogy with \thetag{1.5.12}, by Lemma 4.2 we have
$$
\cases
{\text{for some }} \varepsilon > 0 {\text{ that }} B(\varphi x, 
\varepsilon) \cap \left( {\text{bd}}\,(\varphi K) \right) 
\cap P \\
{\text{and }} 
B(\psi y, \varepsilon) \cap \left( {\text{bd}}\,(\psi L) \right) \cap 
P {\text{ are also centrally}}\\
{\text{symmetric images of each other with respect to }} O .
\endcases
\tag 4.3.1 
$$

\newpage

Here, by the hypotheses of Theorem 4, 
the two sets in \thetag{4.3.1} are smooth curves.

Again as in the proof of Lemma 1.5, continuation, \thetag{1.5.13}, here
by \thetag{4.3.1}
$$
\cases
{\text{the two curves in \thetag{4.3.1} have, at }} \varphi x  {\text{ and }} 
\psi y, {\text{ the}} \\
{\text{same curvatures {\it{(sectional curvatures}}), if one of}} \\
{\text{them exists, or they do not have curvatures there.}} 
\endcases
\tag 4.3.2
$$

Recall that Theorem 4 has as hypothesis the existence of support spheres at
any boundary point of $K$ and $L$, for $S^d$ of radius less than $ \pi /2$. 
Take into account that spheres, for $S^d$ of radius
less than $ \pi /2$, have positive sectional
curvatures for $S^d$ and ${\Bbb{R}}^d$, 
and have sectional curvatures greater than $1$ for $H^d$. 
Therefore each (existing) sectional
curvature of $K$ and $L$ has to be at least the sectional
curvature of some sphere.
Now recall that the sectional curvatures of a sphere 
have the strict lower bounds stated in this Lemma. 
Therefore
$$
\cases
{\text{any existing sectional curvature of }} \varphi K {\text{ at }} \varphi
x \in {\text{bd}}(\varphi K) {\text{ and of}}\\
\psi L {\text{ at }} \psi y
\in {\text{bd}}(\psi L) {\text{ is positive, and for }} H^d {\text{ is greater 
than $1$.}}
\endcases
\tag 4.3.3
$$

Analogously as in \thetag{1.5.6}, here we have \thetag{4.4}, which implies
the analogue of \thetag{1.5.14}, namely
$$
\cases
{\text{for some }} \varepsilon > 0,\,\,B(\varphi x, \varepsilon ) \cap \left(
{\text{bd}} (\varphi K) \right) {\text{ and}} \\
B(\psi y, \varepsilon ) \cap \left( {\text{bd}}\,(\psi L) \right)
{\text{ have }} g {\text{ as an axis of rotation.}}
\endcases
\tag 4.3.4
$$
Now observe that $g$ is normal to $\varphi K$ at $\varphi x$, and to $\psi L$
at $\psi y$, by \thetag{4.3.4} and smoothness of $K$ and $L$.
This proves (2) of Lemma 1.5, as stated in Lemma 4.3. 

Then, by \thetag{4.3.1} and \thetag{4.4} (as in \thetag{1.5.15})
$$
\cases
{\text{either all sectional curvatures of }} \varphi K {\text{ and }} \psi L, 
{\text{ at the points}} \\
\varphi x \in {\text{bd}}\,(\varphi K) {\text{ and }} 
\psi y \in {\text{bd}}\,(\psi L) \,\,\left( {\text{i.e.,}} 
{\text{ the curvatures}} \right. \\
{\text{of all above curves in \thetag{4.3.1}}}, {\text{at the points }} \varphi 
x \in {\text{bd}}\,(\varphi K) \\
{\text{and }} \psi y \in {\text{bd}}\,(\psi L), {\text{ for all }} 
2{\text{-planes }} P {\text{ containing}} \\
\left. g \right) , {\text{ exist and are equal, or all of them do not exist.}}
\endcases
\tag 4.3.5
$$

From now on we turn to the original sets $K$ and $L$.

Again as in the proof of Lemma 1.5, continuation, varying $x $ in 
${\text{bd}}\,K$ and $y$ in ${\text{bd}}\,L$, independently of each other,
we have either that

a) all sectional curvatures of both $K$ and $L$ exist, at each boundary point
$x$ of $K$ and $y$ of $L$, and are equal --- namely to some number $\kappa > 0$
for $S^d$ and ${\Bbb{R}}^d$, and to some number
$\kappa > 1$ for $H^d$, by \thetag{4.3.3}, or that

b) they do not exist anywhere.  

\newpage

Then, as in the last paragraph of the proof of Lemma 1.5, almost everywhere
twice differentiability of convex surfaces rules out possibility b), so
possibility a) holds, as stated in this lemma.
This proves (1) of Lemma 1.5, as stated in Lemma 4.3. 

The statement of Lemma 4.3
about the inequalities for the sectional curvatures follows from \thetag{4.3.3}.
This ends the proof of Lemma 4.3. 
$ \blacksquare $
\enddemo


Now we can turn already to the proofs of Lemmas 1.8 and 1.9. These lemmas will
be common tools for the proofs of Theorem 1 and Theorem 4. 
Therefore the hypotheses of Lemmas 1.8 and 1.9 will
be alternatively those of Theorem 1, and those of Theorem 4. 


\proclaim{Lemma 1.8}
Let $X$ be $S^d$, ${\Bbb{R}}^d$ or $H^d$, and let $K,L$ and $\varphi , \psi $
be as in (*). 
Suppose {\rm{(1)}} of Theorem 1, or {\rm{(1)}} of Theorem 4. 
Suppose that both conclusions {\rm{(1)}} and {\rm{(2)}} of Lemma 1.5 hold 
(in particular, that the second sentence of Lemma 4.3 holds). 

Then any 
$x \in {\text{\rm{bd}}}\,K$ and any $y \in {\text{\rm{bd}}}\,L$ have some open
neighbourhoods
relative to ${\text{\rm{bd}}}\,K$ and to ${\text{\rm{bd}}}\,L$, 
which are congruent to relatively
open geodesic $(d-1)$-balls on a fixed 
sphere, this fixed sphere having a radius at most $\pi /2$, for $S^d$, 
on a fixed sphere for ${\Bbb{R}}^d$, and on a fixed
sphere, parasphere or hypersphere for $H^d$. (Fixed means: we have the same
sphere, parasphere or hypersphere for all $x \in {\text{\rm{bd}}}\,K$ and all 
$y \in {\text{\rm{bd}}}\,L$.) 
Moreover, the congruences carrying these relatively open neighbourhoods of
$x$ and $y$ to these relatively
open geodesic $(d-1)$-balls carry $x$ and $y$ to the
centres of these relatively open geodesic $(d-1)$-balls.
(For ${\Bbb{R}}^d$ and $H^d$ hyperplanes are excluded  --- thus for $H^d$
hyperspheres cannot degenerate to hyperplanes).
\endproclaim


\demo{Proof}
{\bf{1.}}
For the case of the proof of Theorem 1 recall \thetag{1.5.6}, \thetag{1.5.9},
\thetag{1.5.111} and \thetag{1.5.12}. 
Then, for any $2$-planes $P_{\varphi K}$ and
$P_{\psi L}$ containing $g$, we have for some $\varepsilon > 0$ that
$$
\cases
B(O_{\varphi K}, \varepsilon ) \cap \left(
{\text{bd}}\,(\varphi K) \right) \cap P_{\varphi K} {\text{ and }} 
B(O_{\psi L}, \varepsilon ) \cap \left(
{\text{bd}}\,(\psi L) \right) \cap P_{\psi L} \\
{\text{are congruent, contain }} O_{\varphi K} {\text{ and }} O_{\psi L} ,
{\text{with }} e {\text{ being an outer unit}}
\\
{\text{normal of }} \varphi K {\text{ and }} f 
{\text{ being an outer unit normal of }}
\psi L, {\text{ at }} O_{\varphi K} 
\\
{\text{and }} O_{\psi L}, {\text{ respectively.}}
\endcases
\tag 1.8.1
$$

Observe that \thetag{1.8.1} implies congruence of
$B(O_{\varphi K}, \varepsilon ) \cap 
\left( {\text{bd}}\,(\varphi K) \right) $ and 
$B(O_{\psi L}, \varepsilon ) $
\newline
$\cap \left( {\text{bd}}\,(\psi L) \right) $, both
having $g$ as an axis of rotation.

{\bf{2.}}
For the case of the proof of Theorem 4 recall \thetag{4.3.1}
and \thetag{4.3.2}.
By these we have, for some $\varepsilon > 0$, and for any
$2$-planes $P_{\varphi K}$ and $P_{\psi L}$ containing $g$, that 
$$
\cases
B( \varphi x, \varepsilon) \cap \left( {\text{bd}}\,(\varphi K) \right) \cap
P_{\varphi K} {\text{ and }} B( \psi y, \varepsilon) \cap \left( 
{\text{bd}}\,(\psi L) \right) \cap P_{\psi L} \\
{\text{are congruent, contain }} \varphi x {\text{ and }} \psi y,
{\text{ with }} g {\text{ being a normal to}} \\
\varphi K {\text{ at }} \varphi x {\text{ and to }} \psi L {\text{ at }} \psi y.
\endcases
\tag 1.8.2
$$

\newpage

Observe that \thetag{1.8.2} implies congruence of
$B(\varphi x, \varepsilon ) \cap 
\left( {\text{bd}}\,(\varphi K) \right) $ and 
$B(\psi y, \varepsilon ) \cap \left( {\text{bd}}\,(\psi L) \right) $, both
having $g$ as an axis of rotation.

{\bf{3.}}
Here the notations are different. To exclude this, we rewrite
both \thetag{1.8.1} and \thetag{1.8.2} for $K$ and $L$, by
using \thetag{1.5.6} and its analogue \thetag{4.3.4}, and Lemma 1.5, (2)
and its analogue Lemma 4.3 (about (2) of Lemma 1.5), as 
$$
\cases
B(x, \varepsilon ) \cap \left( {\text{bd}}\,K \right) {\text{ and }} B(y, 
\varepsilon ) \cap \left( {\text{bd}}\,L \right) {\text{ are congruent}}\\
{\text{surfaces of revolution, with axes of rotation their outer}}\\
{\text{unit normals }} n {\text{ and }} m 
{\text{ at }} x \in {\text{bd}}\,K {\text{ and }} y \in {\text{bd}}\,L, 
{\text{ respectively.}}
\endcases
\tag 1.8.3
$$

Yet we do not know the (congruent) shapes 
of the $2$-dimensional normal sections of 
$B(x, \varepsilon ) \cap ( {\text{bd}}\,K ) $ at $x \in {\text{bd}}\,K$
and of $B(y, \varepsilon ) \cap ( {\text{bd}}\,L ) $ at $y \in
{\text{bd}}\,L$. However, the surfaces mentioned in \thetag{1.8.3}
are surfaces of revolution about $n$
and $m$. Therefore they are also symmetric with respect to 
any hyperplane
containing their respective axes of rotation. This however implies that the
normals at any points $x^* \in B(x, \varepsilon ) \cap ( {\text{bd}}\,K ) $ 
or $y^* \in B(y, \varepsilon ) \cap ( {\text{bd}}\,L ) $ 
of a $2$-dimensional normal section of $K$ or $L$ with $2$-planes $P_K$ or
$P_L$ 
containing $n$ or $m$, respectively, lie in $P_K$ or $P_L$, respectively.
(Normals exist for Theorem 1 by hypothesis (A), and for Theorem 4 by its
hypotheses). Then the $2$-dimensional
normal sections of ${\text{bd}}\,K$ at $x$ and of ${\text{bd}}\,L$ at $y$,
containing $n$ and $m$,
are normal sections of ${\text{bd}}\,K$ at $x^*$ and of ${\text{bd}}\,L$ at
$y^*$, respectively.

Now applying  conclusion (2) of Lemma 1.5 
(the proof of Lemma 1.5 already has been completed), or its
analogue, conclusion of Lemma 4.3 about (2) of Lemma 1.5
(also already proved),
we get the following. For the above $x,x^*,y,y^*$ all sectional
curvatures of $K$ or $L$, respectively,
are equal to some $\kappa \ge 0$, where for Theorem 1 for 
${\Bbb{R}}^d$ and $H^d$ actually $\kappa > 0$, while for Theorem 4 $\kappa > 0$.
Fixing $x,y$ and varying $x^*,y^*$
we get that the $2$-dimensional normal sections
$B(x, \varepsilon ) \cap ( {\text{bd}}\,K ) \cap P_K$ and 
$B(y, \varepsilon ) \cap ({\text{bd}}\,L) \cap P_L $ 
have constant curvature $\kappa $. That is, they are relatively
open arcs of congruent circles, paracycles or hypercycles in $P_K$ and $P_L$,
with midpoints $x$ and $y$. 
For Theorem 1, for ${\Bbb R}^d$ and $H^d$ by $\kappa > 0$
they cannot be straight line 
segments, while for Theorem 4 for $S^d$
they cannot be large-circles, for ${\Bbb{R}}^d$ they cannot be straight lines, 
and for $H^d$ they cannot be straight lines, hypercycles and paracycles, by the
hypotheses of Theorem 4.

Last, we obtain $B(x, \varepsilon ) \cap ( {\text{bd}}\,K ) $ and 
$B(y, \varepsilon ) \cap ( {\text{bd}}\,L ) $ by rotation of 
$B(x, \varepsilon ) \cap  ( {\text{bd}}\,K ) \cap P_K  $ and  
$B(y, \varepsilon ) \cap ( {\text{bd}}\,L ) \cap P_L$ 
about the axes $n$ and $m$.
Therefore these sets are exactly such as stated in this lemma.
$\blacksquare $
\enddemo


\newpage

\proclaim{Lemma 1.9}
Let $X$ be $S^d$, ${\Bbb{R}}^d$ or $H^d$, and let $K,L$ and $\varphi , \psi $
be as in (*).
Suppose {\rm{(1)}} of Theorem 1, or {\rm{(1)}} of Theorem 4. 
Then the  conclusion of Lemma 1.8 implies {\rm{(2)}} of Theorem 1. 
In particular, any of {\rm{(1)}} of Theorem 1 and {\rm{(1)}} of Theorem 4
implies {\rm{(2)}} of Theorem 1.
\endproclaim


\demo{Proof}
By the conclusion of Lemma 1.8, locally, any of ${\text{bd}}\,K$ and 
${\text{bd}}\,L$ is an analytic surface (namely, sphere, parasphere or
hypersphere), given up to congruence. (I.e., we have congruent spheres,
paraspheres or hyperspheres for any points of ${\text{bd}}\,K$ and
${\text{bd}}\,L$.)

Now let $x \in {\text{bd}}\,K$ be arbitrary. By the conclusion of
Lemma 1.8, for some relatively
open geodesic $(d-1)$-ball $B_x$
on ${\text{bd}}\,K$, with centre $x$, we have that $B_x$
is a subset of an above analytic hypersurface, given up to congruence. Then
$$
\cases
{\text{for }} x_1,x_2 \in {\text{bd}}\,K, {\text{ with }} B_{x_1} \cap B_{x_2} 
\ne \emptyset , {\text{ we have that }} B_{x_1},\,\,  B_{x_2} \\
\subset {\text{bd}}\,K 
{\text{ are subsets of the same analytic 
hypersurface, i.e., they}} \\
{\text{are open subsets of the same sphere, parasphere or hypersphere.}} 
\endcases
\tag 1.9.1
$$
This follows by simple geometry (recall that in the conformal model these
surfaces are portions of spherical surfaces inside the model $S^d$), 
or by analytic continuation. 

Now let us introduce an equivalence relation $\sim $
on the points $x$ of ${\text{bd}}\,K$. 
$$
\cases
{\text{Two points }} x',x'' \in {\text{bd}}\,K {\text{ are called }}
equivalent, {\text{ written }} x' \\
\sim x'', {\text{ if there exists a finite sequence }} 
x'=x_1, \ldots , x_N = x'' \in \\
{\text{bd}}\,K, {\text{ such that }} B_{x_i} \cap B_{x_{i+1}} \ne \emptyset , 
{\text{ for each }} i=1,\ldots , N-1.
\endcases
\tag 1.9.2
$$
It is standard to show 
that $\sim $ is in fact an equivalence relation.
Let the equivalence classes with respect to $ \sim $ be denoted by
$C_{\alpha }$, the $\alpha $'s forming an index set $A$. 
(It is easy to see that $A$ is at most countably infinite, but this
is not necessary for us.) 

By the conclusion of
Lemma 1.8 and from \thetag{1.9.1}, by using induction with respect to $N$, 
we get that
$$
\cases
{\text{each equivalence class }} C_{\alpha } \subset {\text{bd}}\,K 
{\text{ is a relatively open subset}} \\
{\text{of a sphere, parasphere or hypersphere, this surface being}} \\
{\text{given up to congruence, 
and also is relatively open in }} {\text{bd}}\,K.
\endcases
\tag 1.9.3
$$ 
Two different sets $C_{\alpha } \subset {\text{bd}}\,K
$ are disjoint, since else their union
would be a subset of some equivalence class, a contradiction.
Thus $\{ C_{\alpha } \mid \alpha \in A \} $
forms a relatively open partition of ${\text{bd}}\,K$ (by \thetag{1.9.3}), 
which implies that it
forms a relatively open-and-

\newpage

closed partition of ${\text{bd}}\,K$. 
Now observe that a connected component of
${\text{bd}}\,K$ cannot intersect a relatively 
open-and-closed subset $C_{\alpha }$ of ${\text{bd}}\,K$, and also
its complement in ${\text{bd}}\,K$, which implies that
$$
{\text{each }} C_{\alpha } {\text{ (for }} \alpha \in A )
{\text { is the union of some connected components of }} {\text{bd}}\,K. 
\tag 1.9.4
$$

On the other hand, each $B_x$ is connected, and thus 
no $B_x$ can intersect different connected components of ${\text{bd}}\,K$.
Hence, by the definition of $\sim $ and by induction for $N$, we get that 
$$
\cases
{\text{the sets }} C_{\alpha } {\text{ for }} \alpha \in A
{\text{ are also subsets}} \\
{\text{of some connected components of }} {\text{bd}}\,K.
\endcases
\tag 1.9.5
$$
Now observe that both $\{ C_{\alpha } \mid \alpha \in A \} $ and the
connected components of ${\text{bd}}\,K$ form partitions of ${\text{bd}}\,K$.
Then \thetag{1.9.4} and \thetag{1.9.5} imply that 
$$
{\text{the sets }} C_{\alpha } {\text{ for }} \alpha \in A
{\text{ are exactly the connected components of }} {\text{bd}}\,K.
\tag 1.9.6
$$

Up to now, we know the following. By \thetag{1.9.6} and \thetag{1.9.3},
$$
\cases
{\text{the connected components }} C_{\alpha } {\text{ of }}
{\text{bd}}\,K {\text{ (for }} \alpha \in A ) 
{\text{ are}} \\
{\text{relatively
open subsets of some congruent spheres, para-}} \\
{\text{spheres or hyperspheres,
and also are relatively open in }} {\text{bd}}\,K . 
\endcases
\tag 1.9.7
$$
Since ${\text{bd}}\,K$ is closed
in $X$, its connected components $C_{\alpha }$, 
being relatively closed in ${\text{bd}}\,K$,
are closed in $X$ as well. Therefore 
$$
\cases
{\text{the connected components }} C_{\alpha } {\text{ of }}
{\text{bd}}\,K {\text{ are also closed}} \\
{\text{in the above congruent spheres, paraspheres or}}\\
{\text{hyperspheres containing them (from \thetag{1.9.3}).}}
\endcases
\tag 1.9.8
$$
By \thetag{1.9.7} and \thetag{1.9.8}
$$
\cases
{\text{the connected components }} C_{\alpha } {\text{ of }} 
{\text{bd}}\,K {\text{ are non-empty,}} \\
{\text{relatively open-and-closed subsets of some congruent}} \\
{\text{spheres, paraspheres or hyperspheres.}}
\endcases
\tag 1.9.9
$$

\newpage

However, spheres, paraspheres and
hyperspheres are connected, i.e., have no non-empty, 
relatively open-and-closed proper 
subsets. Therefore, taking in consideration that by the conclusion of
Lemma 1.8 the congruent
spheres, paraspheres or hyperspheres for ${\text{bd}}\,K$ are congruent to those
for ${\text{bd}}\,L$, we have that
$$
\cases
{\text{the connected components of }} {\text{bd}}\,K, {\text{ and, similarly, 
of}} \\ 
{\text{bd}}\,L, {\text{ are congruent spheres, paraspheres or hyperspheres.}}
\endcases
\tag 1.9.10
$$
This shows that (1) of Theorem 1 implies the first sentence of (2) of Theorem 1.
(For $S^d$ we have radius of the sphere
at most $\pi /2$, by hypothesis $(*)$ of Theorem 1.)

The second sentence of (2) of Theorem 1 follows for $K$ (and analogously for
$L$) like this. For
the case that in \thetag{1.9.10} we have
one sphere or one parasphere ($= {\text{bd}}\,K$), 
its convex hull $K'$ is the ball
or paraball bounded by the sphere or parasphere. Therefore $K' \subset K$.
If we had $K' \subsetneqq K$, then $k \in K \setminus K'$ and $k' \in
{\text{int}}\,K'$ would imply that for $k'' \in [k,k'] \cap 
( {\text{bd}}\,K ) $ we would have
$k'' \in {\text{int}}\,K$, a contradiction. (For $X=S^d$ we can choose $k,k'$
not antipodal.) Hence $K'= K$.

However, if in \thetag{1.9.10} we have several 
spheres or paraspheres, then by the conclusion of 
Lemma 1.8 they are necessarily disjoint, and their
closed
convex hull contains the balls and paraballs bounded by these spheres or 
paraspheres. Moreover, their closed convex hull contains a segment $[x,y]$
with $x$ and $y$ some
interior points of two different above balls or paraballs $B(x)$ and $B(y)$,
respectively, and even contains a small
neighbourhood of this segment. 
(For $X=S^d$ we may choose $x,y$ not antipodal in $S^d$.)
Then $[x,y]$ intersects the boundaries of 
$B(x)$ and $B(y)$
at points $x',y'$, with order $x,x',y',y$ on
$[x,y]$. Then $x'$ lies in the interior of the (closed)
convex hull of $B(x)$ and $B(y)$,
hence in ${\text{int}}\,K$. However,
$x'$ lies in a connected component of ${\text{bd}}\,K$, thus also in
${\text{bd}}\,K$,
a contradiction. This proves that (1) of Theorem 1 implies also
the second sentence of (2) of Theorem 1. This ends the proof of 
$(1) \Longrightarrow (2)$ in Theorem 1.
$ \blacksquare $
\enddemo


\demo{Proof of Theorem 1} 
Recall that we have \thetag{*}.
By Lemma 1.3, (2) of Theorem 1 implies (1) of Theorem 1. By Lemma 1.9, (1) of
Theorem 1 implies (2) of Theorem 1.
$\blacksquare $
\enddemo


Before the proof of Theorem 2 we need a lemma.


For a ($(d-1)$-dimensional)
spherical cap $C$ in ${\Bbb{R}}^d$ we write $S(C)$ for the sphere
containing $C$, and $B(C)$ for the ball bounded by $S(C)$. 
We write $c( \cdot )$ for the centre of a ball.
In the next lemma
we will use the conformal model for $H^d$, and we will consider this model as a
subset of the Euclidean space ${\Bbb{R}}^d$ in the usual way. 


\proclaim{Lemma 2.1}
Using the notations of Lemma 1.2 and those introduced just before this lemma,
but supposing the second possibility in Lemma 1.2 (1), 
we write $u_0$ for the unique common 
infinite point of $\varphi K_{01}$ and $\psi L_{01}$ 
(and also of $\varphi K_1$ and $\psi L_1$). 

\newpage

Further
let us suppose that the hyperplane $H$, with respect to
which $\varphi K_{01}$ and $\psi L_{01}$, and also $\varphi K_1$ and 
$\psi L_1$ are symmetric images of
each other, contains the centre $0$ of the conformal model in
$B^d\,\,( \subset {\Bbb{R}}^d)$. Further we consider everything in the
Euclidean geometry of ${\Bbb{R}}^d$. Then 
$$
\cases
B(\varphi K_1) {\text{ and }} B(\psi L_1) 
{\text{ are different congruent balls in}} \\
{\Bbb{R}}^d, {\text{and }} 
u_0 \in B(\varphi K_1) \cap B(\psi L_1) 
\subset ({\text{\rm{int}}}\,B^d) \cup \{ u_0 \} .
\endcases
$$
\endproclaim


\demo{Proof}
We may suppose that $u_0=(0, \ldots , 0, -1)$ and $(u_0 \in)\,\, 
H \subset H^d$ is the
hyperplane $\xi _1=0$.

By $0 \in H$ the symmetry with respect to the hyperplane $H$ in $H^d$ is just
the restriction of the symmetry with respect to ${\text{aff}}(H)$ in
${\Bbb{R}}^d$, i.e., with respect to
the hyperplane $\xi _1=0$ in ${\Bbb{R}}^d$.

We have that each of $\varphi K_{01}$, $\psi L_{01}$, $\varphi K_1$ and 
$\psi L_1$ are 
spherical caps in ${\Bbb{R}}^d$, which are relatively open in 
$S(\varphi K_{01})$, $S(\psi L_{01})$, $S(\varphi K_1)$ and
$S(\psi L_1)$, respectively.

The spherical caps $\varphi K_{01}$ and $\psi L_{01}$  
are symmetric images of each
other with respect to ${\text{aff}}(H)$ in ${\Bbb{R}}^d$.
Say, $\varphi K_{01}$ and $\psi L_{01}$ lie in the open halfspaces 
$\xi _1 < 0$ and
$\xi _1 > 0$, respectively.
Also $\varphi K_{01}$ and $\psi L_{01}$ (as hyperplanes in $H^d$) 
intersect $S^{d-1}$ orthogonally, and 
touch each other at $(0, \ldots , 0, -1)$. 
This implies that 
$c\left( B(\varphi K_{01}) \right) $ and $c \left( B(\psi L_{01}) \right) $ 
lie on the line parallel to the $\xi _1$-axis and
passing through $(0, \ldots , 0, -1)$, and are symmetric images of each other
with respect to ${\text{aff}}(H)$ in ${\Bbb{R}}^d$,
with $c \left( B(\varphi K_{01}) \right) $
and $c \left( B(\psi L_{01}) \right) $
lying in the open halfspaces $\xi _1 < 0$ and $\xi _1 > 0$, respectively.

The spherical caps
$\varphi K_1$ and $\psi L_1$ 
are also symmetric images of each other with respect to
${\text{aff}}(H)$ in ${\Bbb{R}}^d$. Therefore also 
$$
\cases
B(\varphi K_1) {\text{ and }} B(\psi L_1) {\text{ are the 
symmetric images of each other with}} \\ 
{\text{respect to }} 
{\text{aff}}\,(H), {\text{ hence are (different) congruent balls in }}
{\Bbb{R}}^d, 
\endcases
\tag 2.1.1
$$
proving the first statement of the lemma.

$$
\cases
{\text{The intersection }} B(\varphi K_1) \cap B(\psi L_1) {\text{ is thus 
bounded by two}} \\
{\text{congruent spherical caps }} \varphi C_K {\text{ on }} S(\varphi K_1)  
{\text{ and }} \psi C_L {\text{ on }} S(\psi L_1),\\
{\text{respectively, and these spherical caps are smaller than halfspheres.}} 
\endcases
\tag 2.1.2
$$
The common relative boundary
of $\varphi C_K$ and $\psi C_L$ (with respect to $S(\varphi K_1)$ and $S(\psi
L_1)$, respectively) is a
$(d-2)$-sphere lying in the hyperplane $\xi _1 = 0$ in ${\Bbb{R}}^d$,
and having centre 
$[c \left( B(\varphi K_1) \right) + c \left( B(\psi L_1) \right) ]/2$. 

\newpage

Observe that
${\text{relbd}}(\varphi K_{01}) = {\text{relbd}}(\varphi K_1)$ 
(taken in $S(\varphi K_{01})$ and 
$S(\varphi K_1)$, respectively), and also  
${\text{relbd}}(\psi L_{01}) = {\text{relbd}}(\psi L_1)$ (taken in 
$S(\psi L_{01})$ and 
$S(\psi L_1)$, respectively), and all these sets lie in $S^d$.
Therefore $c \left( B(\varphi K_1) \right) $ and $c \left( B(\psi L_1) 
\right) $ lie on
the open segments $\left( c \left( B(\varphi K_{01}) \right) , 0 \right) $ and  
$\left( c \left( B(\psi L_{01}) \right), 0 \right) $, 
and are symmetric images of each other with respect to the hyperplane 
$\xi _1 = 0$ in ${\Bbb{R}}^d$.
In particular, they lie in the $\xi _1 \xi _d$-coordinate-plane, 
and also in the open
slab $-1 < \xi _d <0$, 
with $c \left( B(\varphi K_{01}) \right) $
and $c \left( B(\psi L_{01}) \right) $
lying in the open halfspaces $\xi _1 < 0$ and $\xi _1 > 0$, respectively.

We consider $(0, \ldots , 0,1)$ as a vertical upward vector. 
We have that $B(\varphi K_1) \cap B(\psi L_1)$ is rotationally symmetric 
with respect
to the line $c \left( B(\varphi K_1) \right) c \left( B(\psi L_1) \right) $. 
The lowest and highest points of $B(\varphi K_1) \cap B(\psi L_1)$ 
lie in the $\xi _1 \xi _d$-coordinate plane, and 
they are the points of intersection of $S(\varphi K_1)$ and $S(\psi L_1)$ 
and the
$\xi _1 \xi _d$-coordinate plane. These are uniquely determined points: namely
$(0, \ldots , 0,-1)$ and $c \left( B(\varphi K_1) \right) + c \left( B(\psi
L_1) \right) - (0, \ldots , 0,-1)$. 
In particular, this proves 
$$
u_0 = (0, \ldots , 0, -1) \in B(\varphi K_1) \cap B(\psi L_1) ,
\tag 2.1.3
$$
proving the first half of the second statement of the lemma.

Since $c \left( B(\varphi K_1) \right) $ and 
$c \left( B(\psi L_1) \right) $ lie in the open slab $-1 < \xi _d < 0$, the 
highest
point of $B(\varphi K_1) \cap B(\psi L_1)$
lies on the $\xi _d$-axis, in the open segment $\left( (0, \ldots , 0,
-1), \right. $
\newline
$\left. (0, \ldots , 0,1) \right) $, i.e., is of
the form $(0, \ldots , 0, \beta )$, where $\beta \in (-1, 1)$.

We 
strictly increase $B(\varphi K_1) \cap B(\psi L_1)$ 
if we replace the congruent spherical
caps $\varphi C_K$ and $\psi C_L$ by 
halfspheres with the same relative boundary 
(with respect to $S(\varphi K_1)$, $S(\psi L_1)$, 
and the sphere containing these
halfspheres, respectively), lying on the same side of the hyperplane
$\xi _1=0$ as $\varphi C_K$ and $\psi C_L$, respectively. 
Thus we obtain a sphere bounding a ball $B(\varphi K_1,\psi L_1)$
with the above common relative boundaries of the spherical caps $
\varphi C_K$ and $\psi C_L$ as its equator
(Thales ball). Its lowest point is $(0, \ldots , 0, -1)$ and its highest point
is $(0, \ldots , 0, \beta )$. Then $B(\varphi K_1,\psi L_1)$ 
arises by diminishing $B^d$
from $(0, \ldots , 0, -1)$ in ratio $(1 + \beta )/2 \in (0,1)$. 
Hence 
$$
\cases
B(\varphi K_1) \cap B(\psi L_1) \subset B(\varphi K_1,\psi L_1) \subset \\
({\text{int}}\,B^d) \cup \{ (0, \ldots , 0, -1) \} = ({\text{int}}\,B^d)
\cup \{ u_0 \} ,
\endcases
\tag 2.1.4
$$
proving the second half of the second statement of the lemma.
$ \blacksquare $
\enddemo


\demo{Proof of Theorem 2}
The implication $(2) \Longrightarrow (1)$ of this Theorem
follows since $(\varphi K) \cap
(\psi L)$ has as centre of symmetry the midpoint of the segment connecting
the centres of $\varphi K$ and $\psi L$.

\newpage

Therefore we have to prove only $(1) \Longrightarrow (2)$ of this Theorem.
 
Observe that (1) of Theorem 2 implies (1) of Theorem 1, and (1) of Theorem 1
implies, by Theorem 1, (2) of Theorem 1, i.e.,
that the connected components of the boundaries both of
$K$ and $L$ are either 

(1) congruent spheres (for $X=S^d$ of radius at most $\pi /2$), or 

(2) paraspheres, or 

(3) congruent hyperspheres, 

{\noindent }and in cases (1) and (2) here $K$ and $L$ are congruent balls
(for $X=S^d$ of radius at most $\pi /2$), or they are paraballs, respectively.

In case (1) here $K$ and $L$ are congruent balls (for $X=S^d$ of radius at most
$\pi /2$), hence Theorem 2, (2) is proved.

There remained the cases here when we have $X=H^d$ and 

(2) $K$ and $L$ are two paraballs, or 

(3) the boundary components both of $K$ and $L$ are congruent
hyperspheres, and their numbers are at least $1$, but at most countably
infinite.

We are going to show that neither of these two cases can occur.

In case (2) here $K$ and $L$ are paraballs. We choose
$\varphi $ and $\psi $ so that $\varphi K = \psi L$. Then their intersection
is the paraball $\varphi K = \psi L$, 
which is not centrally symmetric, since it has exactly one point
at infinity. (This shows also the statement in brackets in (1) of Theorem
2 in this case.)
Hence case (2) here cannot occur.

In case (3) here, let all boundary
components $\varphi K_i$ of $\varphi K$, and $\psi L_i$ of $\psi L$
be congruent hyperspheres, with base 
hyperplanes $\varphi K_{0,i}$ and $\psi L_{0,i}$. 
Denote by
$\lambda $ the common value of the distance, for which these hyperspheres are
distance surfaces for their base hyperplanes. 
By the hypothesis $C^2_+$ (or its weakening (A) and (B) of the theorem)
we have $\lambda > 0$.
These base hyperplanes 
bound closed convex sets $\varphi K_0$ and $\psi L_0$, respectively, 
possibly with empty interior, 
and on the other closed side of each $\varphi K_{0i}$ as $\varphi K_i$,
and such that 
the parallel domains of $\varphi K_0$ and $\psi L_0$, with distance $\lambda $,
equal $\varphi K$ and $\psi L$, respectively, by Lemma 1.1. 

Let 
$$
\cases
H',H'' \subset H^d {\text{  be two hyperplanes, having one common infinite}} 
\\ 
{\text{point }} u_0 \,, {\text{ but no other common finite or infinite point.}} 
\endcases
\tag 2.1
$$

They are symmetric images of
each other with respect to a hyperplane $H \subset H^d$, having $u_0$
as an infinite point. As in Lemma 2.1, 
for simplicity we may
assume that $H$ contains the centre of the conformal model, which model we use
also here.
 
Then there exists isometries $\varphi , \psi $ of $H^d$ to itself such that 

\newpage

$$
\cases
\varphi K_{01} = H' {\text{ and }} \psi L_{01} = H'', {\text{ and }} 
{\text{int}}(\varphi K_0) {\text{ lies on}} \\
{\text{the opposite closed side of }} H' {\text{ as }} H'', {\text{ and }} 
{\text{int}}(\psi L_0) \\
{\text{lies on the opposite closed side of }} H'' {\text{ as }} H' . 
\endcases
\tag 2.2
$$
If one or both of these interiors
is/are empty, we consider the last two statements of \thetag{2.2}
as satisfied for the respective interior/s.

Possibly $\varphi $ (or $\psi $) is not orientation preserving. In this case
we apply after $\varphi $ (or $\psi $)
a symmetry with respect to a hyperplane orthogonally
intersecting $H'$ (or $H''$). 
Then the composed isometry satisfies the same properties
which $\varphi $ and $\psi$ in \thetag{2.2} 
satisfied, and additionally it is orientation preserving. So
we may suppose that $\varphi $ and $\psi $ are orientation
preserving.

Then the hypotheses of Lemma 1.2 are satisfied: (1) and (2) of Lemma 1.2 
by 
\thetag{2.1} and \thetag{2.2},
and (3) of Lemma 1.2 
is just a notation. Then Lemma 1.2 gives,
using its notations, that
$$
(\varphi K) \cap (\psi L) = (\varphi K_1^*) \cap (\psi L_1^*) .
\tag 2.3
$$
Here $\psi L_{01}$, $\psi L_1$ and $\psi L_1^*$ 
are symmetric images of $\varphi K_{01}$, $\varphi K_1$ 
and $\varphi K_1^*$ with respect to the hyperplane $H \subset H^d$, 
respectively.

Now we consider the conformal model as embedded in ${\Bbb{R}}^d$ in the usual
way. We will apply Lemma 2.1 together with its notations. Thus,
${\text{int}}\,( \cdot )$ and ${\text{cl}}\,( \cdot )$ 
denote interior and closure in ${\Bbb{R}}^d$, which contains 
the conformal model in the canonical way. We have 
$$
\varphi K_1^* = B(\varphi K_1) \cap ({\text{int}}\,B^d) {\text{ and }} 
\psi L_1^* = B(\psi L_1) \cap ({\text{int}}\,B^d), {\text{ implying}}
\tag 2.4
$$
$$
(\varphi K_1^*) \cap (\psi L_1^*) = B(\varphi K_1) \cap B(\psi L_1) 
\cap ({\text{int}}\,B^d) .
\tag 2.5
$$
By Lemma 2.1, \thetag{2.5} and once more by Lemma 2.1 we have
$$
\cases
B(\varphi K_1) \cap B(\psi L_1) =
[B(\varphi K_1) \cap B(\psi L_1)] \cap [({\text{int}}\,B^d) 
\cup \{ u_0 \} ]  \\
= [B(\varphi K_1) \cap B(\psi L_1) \cap ({\text{int}}\,B^d)] \cup \\
[B(\varphi K_1) \cap B(\psi L_1) \cap \{ u_0 \} ] 
= \left( (\varphi K_1^*) \cap (\psi L_1^*) \right) \cup \{ u_0 \} .  
\endcases
\tag 2.6
$$

Then, also using \thetag{2.3} and \thetag{2.6},
$$
\cases
{\text{the set of infinite points of }} (\varphi K) \cap (\psi L) = 
(\varphi K_1^*) \cap (\psi L_1^*) {\text{ is contained}} \\ 
{\text{in }} [ {\text{cl}} \left( 
(\varphi K_1^*) \cap (\psi L_1^*) \right) ] \cap S^{d-1} \subset [ {\text{cl}} 
\left( (\varphi K_1^*) \right) \cap {\text{cl}} 
\left( (\psi L_1^*) \right) ] \cap S^{d-1} \subset \\ 
[ B(\varphi K_1) \cap B(\psi L_1) ] 
\cap S^{d-1} = [\left( ( \varphi K_1^*) \cap (\psi L_1^*) \right) \cup 
\{ u_0 \} ] \cap S^{d-1} = \{ u_0 \} .
\endcases
\tag 2.7
$$

\newpage

On the other hand, by \thetag{2.1} and \thetag{2.2}
$u_0$ is an infinite point both of $H' = \varphi K_{01}$ and 
$H'' = \psi L_{01}$. Therefore, also using \thetag{2.3} (and meaning $(1
- \varepsilon ) u_0$ in ${\Bbb{R}}^d$), 
$$
\cases
{\text{for sufficiently small }} \varepsilon > 0 {\text{ we have }} 
(1 - \varepsilon) u_0 \in \\
{\text{int}}
[(\varphi K) \cap (\psi L)] = {\text{int}}[(\varphi K_1^*) \cap 
(\psi L_1^*)] , {\text{ hence }}
u_0 {\text{ is}} \\
{\text{an infinite point of }} 
(\varphi K) \cap (\psi L) = (\varphi K_1^*) \cap (\psi L_1^*) .
\endcases
\tag 2.8
$$ 
 
Then \thetag{2.7} and \thetag{2.8} imply, also using \thetag{2.3}, that
$$ 
\cases
(\varphi K) \cap (\psi L) = (\varphi K_1^*) \cap (\psi L_1^*) {\text{ has a 
unique infinite}} \\
{\text{point, namely }} u_0, {\text{ hence it is not centrally symmetric.}}
\endcases
\tag 2.9
$$
(This shows also the statement in brackets in (1) of Theorem
2 in this case as well.)
Hence here 
also case (3) of this proof cannot occur, ending the proof of Theorem 2.
$ \blacksquare $
\enddemo


Before the proof of Theorem 3 we give five simple lemmas.


\proclaim{Lemma 3.1}
Let $X=H^d$, and let $K^*$ and $L^*$ be closed convex sets, bounded by
two congruent hyperspheres $K$ and $L$, respectively.
Suppose that 
the base planes $\varphi K_0$ of $\varphi K$ and $\psi L_0$ 
of $\psi L$ either
have no common finite or infinite point, or have one common infinite point
but no other common finite or infinite point.
Suppose that $\varphi K$ lies on that
side of $\varphi K_0$ as $\psi L$, but $\psi L$ lies on the opposite side
of $\psi L_0$ as $\varphi K$. Then $\varphi K^* \subset 
{\text{\rm{int}}}\,(\psi L^*)$.
\endproclaim


\demo{Proof}
Let $\lambda > 0$ denote 
the common distance for which $K$ and $L$ are distance surfaces.

Let $\varphi x \in \varphi K^*$. 
If $\varphi x$ lies in the same (closed) side of
$\psi L_0$ as $\varphi K_0$, then $\varphi x \in {\text{int}}\,(\psi L^*)$.

If $\varphi x$ lies in the other (open) side of
$\psi L_0$ as $\varphi K_0$, then for some $\varphi x_0 \in \varphi K_0$ we have
$|(\varphi x_0)(\varphi x)| \le \lambda $. Then $[\varphi x_0, \varphi x]$ 
intersects $\psi L_0$ in some point
$\psi y \,\,( \ne \varphi x_0)$. 
Consequently, ${\text{dist}}\,(\psi L_0, \varphi x) \le 
|(\psi y)(\varphi x)| < |(\varphi x_0)( \varphi x)| \le \lambda $, hence  
again $\varphi x \in {\text{int}}\,(\psi L^*)$.
$\blacksquare $
\enddemo


\proclaim{Lemma 3.2}
Let $X=H^d$, and let $K^*$ and $L^*$ be closed convex sets, bounded by
two congruent hyperspheres $K$ and $L$, respectively.
Suppose that
$\varphi K^*$ and $\psi L^*$ have no common infinite points. Then $(\varphi
K^*) \cap (\psi L^*)$ has a centre of symmetry.
\endproclaim


\demo{Proof}
Since $\varphi K^*$ and $\psi L^*$ have
no common infinite points, therefore, by the collinear model,
the base hyperplanes $\varphi K_0$ and $\psi L_0$ of the hyperspheres $\varphi
K$ and $\psi L$ have no common finite or infinite points, moreover
$\varphi K$ lies on the side of $\varphi K_0$ where
$\psi L_0$ lies, and similarly, $\psi L$ lies on the side of $\psi L_0$ where
$\varphi K_0$ lies. Then the symmetry
with respect to the midpoint of the segment realizing the (positive)
distance of $\varphi
K_0$ and $\psi L_0$ (which exists by compactness and by $\varphi K_0$ and
$\psi L_0$ having no common infinite points, and which is orthogonal both to
$\varphi K_0$ and $\psi L_0$)

\newpage

interchanges $\varphi K_0$ and $\psi L_0$, as well as $\varphi K$ and $\psi
L$, and also $\varphi K^*$ and $\psi L^*$.
Hence it is a centre of symmetry of the set $(\varphi K^*) \cap (\psi L^*)$.
$\blacksquare $
\enddemo


\proclaim{Lemma 3.3}
Let $K \subset H^2$ be a closed convex set whose boundary has two connected
components $K_1$ and $K_2$, which are two congruent hypercycles.
\roster
\item
If the total 
number of different infinite points of $K_1$ and $K_2$ is $2$, then $K$ is a
parallel domain of a line, and the centres of symmetry of $K$ form the entire
base line for $K$.
\item
If the total number of different
infinite points of $K_1$ and $K_2$ is $3$, then $K$ has no
centre of symmetry.
\item
If the total number of different
infinite points of $K_1$ and $K_2$ is $4$, then $K$ 
has a unique centre of symmetry, namely 
the midpoint of the (unique)
segment realizing the distance of the base lines of $K_1$
and $K_2$. Moreover, the infinite points of $K_1$ and $K_2$ do not separate
each other on the boundary $S^1$ of the model circle (conformal or collinear).
\endroster
\endproclaim


\demo{Proof}
(1) We need to show only the statement about the centres of symmetry.

The points of the base line
are centres of symmetry of the base line, hence also
of the parallel domain of the base line. 

On the other hand, through any point $k$ of $K$,
not on the base line, we can draw a straight line $l$ 
orthogonal to the base line. Then
$k$ divides the chord of the parallel domain of the base line, lying
on $l$, into two segments, one shorter than 
the distance $\lambda > 0$ for which the hypercycle is a distance line, and
one longer than this distance $\lambda $. Therefore $k$ is not a centre of
symmetry of the parallel domain of the base line.

(2) Suppose that $K_1$ and $K_2$ have one common infinite point, but their
other infinite points are different. Then any symmetry of $K$ preserves
$K_1 \cup K_2 = {\text{bd}}\,K$, 
hence also the set of all different infinite points of
$K_1$ and $K_2$. Then for $K$ centrally symmetric 
the total number of different infinite points
of $K_1$ and $K_2$ has to be even, a contradiction.

(3) We have that $K_1$ and $K_2$, as well as their base lines,
are interchanged by the central symmetry with respect
to the mid-point $m$ of the segment realizing the distance of the base lines
(cf. the proof of Lemma 3.2). This segment is unique, cf. \cite{AVS}, Ch. 1,
Theorem 4.2, and Ch. 4, 1.7. Hence this segment, as well as its mid-point 
$m$ are invariant under any symmetry of $K$. Hence a central
symmetry of $K$ has symmetry centre $O$, say, which is the midpoint of 
the segment with end-points $m$ and the image of $m$
under this central symmetry, which is $m$ (by the last sentence). 
That is, $O$ is the midpoint of the degenerate
segment $[m,m]$, i.e., $O = m$.

\newpage

Then the straight line through $m$, orthogonal to 
the segment realizing the distance of the base lines of $K_1$ and $K_2$,
strictly separates these base lines, and
also their infinite points. Therefore the infinite points of $K_1$ and of
$K_2$ cannot separate each other on the boundary $S^1$ of the model circle.
$\blacksquare $
\enddemo


\proclaim{Lemma 3.4}
Let $P = (\varphi K) \cap (\psi L)$ be a compact convex
hypercycle-arc polygon with the hypercycles containing its arc-sides being
congruent (and $P$ has non-empty interior, has
finitely many arc-sides, and has angles in $(0, \pi )$). Let the arc-sides 
of $P$ lie alternately on ${\text{\rm{bd}}}\,(\varphi K)$ and on 
${\text{\rm{bd}}}\,(\psi L)$. Suppose that 
${\text{\rm{bd}}}\,P$ does not
consist of two finite hypercycle arcs. Let $s_1$ and $s_2$ be two
neighbouring arc-sides 
of $P$, following each other in the positive sense. Then the
total number of the infinite points of the hypercycles $H_1$  
and $H_2$ containing the 
hypercycle-arc-sides $s_1$  and $s_2$ is $4$, and the infinite points of $H_1$
and $H_2$ separate each other on the boundary $S^1$ of the model
circle (conformal, or collinear).

More exactly, let us orient ${\text{\rm{bd}}}\,P$ in the positive sense, and 
let us orient
$H_1$ and $H_2$ coherently with the orientations of $s_1$ and $s_2$. Let us
denote by $h_{11},h_{12} \in S^1$ (or $h_{21},h_{22} \in S^1$) 
the first and last infinite points of $H_1$ (or $H_2$) on $S^1$. Then these
points have the following cyclic order on the positively oriented $S^1$:
$h_{11},h_{21},h_{12},h_{22}$.
\endproclaim


\demo{Proof}
Let, e.g., $s_1$ lie on ${\text{bd}}\,(\varphi K)$ and $s_2$ lie on
${\text{bd}}\,(\psi L)$.
Let us suppose that the common vertex of the arc-sides $s_1$ and $s_2$ is the
centre of the model circle, and that $H_1 \setminus \{ 0 \} $ lies in
the open 
upper half of the model circle (thus has as tangent at $0$ the horizontal
axis). (Observe that the statement of the Lemma is invariant under the choice
of the centre of the model.)
Then $H_2$ is obtained from $H_1$ by a rotation about $0$, through some
angle $\beta \in (0, \pi )$, in the positive sense. 
This rotation has centre $0$, therefore it is a rotation
in the Euclidean sense as well. In particular, $H_1$ and $H_2$ are congruent
in the Euclidean sense as well.
We orient the base lines of $H_1$ (and of $H_2$) from $h_{11}$ to $h_{12}$ (and
from $h_{21}$ to $h_{22}$).

Let the central angle at $0$ of the base lines of 
both hypercycles be $\alpha \in (0, \pi )$ (for this observe that 
$H_1 \setminus \{ 0 \} $ lies in the open upper half of the model circle). 

Then for $\beta \in
(0, \alpha )$ we have that  $h_{11},h_{12}$ and $h_{21},h_{22}$ 
are all different and separate each other on $S^1$, and follow each other 
in the positive cyclic order on $S^1$, as asserted by the lemma. 

For $\beta = \alpha $ and for $\beta \in (\alpha ,\pi )$ we obtain by Lemma
1.2 (applied with $H_1$ and $H_2$ as $\varphi K_1$ and $\psi L_1$ in Lemma 1.2)
that $P = (\varphi K) \cap (\psi L)$ 
equals the intersection of two closed convex sets, bounded by $H_1$ and
$H_2$, respectively. Therefore ${\text{bd}}\,P$ consists of two
hypercycle-arcs, one on $H_1$ and the other on $H_2$. 

Therefore,
for $\beta = \alpha $ we have that $P$ 
is bounded by two semi-infinite arcs on $H_1$ and $H_2$, hence is not
compact, contrary to the hypothesis of the lemma.

\newpage

For $\beta \in (\alpha ,\pi )$ we have that $P$ 
is bounded by two finite arcs on $H_1$ and $H_2$, again 
contrary to the hypothesis of the lemma. This proves the lemma.
$\blacksquare $
\enddemo


\proclaim{Lemma 3.5}
Let $P_0 = (\varphi _0 K) \cap (\psi _0 L)$ satisfy all hypotheses of Lemma
3.4 (written there 
for $P = (\varphi K) \cap (\psi L)$), with each vertex lying only on the
boundary component hypercycles of ${\text{\rm{bd}}}\,( \varphi _0 K)$ and of 
${\text{\rm{bd}}}\,(\psi _0 L)$ 
which contain the arc-sides incident to the vertex.
Moreover, let $P_0$ have a centre of
symmetry $O_0$. Then $O_0$ is uniquely determined. 

Moreover, for all sufficiently small
perturbations $\varphi $ of $\varphi _0$ and $\psi $ of $\psi _0$, satisfying
that $P = (\varphi K) \cap (\psi L)$  has a centre of symmetry $O$, we have the
following. Any pair of opposite arc-sides 
of $P_0$ (i.e., images of each other by
the central symmetry with respect to $O_0$) remains by the small perturbation
a pair of opposite arc-sides
of $P$ (i.e., are the images of each other by
the central symmetry with respect to $O$), the arc-sides of $P_0$ and $P$
being identified via the small perturbation.
\endproclaim


\demo{Proof}
{\bf{1.}} 
First we prove the first statement. Even, we prove that
any compact set $\emptyset \ne C \subset H^d$ has at most one centre of
symmetry. In fact, we can copy the proof for ${\Bbb{R}}^d$. We may suppose that
$C$ has at least two points, else the statement is immediate (observe that now
$X = H^d$, and we can have two centres of symmetry only for $X = S^d$). Then
consider a ball $B \subset H^d$ of minimal (positive) 
radius, containing $C$ (existing by a compactness argument). 
Then $B$ is uniquely determined. In
fact, if there existed two such balls $B_1$ and $B_2$, then their
intersection also would contain $C$, and this intersection would be 
contained in the Thales ball with equator $({\text{bd}}\,B_1) \cap
({\text{bd}}\,B_2)$, which has a smaller radius than those of $B_1$ and $B_2$. 

Then a centre of symmetry of
$C$ coincides with the centre of this unique ball $B$. That is, $O_0$ is
uniquely determined, and hence the first statement of the Lemma is proved.
 
{\bf{2.}}
We may suppose that also $P$ is compact. Then the topological type of $P$
is the same as that of $P_0$ (identifying the arc-sides of $P_0$ and $P$ via the
small perturbations), including also which arc-side lies on
${\text{bd}}\,(\varphi K)$ (on ${\text{bd}}\,(\varphi _0 K)$)
and which arc-side lies on ${\text{bd}}\,(\psi L)$ 
(on ${\text{bd}}\,(\psi _0 L)$).

We already know that $O_0$ and $O$ are uniquely determined.

Let us suppose the contrary of the second statement of the Lemma. 
I.e., there are arbitrarily small perturbations
$\varphi $ of $\varphi _0$ and $\psi $ of $\psi _0$, such that the
``oppositeness relation'' for the arc-sides of $P$ is not the same, as the
oppositeness for the arc-sides of $P_0$, when we identify the arc-sides 
of $P_0$ and
of $P$ via the small perturbations. Observe that the oppositeness relation is
a cyclic perturbation of the arc-sides, and there are exactly $n$ such cyclic
perturbations, where $n$ is the number of arc-sides of $P_0$ (and of $P$). 
Therefore we may suppose that among the arbitrarily small

\newpage

perturbations we consider only such ones, for which these cyclic permutations
are a fixed cyclic permutation, which is different from the cyclic
permutation given by the oppositeness relation for $P_0$ (and for $O_0$). 

Then choosing a suitable subsequence of these small perturbations, we obtain in
its limit the same arc-polygon $P_0$, but with another oppositeness relation,
than that via $O_0$. That is, for some arc-side $s_0$ of $P_0$ its opposite
arc-side with respect to $O_0$ 
is $s_0'$, and in this limit situation the opposite arc-side of
$s_0$ is some other arc-side 
$s_0'' \,\,(\ne s_0')$ of $P_0$. Then the centres of
symmetry cannot be the same, and thus $P_0$ has two different centres of
symmetry. This however contradicts the first statement of this lemma,
and hence the second statement of the Lemma is proved.
$\blacksquare $
\enddemo


\demo{Proof of Theorem 3}
{\bf{1.}}
{\it{Let $X$ be $S^d$ or ${\Bbb{R}}^d$.}} 
Then by the already proved Theorems 1 and
2, we have
$$
\cases
{\text{Theorem }} 3,\,\, (1) \Longrightarrow {\text{ Theorem }} 1, (1) 
\Longrightarrow {\text{ Theorem }} 1,\,\, (2) \\
\Longleftrightarrow {\text{ Theorem }} 3, (2)\,\, (a) \,\,(\Longleftrightarrow 
{\text{ Theorem }} 3, (2)) \Longleftrightarrow \\
{\text{ Theorem }} 2,\,\, (2) \Longrightarrow {\text{ Theorem }} 2,\,\, (1) 
\Longrightarrow {\text{ Theorem }} 3,\,\, (1).
\endcases
\tag 3.1
$$
In particular,
$$
{\text{For }} X = S^d, \,\,{\Bbb{R}}^d {\text{ we have that Theorem 3,}} \,\,
(1) \Longleftrightarrow {\text{ Theorem }} 3, \,\, (2).
\tag 3.2
$$

{\bf{2.}}
{\it{There remained the case $X=H^d$.}} 

{\it{First we prove Theorem}} 3, (2) $\Longrightarrow $ {\it{Theorem}} 3, (1).

By Theorem 2 we have 
$$
\cases
{\text{Theorem }} 3, (2)\,\, (a) \Longleftrightarrow {\text{ Theorem }} 2, (2)
\Longrightarrow \\
{\text{Theorem }} 2, (1) \Longrightarrow {\text{ Theorem }} 3, (1). 
\endcases
\tag 3.3
$$

We have
$$
{\text{Theorem }} 3,\,\, (2)\,\, (b) \Longrightarrow {\text{ Theorem }} 3, (1)
\tag 3.4
$$
by the proof of the implication Theorem 1, (2) $ \Longrightarrow $ Theorem 1,
(1), in Lemma 1.3.  

For the following recall that by Theorem 1 we have
$$
{\text{Theorem 3, }} (1) \Longrightarrow {\text{ Theorem }} 1, (1) 
\Longrightarrow {\text{ Theorem }} 1, (2).
\tag 3.5
$$
Recall that the case of congruent balls, or of two paraballs already were 
settled above,
at Theorem 3, (2), (a) or (b) $\Longrightarrow $ Theorem 3, (1). 

\newpage

{\it{There remained to prove}} 
$$
{\text{Theorem 3, (2) (c) }} \Longrightarrow {\text{ Theorem 3, (1).}}
\tag 3.6
$$
Observe that the
connected components of ${\text{bd}}\,K$ and of ${\text{bd}}\,L$ are disjoint.
Therefore at proving \thetag{3.6}, by Theorem 3, (2), (c) we may suppose that 
$$
\cases
{\text{the connected components of the boundaries of}} \\
{\text{both }} K {\text{ and }} L {\text{ are disjoint congruent 
hyperspheres}} \\
{\text{(degeneration to hyperplanes being not admitted).}}
\endcases
\tag 3.7
$$
We are going to 
prove \thetag{3.6} in each of the cases listed in Theorem 3, (2) (c).

{\bf{2.1.}}
First we prove that 
$$
{\text{Theorem }} 3,\,\, (2)\,\, (c), {\text{ and }} d \ge 3 \,\, 
\Longrightarrow {\text{ Theorem }} 3,\,\, (1)
\tag 3.8
$$
(i.e., case (2), (c) $(\alpha )$ in Theorem 3).

By \thetag{3.7}, the infinite points of all connected components of
${\text{bd}}\,K$ or of ${\text{bd}}\,L$ are sub-$(d-2)$-spheres of the
boundary of the model (either conformal, or collinear). They bound {\it{open
spherical caps on the boundary of the model,}} called {\it{associated to $K$ 
and $L$, such that the convex hulls of these open spherical
caps}} (meant in ${\Bbb{R}}^d$, containing the collinear model of $H^d$ in the
canonical way) {\it{contain the 
respective connected component of \,${\text{\rm{bd}}}\,K$ or of \,
${\text{\rm{bd}}}\,L$.}}

We will show that these open spherical caps are disjoint (but may have
common boundary points). We use the notation $K_0$ from Lemma 1.1 (and
analogously we use the notation $L_0$). 
In the collinear model, $K_0$ or $L_0$ can be obtained
from the model (open) unit ball by cutting off the interiors (in
${\Bbb{R}}^d$) of the
convex hulls of these open spherical caps. (Thus we obtain a set
``like a polytope'', with possibly infinitely many facets, and with other
boundary points on the boundary of the model.)
This implies that, for any of $K$ and $L$, 
no such open spherical cap can contain another such open spherical cap (recall
Lemma 3.1), and 
also that no two such open spherical caps can have a partial overlap (else
${\text{bd}}\,K_0$, meant in ${\Bbb{R}}^d$, 
would have points on both sides of an above ``facet'' of
it --- 
which is the convex hull of the infinite points of an above $(d-2)$-sphere, in
the collinear model).

Let us consider the infinite points of $K_0$ (or of $L_0$), denoted by
$({\text{cl}}\,K_0) \cap S^{d-1}$ (or by $({\text{cl}}\,L_0) \cap S^{d-1}$), 
{\it{where we mean closure}} ${\text{cl}}$ {\it{in}} ${\Bbb{R}}^d$ and
$S^{d-1}$ is the boundary of the model. These can be
obtained from $S^{d-1}$, by deleting all above disjoint open
spherical caps, associated to $K$ (or to $L$).

We use on the boundary $S^{d-1}$ of the model (conformal or collinear) the
geodesic metric inherited from its superset ${\Bbb{R}}^d$.

We are going to show that 
$$
({\text{cl}}\,K_0) \cap S^{d-1} {\text{ (and also }} 
({\text{cl}}\,L_0) \cap S^{d-1}) {\text{ is connected.}} 
\tag 3.9
$$

\newpage

In fact, any two of the points of the set(s) in \thetag{3.9}
can be connected by a geodesic segment $S$ on $S^{d-1}$. This segment $S$
may have relatively open subsegments $S \cap C$ lying in some
above open spherical caps $C$, associated to $K$ (or to $L$), but 
$$
\cases
{\text{then these subsegments }} S \cap C {\text{ will be replaced by the}} \\
{\text{shorter (or some equal) geodesic segments }} S(C) {\text{ on the}} \\ 
{\text{relative boundaries of these spherical caps }} C, {\text{ where}} \\
{\text{the endpoints of }} S \cap C {\text{ and those of }} S(C) 
{\text{ coincide.}}
\endcases
\tag 3.10
$$
(Observe that now $d-2 \ge 1$, therefore these sub-$(d-2)$-spheres are
connected.) Doing this
simultaneously for all these relatively open subsegments, we claim that
$$
\cases
{\text{we obtain a continuous path connecting the arbitrarily chosen}} \\
{\text{points of }} ({\text{cl}}\,K_0) \cap S^{d-1} {\text{ (and also of }} 
({\text{cl}}\,L_0) \cap S^{d-1}), {\text{ in}} \\
({\text{cl}}\,K_0) \cap S^{d-1} {\text{ (in }} ({\text{cl}}\,L_0) \cap
S^{d-1}), {\text{ proving arcwise}} \\ 
{\text{connectedness of }} ({\text{cl}}\,K_0) \cap S^{d-1} {\text{ (and also 
of }} ({\text{cl}}\,L_0) 
\cap S^{d-1}).
\endcases
\tag 3.11
$$

Now we are going to prove \thetag{3.11}. Actually, this
``perturbation'' (via \thetag{3.10}) of the original geodesic segment $S$
is a continuous image of $S$, which suffices to prove \thetag{3.11}.
We define the continuous function $f$ from $S$ to the path described
in \thetag{3.10} as follows. The function $f$ maps points of $S$ in 
some above
open spherical cap $C$ (i.e., points of some $S \cap C$)
to the smaller (or some equal) geodesic segment $S(C)$ on the
boundary of the spherical cap $C$, connecting the two endpoints of $S \cap C$,
in the following way. 
If a point in $S \cap C$ moves
with constant velocity between the two endpoints of $S \cap C$, then its image
in $S(C)$
moves with constant velocity between the same
two endpoints of the geodesic segment
$S(C)$ (i.e., between the two endpoints of $S \cap C$). All other points of $S$
are mapped by $f$ to themselves. Evidently this function $f$ has a
Lipschitz constant at most $ \pi /2$, hence is in fact continuous.
{\it{This proves our claim \thetag{3.11}.}}

Now let $(\varphi K) \cap (\psi L)$ be compact. Then 
$$
\cases
\varphi K {\text{ and }} \psi L {\text{ cannot have any common infinite point
--- }} \\
{\text{equivalently, }} \varphi K_0 {\text{ and }} \psi L_0 {\text{ cannot have 
any common}} \\
{\text{infinite point. (In fact, this holds for $d \ge 2$.)}}
\endcases
\tag 3.12
$$
In fact, else, using the collinear model, and ${\text{int}} \left( ( \varphi
K) \cap (\psi L) \right)  \ne \emptyset $, we obtain a contradiction to
compactness of $ ( \varphi K) \cap (\psi L)$. 

Therefore $\left( {\text{cl}}\,(\varphi K_0) \right) \cap S^{d-1}$ 
and $\left( {\text{cl}}\,(\psi L_0) \right) \cap S^{d-1}$ are 
disjoint. Then $\left( {\text{cl}}\,(\varphi K_0) \right) \cap S^{d-1}$ 
is, by \thetag{3.11}, a connected subset of $S^{d-1} \setminus 
\left( {\text{cl}}\,(\psi L_0) \right) $, therefore 

\newpage

$$
\cases
\left( {\text{cl}}(\varphi K_0) \right)
\cap S^{d-1} {\text{ is contained in a connected component}} \\
{\text{of }} S^{d-1} \setminus 
\left( {\text{cl}}\,(\psi L_0) \right), 
{\text{ which component is the image by}} \\
\psi {\text{ of some of the open spherical caps associated to }} L.
\endcases
\tag 3.13
$$
Now we change the roles of $\varphi K$ and $\psi L$. Therefore also
$$
\cases
\left( {\text{cl}}\,(\psi L_0) \right)
\cap S^{d-1} {\text{ is contained in a connected component}} \\
{\text{of }} S^{d-1} \setminus \left( {\text{cl}}\,(\varphi K_0) \right), 
{\text{ which component is the image by}} \\
\varphi {\text{ of some of the open spherical caps associated to }} K.
\endcases
\tag 3.14
$$

By \thetag{3.13} and \thetag{3.14}
we are in the situation of Lemma 1.2. Therefore we have, with the
notations of Lemma 1.2, that 
$$
(\varphi K) \cap (\psi L) = (\varphi K_1^*) \cap (\psi L_1^*).
\tag 3.15
$$
Then $\varphi K_1^*$ and $\psi L_1^*$ have no common infinite points,
by \thetag{3.13} and \thetag{3.14}. Then 
by Lemma 3.2 the set in \thetag{3.15} has a centre of symmetry, {\it{and
thus \thetag{3.8} is proved.}} 

{\bf{2.2.}}
Second we prove that 
$$
\cases
{\text{Theorem }} 3,\,\, (2)\,\, (c), {\text{ and }} d = 2  {\text{ and 
one of }} K {\text{ and }} L \\
{\text{is bounded by one hypercycle }} \Longrightarrow {\text{ Theorem }} 
3,\,\, (1)
\endcases
\tag 3.16
$$
(i.e., case (2) (c) $( \beta )$ $(\beta ')$ in Theorem 3).

Suppose, e.g., that ${\text{bd}}\,K$
has only one connected component $K_1$. Then, with the notations of Lemma 1.2,
we have  $K_1^* = K$. Then 
$\left( {\text{cl}}\,(\varphi K_1^*) \right) \cap S^1 = 
\left( {\text{cl}}\,(\varphi K) \right) \cap S^1$ 
is a closed 
subarc of the boundary $S^1$ of the model, of length in $(0,
2 \pi )$. Therefore, by \thetag{3.12}, 
$\left( {\text{cl}}\,(\psi L) \right) \cap S_1$ 
is contained in the complementary open subarc, of length in $(0,
2 \pi )$,
which is the complement of the above closed subarc in $S^1$. Then
again we are in the situation of Lemma 1.2, with $K_1$ from above, and with a
suitable connected component $\psi L_1$ of ${\text{bd}}\,(\psi L)$. 
Then by Lemma 1.2 and using its notations we have that \thetag{3.15} holds
once more. Then, as in the end of {\bf{2.1}}, by Lemma 3.2 
the set in \thetag{3.15} has a centre of symmetry, {\it{and
thus \thetag{3.16} is proved}}. 

{\bf{2.3.}}
Third we prove that 
$$
\cases
{\text{Theorem }} 3,\,\, (2)\,\, (c), {\text{ and }} d = 2  {\text{ and }}
K {\text{ and }} L {\text{ are}} \\
{\text{congruent parallel domains of straight lines }} K_0 \\
{\text{and }} L_0, {\text{respectively }} \Longrightarrow {\text{ Theorem }} 
3,\,\, (1)
\endcases
\tag 3.17
$$

\newpage

(i.e., case (2) (c) $(\beta )$ $(\beta '')$ in Theorem 3).

Suppose that 
$\varphi K_0$ and $\psi L_0$ have a common finite point. Any of these finite
points is a centre of symmetry both of $\varphi K$ and $\psi L$ (cf. Lemma
3.3, Proof, (1)), hence also of $( \varphi K) \cap (\psi L)$. 

If $\varphi K_0$ and $\psi L_0$ have no
common finite point, but have a common infinite point, then
we obtain a contradiction to \thetag{3.12}.

Let $\varphi K_0$ and $\psi L_0$
have no common finite or infinite point. Then 
$$
\cases
(\varphi K) \cap (\psi L) {\text{  is the intersection of four closed convex 
sets}} \\
\varphi K_1^*, \,\,\varphi K_2^* {\text{ and }} \psi L_1^*, \,\, \psi L_2^*, 
{\text{ bounded by the hypercycles}} \\
\varphi K_1, \,\,\varphi K_2 {\text{ and }} \psi L_1,\,\, \psi L_2, 
{\text{ which are the connected}} \\
{\text{components of }}
{\text{bd}}\,(\varphi K) {\text{ and of }} {\text{bd}}\,(\psi L), 
{\text{ respectively.}} 
\endcases
\tag 3.18
$$
Let $\varphi K_1$ (and $\psi L_1$)
lie on that side of $\varphi K_0$ (and $\psi L_0$) where $\psi L_0$ (and 
$\varphi K_0$) lies, and then $\varphi K_2$ 
(and $\psi L_2$) lies on the other side of $\varphi K_0$ (and $\psi L_0$). 
Then Lemma 3.1 implies
$$
\cases
(\varphi K) \cap (\psi L) = \left( (\varphi K_1^*) \cap (\varphi
K_2^*) \right) \cap \left( (\psi L_1^*) \cap (\psi L_2^*) \right) = \\
\left( (\varphi K_1^*) \cap (\psi L_2^*) \right) \cap \left( (\psi L_1^*) \cap 
(\varphi K_2^*) \right) = (\varphi K_1^*) \cap (\psi L_1^*).
\endcases
\tag 3.19
$$
The (last) set in \thetag{3.19} has by Lemma 3.2 a centre of symmetry, {\it{and
thus \thetag{3.17} is proved}}.

{\bf{2.4.}}
Fourth we prove that
$$
\cases
{\text{Theorem 3, (2) (c), and }} d = 2 {\text{ and there are no more}} \\
{\text{compact intersections }} (\varphi K) \cap (\psi L) {\text{ than 
those bounded}} \\
{\text{by two finite hypercycle arcs}} \Longrightarrow 
{\text{ Theorem }} 3, 
(1).
\endcases
\tag 3.20
$$
(i.e., case (2) (c) $(\beta )$ $(\beta ''')$ in Theorem 3).
 
If $(\varphi K) \cap (\psi L)$ is bounded by two finite
hypercycle arcs, then it has
a centre of symmetry by Lemma 3.2.
Else we have Theorem 3, (1)
vacuously. {\it{Thus \thetag{3.20} is proved.}}

{\bf{3.}}
Summing up: we have shown Theorem 3, (1) $\Longleftrightarrow $ 
Theorem 3, (2) for $X = S^d, {\Bbb{R}}^d$ (cf. \thetag{3.2}). 
Further, we have shown Theorem 3, 
(2) $\Longrightarrow $ Theorem 3, (1), in all cases: for $X = H^d$:
for case (a) (cf. (\thetag{3.3}); for case (b) (cf. \thetag{3.4}); for case (c) 
$(\alpha )$ (cf. \thetag{3.8}, and last sentence of {\bf{2.1}}); for case (c)
$(\beta )$ $(\beta ' )$ (cf. \thetag{3.16}, and last sentence of {\bf{2.2}});
for case (c) $(\beta )$ $(\beta '')$ (cf. \thetag{3.17}, and last sentence of
{\bf{2.3}}); for case (c) $(\beta )$ $(\beta ''')$ (cf. \thetag{3.20}, and
last sentence of {\bf{2.4}}).

\newpage

There remained the case $X=H^d$, and still we have to prove (1) 
$\Longrightarrow $ (2). We will show the equivalent implication
$\lnot (2) \Longrightarrow \lnot (1)$.
In other words, if neither of Theorem 3, (2), (a), (b), (c) ($\alpha $),
($\beta $) holds, then we have the negation of Theorem 3 (1). In formula
(taking in account the last sentence of Theorem 1, (2) and \thetag{3.7}):
$$
\cases
{\text{we have }} X = H^d,\,\,d=2, {\text{ and the connected components of}} \\
{\text{the boundaries of both }} K {\text{ and }} L {\text{ are congruent 
hypercycles}} \\
{\text{(degeneration to straight lines being not admitted), both }} 
{\text{bd}}\,K \\
{\text{and }} {\text{bd}}\,L {\text{ have at least two connected 
components, and at}} \\
{\text{least one of them is not a parallel domain of a straight line}} \\
{\text{(observe that at the negation of Theorem 3, (c), }} (\beta ), (\beta
'') {\text{ we}} \\
{\text{cannot have incongruent parallel domains of straight lines, by}} \\
{\text{congruence of the boundary components of }} K {\text{ and }} L),
{\text{ and}} \\
{\text{there exists a compact intersection }} (\varphi K) \cap 
(\psi L) {\text{ not bounded}} \\
{\text{by two finite hypercycle arcs }} \Longrightarrow \lnot 
{\text{ Theorem 3, (1)}}. 
\endcases
\tag 3.21
$$
That is, under the hypotheses of \thetag{3.21} 
we have to give $\varphi , \psi $ so that 
$$
\cases
(\varphi K) \cap (\psi L) {\text{ is compact, but is not centrally symmetric
(in particular,}} \\
{\text{is not bounded by two finite hypercycle arcs, cf. Lemma 3.2).}} 
\endcases
\tag 3.22
$$

We begin with choosing a compact intersection $P_0 := 
(\varphi _0 K) \cap (\psi _0
L)$ not bounded by two finite hypercycle arcs. We may suppose that 
$$
\cases
P_0 = (\varphi _0 K) \cap (\psi _0 L) {\text{ (which is compact and is not
bounded}} \\
{\text{by two finite hypercycle arcs) is centrally symmetric,}} 
\endcases
\tag 3.23
$$
else we are done.

First we observe that 
$$
\cases
{\text{for }} (\varphi K) \cap ( \psi L) {\text{ compact, only finitely
many connected}} \\
{\text{components of }} 
\varphi K {\text{ and of }} \psi L {\text{ can contribute to }} {\text{bd}} 
\left( (\varphi K) \cap ( \psi L) \right) \\ 
{\text{(even, }} {\text{bd}} \left( (\varphi K) \cap ( \psi L) \right) 
{\text{ has an empty intersection with all other}} \\
{\text{connected components of }} {\text{bd}}\,(\varphi K) {\text{ and of }} 
{\text{bd}}\,(\psi L) ) . 
\endcases
\tag 3.24
$$
In fact, using the
collinear model, suppose that $(\varphi K) \cap ( \psi L)$ lies in a closed 
circle 
of radius $1 - \varepsilon $ about the centre $0$ of the model, where
$\varepsilon $ is small. Then the 

\newpage

base
lines of those hypercycles, which are 
connected components of ${\text{bd}}\,(\varphi K)$ (and
also of ${\text{bd}}\,(\psi L)$) whose some points are on 
${\text{bd}} \left( (\varphi K) \cap ( \psi L) \right) $, are disjoint, hence
span disjoint open
angular domains with vertex at the centre $0$ of the model, with spanned
angles at least $2 \arccos (1 - \varepsilon )$.  (If $0$ is strictly 
separated by one of the base lines of the hypercycle boundary components 
of $\varphi K$ (or of $\psi L$) 
from the base lines of the remaining hypercycle boundary components of
$\varphi K$ (or of $\psi L$), 
then the corresponding angle is considered as greater than $\pi $.)
Hence the total number of the hypercycle boundary components, 
taken together for $\varphi K$ and $\psi L$,
whose some points are on ${\text{bd}} \left( (\varphi K) \cap (\psi L) \right)
$, is at most $2 \cdot 2 \pi
/ \left( 2 \arccos (1 - \varepsilon ) \right) $, hence is finite.
{\it{This proves}} \thetag{3.24}.

Therefore we may say that 
$$
(\varphi K) \cap ( \psi L) {\text{ is a {\it{hypercycle-arc-polygon}} (later 
{\it{arc-polygon}}).}}
\tag 3.25
$$
$$
\cases
(\varphi K) \cap ( \psi L) {\text{ cannot have two neighbourly sides,
belonging to}} \\
{\text{different connected components of either }} {\text{bd}}\,(\varphi K) 
{\text{ or }} {\text{bd}}\,\psi L \\
{\text{(these components being disjoint), neither belonging to the}} \\
{\text{same component of either }} {\text{bd}}\,(\varphi K) {\text{ or }} 
\psi L  {\text{ (else the union of}} \\
{\text{these sides would form a single side of }} (\varphi K) \cap ( \psi L) 
). {\text{ By the}} \\
{\text{same reason, it is impossible that two different neighbourly}} \\ 
{\text{arc-sides would lie on the same hypercycle, which is a common}} \\
{\text{boundary component of }} \varphi K {\text{ and }} \psi L . {\text{ (This 
excludes angles  $\pi $.)}} \\
{\text{Therefore }} (\varphi K) \cap ( \psi L) {\text{ has alternately 
arc-sides on hypercycles}} \\ 
{\text{(only) in }} {\text{bd}}\,(\varphi K) {\text{ and (only) in }}
{\text{bd}}\,(\psi L) {\text{ --- in particular, it has}} \\
{\text{an even number of arc-sides. Also by the same reason, through}} \\
{\text{any vertex of our arc-polygon there cannot pass a third}} \\
{\text{boundary component either of }} \varphi K {\text{ or of }} \psi L. 
\endcases
\tag 3.26
$$

Therefore (cf. \thetag{3.23}, and applying \thetag{3.25}, \thetag{3.26} 
for $\varphi _0$ and $\psi _0$ rather than $\varphi $ and $\psi $), 
$$
\cases
P_0 = (\varphi _0 K) \cap (\psi _0 L) {\text{ is a centrally symmetric
hypercycle-arc-polygon,}} \\
{\text{whose arc-sides lie alternately (only) on }} 
{\text{bd}}\,(\varphi _0 K) {\text{ and (only) on}} \\
{\text{bd}}\,(\psi _0 L), {\text{ and whose centre of symmetry will be denoted 
by }} O_0 . 
\endcases
\tag 3.27
$$
Since $P_0$ is not bounded by two finite hypercycle arcs (cf. \thetag{3.23}),
by \thetag{3.26} 

\newpage

$$
{\text{ the (even) number }} n {\text{ of the arc-sides of }} P_0 {\text{ is 
at least $4$.}}
\tag 3.28
$$
We are going to show that 
$$
\cases
{\text{for some small perturbations }} \varphi  {\text{ and }} \psi 
{\text{ of the original }} \varphi _0 {\text{ and}} \\
\psi _0, {\text{we have that }} P := (\varphi K) \cap ( \psi L) {\text{ is not 
centrally symmetric.}} 
\endcases
\tag 3.29
$$
Observe that since no vertex of $P_0$ lies on any boundary component either
of $\varphi _0 K$ or $\psi _0 L$ other than the (only)
boundary components containing the arc-sides at this vertex, therefore 
$$
\cases
{\text{by small perturbations }} \varphi  {\text{ of }} \varphi _0 
{\text{ and }} \psi  {\text{ of }} \psi _0, {\text{ the topological}} \\
{\text{type of }} P {\text{ remains the same as that of }} P_0, {\text{
including also that}} \\ 
{\text{which arc-sides lie only on }} {\text{bd}}\,(\varphi K)
{\text{ (respectively only on }} \\
{\text{bd}}\,(\varphi _0 K)) {\text{ and only on }} 
{\text{bd}}\,(\psi L) {\text{ (respectively only on }} 
{\text{bd}}\,(\psi _0 L)), \\
{\text{the arc-sides of }} P_0 {\text{ and }} 
P {\text{ identified via the small perturbations.}}
\endcases
\tag 3.30
$$

So each arc-side of $(\varphi _0 K) \cap (\psi _0 L)$ has {\it{an opposite
arc-side,
its centrally symmetric image with respect to $O_0$.}} (This is not necessarily
opposite according to the cyclic order of the arc-sides.)
$$
\cases
{\text{Then also the entire hypercycles containing these two opposite arc-}} \\
{\text{sides are centrally symmetric images of each other with respect to 
}} O_0 
\endcases
\tag 3.31
$$
(by the conformal model, choosing $O_0 = 0$, 
and by elementary geometry, or by analytic continuation). Now we are going to
show that
$$
\cases
{\text{the hypercycles containing two opposite arc-sides}} \\ 
{\text{of }} (\varphi _0 K) \cap (\psi _0 L) {\text{ have no common finite 
points.}}
\endcases
\tag 3.32
$$

In fact, if both of these hypercycles are boundary components of $\varphi _0
K$ (or of $\psi _0 L$) then they
are disjoint, i.e., have no common finite points.

Now let one of these hypercycles, $\varphi _0 K_1$, say, be a boundary 
component of $\varphi _0 K$, and the other hypercycle, $\psi _0 L_1$, say, 
be a boundary component of $\psi _0 L$. Let the respective base lines be 
$\varphi _0 K_{01}$ and $\psi _0 L_{01}$. Then also $\varphi _0 K_1$ and 
$\psi _0 L_1$ are symmetric images of each other with respect to $O_0$ (cf. 
\thetag{3.31}), and 
the same holds for their base lines $\varphi _0 K_{01}$ and $\psi _0 L_{01}$
as well. We may suppose provisionally
that $O_0$ is the 

\newpage

centre $0$ of the collinear model circle. 

Then for $O_0 \in \varphi _0 K_{01}$ we have $\varphi _0 K_{01} = 
\psi _0 L_{01}$, and therefore $\varphi _0 K_1$ and $\psi _0 L_1$ have two
common infinite points (those of their common base line), 
but have no common finite point, proving \thetag{3.32} in this case.

Now suppose $0 = O_0 \not\in \varphi _0
K_{01}$. Then, by symmetry, $0 = O_0 
\not\in \psi _0 L_{01}$, and, in the collinear model, these base lines
are two opposite sides of a rectangle inscribed to the model circle
$S^1$. Then  $\varphi _0 K_1$, and by symmetry, also $\psi _0 L_1$, lie
either on the same, or on the other side of their own 
base lines, as the other base line lies. In the second case
$\varphi _0 K_1$ and $\psi _0 L_1$ have no common finite or infinite point,
proving \thetag{3.32} in this case. 
In the first case, by Lemma 1.2 we have that $(\varphi _0 K) \cap (\psi _0 L)$
equals the intersection of two closed convex sets, bounded by $\varphi _0
K_1$ and by $\psi _0 L_1$, respectively. However, now this
intersection is bounded by two finite hypercycle arcs, contradicting
\thetag{3.23}. {\it{This ends the proof of \thetag{3.32}}}.

By \thetag{3.32}, we may apply Lemma 3.3. This yields that 
the hypercycles in \thetag{3.32} either constitute the boundary of a
parallel domain of a line, and then $O_0$ (not fixed at $0$ any more!)
can be any point of this base line, or
these hypercycles have altogether four different infinite points, and then the
closed convex set bounded by them has a unique centre of symmetry $O_0$.

We continue with case distinctions.

{\bf{3.1.}} Let us suppose that, e.g., $K$ is the parallel domain of a straight
line (then $L$ cannot be such, by \thetag{3.21}). 
Then, by the alternance property of the
arc-sides of $P_0$, we have that $P_0$ is an arc-quadrangle, and
that two opposite arc-sides of $P_0$ lie (only) on
${\text{bd}}\,(\varphi K)$, and the other two opposite arc-sides of $P_0$ lie 
(only) on ${\text{bd}}\,(\psi L)$.
 
Then by \thetag{3.31}
the centre of symmetry $O_0$ of $P_0$ is the centre of symmetry of the
union of the entire
hypercycles containing any two opposite arc-sides of $P_0$. Then
by \thetag{3.32} and Lemma 3.3,
on the one hand, $O_0$ lies on
the base line of $\varphi _0 K$, and on the other hand, $O_0$ is a point $O'$ 
uniquely determined by
the hypercycles containing the other two opposite arc-sides 
of $P_0$ (which lie in ${\text{bd}}\,(\psi _0 L)$). By central
symmetry of $P_0$ we have that the base line of $\varphi _0 K$ 
contains the above
uniquely determined point $O'$. Then small generic motions of $\varphi _0 K$
and of $\psi _0 L$ (yielding $\varphi K$ and $\psi L$) preserve the
oppositeness relation for the perturbed arc-sides, by Lemma 3.5, but
destroy this incidence property. Hence,
$(\varphi K) \cap (\psi L)$ will be generically not centrally
symmetric. {\it{This ends the proof for case}} {\bf{3.1}}.

{\bf{3.2.}} Let us suppose that $P_0$ is an arc-quadrangle, such that no two
opposite arc-sides of $P_0$ lie on hypercycles with identical infinite points.
Then, by Lemma 3.3, for both pairs of opposite arc-sides 
of $P_0$ the union of the
hypercycles containing them have altogether four infinite points. Moreover, 
$O_0$ coincides with a point $O'$ uniquely determined  by the union of the
hypercycles containing some two opposite arc-sides of $P_0$ (contained in 
${\text{bd}}\,(\varphi _0 K)$) and also with 
a point $O''$ uniquely determined  by the 

\newpage

union of the
hypercycles containing the other two opposite arc-sides of $P_0$ (contained in 
${\text{bd}}\,(\psi _0 L)$). Once more, a small generic motion of $\varphi _0 K$
and of $\psi _0 L$ (yielding $\varphi K$ and $\psi L$) preserves the
oppositeness relation for the perturbed arc-sides, 
by Lemma 3.5, but destroys this coincidence property. 
Hence, $(\varphi K) \cap (\psi L)$ will be generically not centrally
symmetric. {\it{This ends the proof for case}} {\bf{3.2}}.

{\bf{3.3.}} 
By \thetag{3.23} 
we know that $(\varphi _0 K) \cap (\psi _0 L)$ is not bounded by
two finite 
hypercycle arcs. Moreover, in {\bf{3.1}} and {\bf{3.2}} we have
settled
the case when $(\varphi _0 K) \cap (\psi _0 L)$ was an arc-quadrangle. 
Therefore, in what follows, by \thetag{3.28} we may suppose that 
$$
{\text{the arc-polygon }} P_0 {\text{ has }} n \ge 6 {\text{ arc-sides.}} 
\tag 3.33
$$

Let $s_1'$ and $s_1''$ be two opposite arc-sides of $P_0$ (i.e.,
corresponding to each other by the central symmetry of $P_0$, with respect to
$O_0$). Then, by \thetag{3.32} and 
Lemma 3.3, the hypercycles containing $s_1'$ and $s_1''$
have altogether two, or four infinite points. Accordingly, in the first case,
$O_0$ lies on the common base line of these two hypercycles, and in the second
case, $O_0$ is a point uniquely determined by the union of these hypercycles.

{\bf{3.4.}}
We are going to show that 
$$
\cases
{\text{it is impossible that for any two opposite arc-sides }} s_1' 
{\text{ and }} s_1'' {\text{ of }} \\ 
P_0 {\text{ the first case would hold, i.e., that the hypercycles 
containing }} \\
s_1' {\text{ and }} s_1'' {\text{ would have altogether two infinite points.}}
\endcases
\tag 3.34
$$

$$
{\text{Suppose the contrary, i.e., that always the first case holds.}} 
\tag 3.35
$$
Then the central angles, at the centre $0$ of the collinear model, 
of the base lines of the
hypercycles containing any two opposite arc-sides of $P_0$ sum to $2 \pi $.
We mean the central angle as $\pi $ if
the base line passes through the centre $0$ of the model, and as less than $\pi
$ or greater than $\pi $ according to whether $0$ lies on the side of this
base line, not containing, or containing the respective hypercycle. 
Hence 
$$
{\text{the arithmetic mean of these central angles, for all arc-sides of }} 
P_0, {\text{ is }} \pi . 
\tag 3.36
$$

On the other hand, we are going to show that 
$$
\cases
{\text{the sum of these angles for those arc-sides of }} P_0, {\text{
which}} \\ 
{\text{belong to }} {\text{bd}}\,(\varphi _0 K) {\text{ (or to }} 
{\text{bd}}\,(\psi _0 L)) {\text{ is at most }} 2 \pi .
\endcases
\tag 3.37
$$

\newpage

This will follow if we will have shown that the corresponding open
angular domains are disjoint. Since different such
open angular domains are disjoint, this means that 
$$
\cases
{\text{we have to show only that it is impossible that}} \\
{\text{two different arc-sides of }} P_0 {\text{ would lie on the}} \\
{\text{same boundary component of }} \varphi _0 K {\text{ (or of }} 
\psi _0 L). 
\endcases
\tag 3.38
$$

Let $s_1, \ldots, s_n$ be the 
arc-sides of $P_0$, following each other in the
positive orientation. Then we have $n$ oriented chords $c_i$ 
of the collinear model circle, 
which are the base lines of the boundary components either of 
${\text{bd}}\,(\varphi _0 K)$ or of ${\text{bd}}\,(\psi _0 L)$,
containing the 
arc-sides $s_i$ of $P_0$, respectively, and the orientation is as follows.
The first (last) infinite point of $c_i$ is the first (last) infinite
point of the respective
hypercycle boundary component, according to the positive orientation of
${\text{bd}}\,(\varphi _0 K)$ or of ${\text{bd}}\,(\psi _0 L)$.

We investigate the arc-sides
$s_n,s_1,s_2$. Let, e.g., $s_n \cup s_2 \subset
{\text{bd}}\,(\varphi _0 K)$, and $s_1 \subset {\text{bd}}\,(\psi _0 L)$. 
We measure all angles in ${\Bbb{R}}^2$, containing the collinear model of $H^2$
in the canonical way.
We
may suppose that $c_1$ is of horizontal right direction. Then by Lemma 3.4 we
have that $c_n$ goes from upward to downward, and $c_2$ from downward to
upward. This means that the direction of $c_n$ lies in the angular interval
$(- \pi , 0)$, while the direction of $c_2$ lies in the angular interval
$(0, \pi )$. Therefore the sum of the angles of the positive rotations which 
take the 
direction of $c_n$ to the direction of $c_1$, and the direction of $c_1$ to
the direction of $c_2$, lies in $(0, 2 \pi )$. In other words, 
$$
\cases
{\text{the angle of the positive  rotation, lying in the interval}} \\
(0, 2 \pi ), {\text{ which takes the direction of }} c_n {\text{ to the
direction}} \\
{\text{of }} c_2, {\text{ is the sum of the angles of the positive
rotations,}}\\
{\text{lying in the interval }} (0, 2 \pi ), {\text{ which take the
direction}} \\
{\text{of }} c_n {\text{ to the direction of }} c_1, {\text{ and 
the direction of }} c_1 {\text{ to}} \\
{\text{the direction of }} c_2.
\endcases
\tag 3.39
$$
(Observe that without using Lemma 3.4, the sums of the angles of the positive
rotations, lying in the interval $(0, 2 \pi )$, taking the direction of 
$c_n$ to the direction of $c_1$, and the direction of $c_1$ to the direction 
of $c_2$, could be any angle in $(0, 4 \pi )$. Then the angle of positive 
rotation, lying in the interval $(0, 2 \pi )$, taking the direction of 
$c_n$ to the direction of $c_2$, can be just $2 \pi $ less, than the sum of the
above two angles.)

\newpage

Applying \thetag{3.39} to any three consecutive arc-sides of $P_0$, we obtain
that
$$
\cases
{\text{the total rotation of the directions of the base lines $c_i$ (for}} \\
1 \le i \le n+1, {\text{ where }} c_{n+1} := c_1), {\text{ measured in }} 
{\Bbb{R}}^2, {\text{ containing}} \\
{\text{the collinear model of }} H^2 {\text{ in the canonical way, is equal 
to}} \\ 
{\text{the total rotation of the directions of the base lines }} c_i 
{\text{ for}} 
\\
1 \le i \le n+1, {\text{with }} i {\text{ being even, which is by 
\thetag{3.26} the}} \\
{\text{total rotation of the directions of the base lines }} c_i 
{\text{ for}} \\
1 \le i \le n+1, {\text{ with }} s_i \subset {\text{bd}}\,(\varphi _0 K), 
{\text{ which is the total}} \\
{\text{rotation of }} {\text{bd}}\,(\varphi _0 K), {\text{ which is }} 2  \pi . 
\endcases
\tag 3.40
$$
(Observe that ${\text{bd}}\,(\varphi _0 K)$ may have other boundary
components, not contributing to ${\text{bd}}\,P_0$. However, 
then we may delete them, and this 
does not change the total rotation of the directions of the
above investigated
base lines $c_i$. Also observe that
these base lines $c_i$ can be supposed to be distinct, since we already have
settled the case when $K$ was a parallel domain of a straight line,
cf. {\bf{3.1}}.)
Then \thetag{3.39} and \thetag{3.40} imply that if we suppose, like above, 
that the direction of $c_1$ is horizontal
to the right, then the directions of $c_2, \ldots , c_n$ form a strictly
increasing sequence in $(0, 2 \pi )$. In particular, no two sides $s_i$ of
$P_0$ can lie on the same boundary component either of $\varphi K_0$ or of
$\psi L_0$. {\it{Thus \thetag{3.38} is proved.}}

Hence, by \thetag{3.33}, \thetag{3.37} and \thetag{3.40}, 
$$
\cases
{\text{the arithmetic mean of the central angles from \thetag{3.36},}} \\
{\text{for all arc-sides of }} P_0, {\text{ is at most }} 4 \pi /n \le 
2 \pi /3 .
\endcases
\tag 3.41
$$
Then \thetag{3.36} and \thetag{3.41} lead to a contradiction, {\it{and thus
\thetag{3.34} is proved}}.

{\bf{3.5.}}
By \thetag{3.32}, Lemma 3.3 and \thetag{3.34}, 
$$ 
\cases
{\text{there exist opposite arc-sides }} s_1' {\text{ and }} s_1'' {\text{ of }}
P_0, {\text{ such}} \\
{\text{that the hypercycles containing them have altogether}} \\
{\text{four infinite points, and the infinite points of the}} \\
{\text{hypercycles containing }} s_1' {\text{ and }} s_1'' {\text{ do not 
separate}} \\
{\text{each other on the boundary }} S^1 {\text{ of the model circle.}}
\endcases
\tag 3.42
$$
In this case, by Lemma 3.3, 
the centre of symmetry $O_0$ of $P_0$, being also the centre of
symmetry of the union of the hypercycles containing $s_1'$ and $s_1''$
(cf. \thetag{3.31}), is a point uniquely determined by this union.

{\bf{3.6.}}
Again we make a case distinction. Either $s_1'$ and $s_1''$ from \thetag{3.42}
belongs to the same boundary from among
the boundaries ${\text{bd}}\,(\varphi _0 K)$ and
${\text{bd}}\,(\psi _0 L)$, or they 

\newpage

belong 
to different boundaries. By the alternance
property of the arc-sides 
of $P_0$, if any of these cases holds for some opposite
pair of arc-sides of $P_0$, 
then the same case holds also for all opposite pairs of arc-sides of $P_0$.

{\bf{3.6.1.}}
We begin with the case when  both $s_1'$ and $s_1''$ from \thetag{3.42}
belong, e.g., to ${\text{bd}}\,(\varphi _0 K)$. 

Let $s_2'$ and $s_2''$ be the arc-sides of $P_0$, 
following the arc-sides $s_1'$ and $s_1''$ in the positive sense. 
Then they are also centrally symmetric images of each other, with respect to
the central symmetry with centre $O_0$. Then $O_0$ coincides with $O_1$, 
which is the unique
centre of symmetry of the union of the hypercycles containing $s_1'$ and
$s_1''$. Depending on the fact whether the union of the hypercycles containing 
$s_2'$ and $s_2''$ has altogether two or four infinite points, $O_0$ lies on a
unique line $l_2$, or coincides with a unique point $O_2$ (cf. Lemma 3.3).

Then take some fixed small generic perturbations $\varphi $
and $\psi $ of $\varphi _0$ and $\psi _0$. We may suppose that they
preserve the
oppositeness relation for the arc-sides, by Lemma 3.5.

Then the perturbed arc-sides $s_1'$ and  $s_1''$ are opposite in $P = (\varphi
K) \cap (\psi L)$ as well. Then 
the perturbed point $O_1$ is uniquely determined. The centre of symmetry of
the union of the hypercycles containing the perturbed arc-sides
$s_2'$ and $s_2''$ (these arc-sides being opposite also in $P$, by Lemma 3.5),
either lies on the
unique perturbed line $l_2$, or is the unique perturbed point $O_2$ (cf. Lemma
3.3). Then generically
the perturbed point $O_1$ will not lie on the perturbed line $l_2$, or will
not coincide with the perturbed point $O_2$. Therefore, 
$(\varphi K) \cap (\psi L)$ is generically not centrally symmetric.
{\it{This ends the proof for case}} {\bf{3.6.1}}.

{\bf{3.6.2.}} 
We turn to the other case when, e.g., $s_1'$ belongs to  
${\text{bd}}\,(\varphi _0 K)$, while $s_1''$ belongs to ${\text{bd}}\,(\psi _0
L)$.

We choose $s_2'$ and $s_2''$ like in {\bf{3.6.1}}. Then, by the alternance
property, $s_2'$ belongs to ${\text{bd}}\,(\psi L)$, and $s_2''$ belongs to  
${\text{bd}}\,(\varphi K)$. Then $s_1'$ and $s_1''$, as well as  $s_2'$ and
$s_2''$ are opposite pairs of arc-sides 
for $(\varphi _0 K) \cap (\psi _0 L)$, and
their perturbations are opposite pairs of arc-sides 
for $(\varphi K) \cap (\psi L)$, by Lemma 3.5. 

Then 
$$
\cases
{\text{we take }} \varphi := \varphi _0, {\text{ while }} \psi  {\text{ is 
obtained from }} \psi _0 {\text{ by applying after }} \psi _0 \\
{\text{a small translation along the base line of the hypercycle 
containing}}\\
s_1', {\text{ so that the topological type of }} P {\text{ remains the same 
as that of}} \\ 
P_0, {\text{ including also that which arc-sides lie on }} 
{\text{bd}}\,(\varphi K) {\text{ (respectively}} \\
{\text{on }} {\text{bd}}\,(\varphi _0 K)) {\text{ and on }} 
{\text{bd}}\,(\psi L) {\text{ (respectively on }} {\text{bd}}\,(\psi _0 L)) 
{\text{ --- the}} \\
{\text{sides of }} P_0 {\text{ and of }} P {\text{ identified via the small 
perturbations.}}
\endcases
\tag 3.43
$$ 

\newpage

Then the base line from \thetag{3.43}
is invariant under all (not only small)
such translations,
and in general, the orbits of the points of $H^2$ for all (not only small)
such translations are the hypercycles, i.e., (signed) distance curves, for this 
base line. 

Now consider the base lines $l_1'$ and $l_2'$ 
of the hypercycles containing $s_1'$ and $s_2'$. These have the same
infinite points as the respective hypercycles. Hence, by Lemma 3.4, the 
(altogether four) infinite
points of these base lines separate each other on the boundary $S^1$ of the
model circle (conformal, or collinear). Therefore 
$$
{\text{the base lines }} l_1' {\text{ and }} l_2' {\text{ intersect each 
other.}} 
\tag 3.44
$$

For simplicity, let us assume that $l_1'$ contains $0$ and is horizontal, and
the small translation is to the left hand side. Let $l_1''$ and $l_2''$
denote the base
lines of the hypercycles containing the arc-sides $s_1''$ and $s_2''$. 
Then, by \thetag{3.42}, and using the collinear model, the
straight lines $l_1'$ and $l_1''$ have no common finite or infinite point. We
may suppose that $l_1''$ lies above $l_1'$.
$$
{\text{Let }} d > 0 {\text{ denote the distance of }} l_1' {\text{ and }} l_1''.
\tag 3.45
$$
The translation along $l_1'$ preserves the segment $s$
realizing this distance $d$,
i.e., takes this original segment to the segment realizing
the distance of $l_1'$ and the translated line $l_1''$. Now observe that
$s$ is orthogonal to both $l_1'$ and $l_1''$, and both of $l_1'$ and
the hypercycle at
distance $d$ from $l_1'$ are symmetrical and orthogonal
to the line spanned by $s$.
Hence this hypercycle is orthogonal to $s$ as well.
In other words,
$$
\cases
l_1'' {\text{ moves so that it always touches the (signed)}} \\
{\text{distance line at distance }} d {\text{ for the base line }} l_1'. 
\endcases
\tag 3.46
$$
  
Let the intersecting straight lines $l_1'$ and $l_2'$ (cf. \thetag{3.44}) 
enclose an angle of size $\alpha '$. We mean
the angle of the open angular domain, disjoint to the convex hulls of the
hypercycles containing the arc-sides $s_1'$ and $s_2'$ (``inner angle''). 

Recall that the arc-sides 
$s_1'$ and $s_2''$ belong to ${\text{bd}}\,(\varphi K)$, and the arc-sides 
$s_2'$ and $s_1''$ belong to ${\text{bd}}\,(\psi L)$. Therefore
the hypercycles containing the arc-sides $s_1'$ and $s_2''$, as well as their  
base lines $l_1'$ and $l_2''$ are not moved by our translation, 
but the hypercycles
containing the arc-sides $s_2'$ and $s_1''$, as well as their  
base lines $l_2'$ and $l_1''$ are moved by our translation. 
However, this translation is
a congruence, hence preserves the above described ``inner''
angle $\alpha '$ of the fixed $l_1'$ and the
moving $l_2'$. Now
let us investigate the ``opposite'' angle $\alpha ''$ of $l_1''$ and
$l_2''$, again meant as 
the angle of the open angular domain, disjoint to the convex hulls of the
hypercycles containing the arc-sides $s_1''$ and $s_2''$ (``inner angle'').

Recall \thetag{3.46}. Let $(l_1'')_{\text{new}}$ denote the translated position
of the straight line $l_1''$. 

\newpage

For a sufficiently small translation we have that 
$l_1''$ and $(l_1'')_{\text{new}}$ intersect each other (even their 
directions ``to the left'' are close to each other, in the collinear model,
in the Euclidean sense), 
and both intersect the fixed $l_2''$. Thus 
$$
l_1'',\,\, (l_1'')_{\text{new}} {\text{ and }} l_2'' {\text{ bound a 
triangle }} T , {\text{ and}}
\tag 3.47
$$
$$
\cases
{\text{the ``inner'' angle }} \alpha '' {\text{ gets moved to the analogously 
defined}} \\
{\text{``inner angle'' of }}
(l_1'')_{\text{new}} {\text{ and }} l_2'', {\text{ denoted by }} 
\alpha ''_{\text{new}}.
\endcases
\tag 3.48
$$
Then 
$$
\cases
\alpha '' {\text{ is an inner angle of }} T, {\text{ at its vertex }}
l_1'' \cap l_2'', {\text{ and}} \\
\pi - \alpha ''_{\text{new}} {\text{ is the inner angle of }} T, {\text{ at 
its vertex }} (l_1'')_{\text{new}} \cap l_2''. 
\endcases
\tag 3.49
$$ 
Let the angle of $T$ at its vertex $l_1'' \cap (l_1'')_{\text{new}}$ be 
$\beta ''$. Then
$$
\alpha '' + (\pi - \alpha ''_{\text{new}}) < 
\alpha '' + (\pi - \alpha ''_{\text{new}}) + \beta '' < \pi , {\text{ thus }} 
\alpha '' < \alpha ''_{\text{new}} .
\tag 3.50
$$
However, by an eventual 
central symmetry of $(\varphi K) \cap (\psi L)$, the moved
arc-sides $s_1'$ and $s_2'$
should be taken over to the moved arc-sides 
$s_1''$ and $s_2''$ (cf. Lemma 3.5), 
therefore the angle 
$\alpha ' \,\,( = \alpha '')$  of the base lines of the hypercycles
containing the moved arc-sides $s_1'$ and $s_2'$
should be taken over by this central symmetry 
to the angle $\alpha ''_{\text{new}}$, and then
$$
\alpha '' = \alpha ' = \alpha ''_{\text{new}} .
\tag 3.51
$$
Then \thetag{3.50} and \thetag{3.51} yield a contradiction. {\it{This ends the
proof of case}} {\bf{3.6.2}}, {\it{and thus the proof of \thetag{3.21},
and thus the proof of Theorem}} 3.
$\blacksquare $
\enddemo


\demo{Proof of Theorem 4, {\bf{continuation}}}
Recall that we already have to prove only the implication $(1) \Rightarrow
(2)$ of this Theorem, cf. {\bf{1}} of the proof of this Theorem.

By Lemma 4.3 both conclusions (1) and (2) of Lemma 1.5 hold, and moreover,
$$
\cases
{\text{the constant sectional curvatures in Lemma 1.5 (1) are}} \\
{\text{positive for }} S^d {\text{ and }} {\Bbb{R}}^d, {\text{ and are greater 
than }} 1 {\text{ for }} H^d.
\endcases
\tag 4.7
$$
By Lemma 1.8, (1) of Theorem 4 and conclusions
(1) and (2) of Lemma 1.5 imply the conclusions of Lemma 1.8.
By Lemma 1.9, (1) of Theorem 4 and 
the conclusions of Lemma 1.8 imply the conclusions of Lemma 1.9, namely 
(2) of Theorem 1. In 

\newpage

particular,
(1) of Theorem 4 implies (2) of Theorem 1. 

Then \thetag{4.7} ensures that, under the hypotheses of Theorem 4 and (1) of
Theorem 4,
in (2) of Theorem 1 parasphere or (congruent) hypersphere 
connected components of the boundaries of $K$ and $L$ cannot
occur. (Also recall 
that (2) of Theorem 1 excluded hyperplane connected components
for ${\Bbb{R}}^d$ and $H^d$, which is now a consequence of \thetag{4.7}.)
That is, by (2) of Theorem 1, $K$ and $L$ are congruent balls, and,
by \thetag{4.7}, for the case of $S^d$ they have radius less than $ \pi /2$.

This proves that (1) of Theorem 4, without the statement in brackets, 
implies (2) of
Theorem 4. Since in the beginning of the proof we could assume \thetag{4.2},
and we have
\thetag{4.6}, we have that (1) of Theorem 4, even when only
taken with the statement in brackets, also implies (2) of Theorem 4. 

{\it{This ends the proof of Theorem 4.}}
$ \blacksquare $
\enddemo


\definition{Acknowledgements} 
The authors express their gratitude to I. B\'ar\'any for
carrying the problem and bringing the two authors
together.
We also thank the anonymous referee, whose suggestions have greatly improved
the presentation of the material.
\enddefinition




\Refs

\widestnumber\key{WWW}


\ref 
\key 
\book 
\by  
\publ 
\publaddr 
\yr 
\endref 

\ref 
\key  
\by 
\paper  
\jour 
\pages  
\endref   

\ref
\key 
\by 
\paper 
\jour 
\vol 
\yr 
\pages  
\endref 

\ref 
\key AVS 
\book Geometry of spaces of constant curvature
\by D. V. Alekseevskij, E. B. Vinberg, A. S. Solodovnikov
\publ Geometry II (Ed. E. B. Vinberg), Enc. Math. Sci. {\bf{29}}, 
1-138, Springer 
\publaddr Berlin
\yr 1993
\MR {\bf{95b:}}{\rm{53042}}
\endref 

\ref  
\key Ba
\book Nichteuklidische Geometrie. Hyperbolische Geometrie der Ebene 
\by R. Baldus
\publ 4-te Aufl. Bearb. und erg\"anzt von F. L\"obell, Sammlung G\"oschen
{\bf{970/970a}}, de Gruyter
\publaddr Berlin 
\yr 1964
\MR {\bf{29\#}}{\rm{3936}}
\endref 

\ref 
\key BF 
\book Theorie der konvexen K\"orper, {\rm{Berichtigter Reprint}}
\by T. Bonnesen, W. Fenchel
\publ Springer
\publaddr Berlin-New York
\yr 1974.
{\rm{English transl.:}} 
{\it{Theory of convex bodies}}, 
{\rm{Translated from the 
German and edited by L. Boron, C. Christenson, B. Schmidt, BCS
Associates, Moscow, ID, 1987}} 
\MR {\bf{49\#}}{\rm{9736}}, {\bf{88j:}}{\rm{52001}}
\endref 

\ref 
\key Bo
\book Non-Euclidean geometry, a critical and
historical study of its developments, {\rm{With a Supplement containing the
G. B. Halstead translations of}} ``The science of absolute space'' {\rm{by 
J. Bolyai and}} ``The theory of parallels'' {\rm{by N. Lobachevski}}
\by R. Bonola
\publ Dover Publs. Inc.
\publaddr New York, N.Y.
\yr 1955
\MR {\bf{16-}}{\rm{1145}}
\endref 

\ref 
\key C
\book Non-Euclidean Geometry, 6-th ed. 
\by H. S. M. Coxeter 
\publ Spectrum Series, The Math. Ass. of America
\publaddr Washington, DC
\yr 1998
\MR {\bf{99c:}}{\rm{51002}}
\endref 

\ref 
\key HM
\by E. Heil, H. Martini
\paper Special convex bodies 
\jour 
In: Handbook of Convex
Geometry, (eds. P. M. Gruber, J. M. Wills), North-Holland,
Amsterdam etc., 1993, Ch. 1.11
\pages 347-385 
\MR {\bf{94h:}}{\rm{52001}}
\endref   
 
\ref
\key H
\by R. High
\paper Characterization of a disc, Solution to problem 1360 (posed by
P. R. Scott)
\jour Math. Magazine
\vol 64 
\yr 1991
\pages 353-354 
\endref 

\newpage

\ref 
\key L
\book Nichteuklidische Geometrie, {\rm{3-te Auflage}}
\by H. Liebmann
\publ de Gruyter
\publaddr Berlin
\yr 1923.
Jahr\-buch Fortschr. Math. {\bf{49}}, 390
\endref 

\ref 
\key P
\book Nichteuklidische Elementargeometrie der Ebene, {\rm{Math. Leitf\"aden}}
\by O. Perron 
\publ Teubner
\publaddr Stuttgart
\yr 1962
\MR {\bf{25\#}}{\rm{2489}}
\endref 

\ref 
\key Sch
\book Convex bodies: the Brunn-Minkowski theory, {\rm{Encyclopedia of Math. 
and its Appl.,}} {\bf{44}}
\by R. Schneider 
\publ Cambridge Univ. Press
\publaddr Cambridge
\yr 1993.
{\rm{Second expanded ed., Encyclopaedia of Math. and
its Appl.,}} {\bf{151}}, {\rm{Cambridge Univ. Press, Cambridge, 2014}} 
\MR {\bf{94d:}}{\rm{52007}}, {\bf{3155183}} 
\endref 

\ref
\key So05
\by V. Soltan
\paper Pairs of convex bodies with centrally symmetric intersections of
translates
\jour Discr. Comput. Geom.
\vol 33
\yr 2005
\pages 605-616
\MR {\bf{2005k:}}{\rm{52012}} 
\endref

\ref
\key So06 
\by V. Soltan
\paper Line-free convex bodies with centrally symmetric intersections of
translates
\jour Revue Roumaine Math. Pures Appl.
\vol 51
\yr 2006
\pages 111-123
\MR {\bf{2007k:}}{\rm{52010}}
{\rm{Also in: Papers on Convexity and Discrete Geometry, Ded. to T. Zamfirescu
on the occasion of his 60-th birthday, Editura Acad. Rom\^ane, Bucure\c sti,
2006, 411-423}}
\endref

\ref 
\key St
\book Differential Geometry
\by J. J. Stoker 
\publ New York Univ., Inst. Math. Sci.
\publaddr New York
\yr 1956
\endref 

\ref
\key V
\by I. Vermes
\paper \"Uber die synthetische Behandlung der
Kr\"ummung und des Schmieg\-zy\-kels der ebenen Kurven in der
Bolyai-Lobatschefskyschen Geometrie
\jour Stud. Sci. Math. Hungar.
\vol 28
\yr 1993
\pages 289-297 
\MR {\bf{95e:}}{\rm{51030}} 
\endref 

\endRefs
\enddocument